\documentclass[a4paper,11pt]{amsart}
\usepackage[utf8]{inputenc}
\usepackage[left=2.7cm,right=2.7cm,top=3.5cm,bottom=3cm]{geometry}
\usepackage{amsthm,amssymb,amsmath,amsfonts,mathrsfs,amscd}
\usepackage{bm}
\usepackage[dvipsnames]{xcolor}
\usepackage{tikz-cd}
\usepackage[pagebackref]{hyperref}
\usepackage{latexsym}
\usepackage[english]{babel}
\usepackage{mathtools,braket}
\usepackage[all]{xy}
\usepackage{booktabs}
\usepackage{comment}
\usepackage{verbatim}
\usepackage{enumitem}
\theoremstyle{plain}
\newtheorem{theorem}{Theorem}[section]
\newtheorem*{thmintro}{Theorem}
\newtheorem{lemma}[theorem]{Lemma}
\newtheorem{maintheorem}{Theorem}
\newtheorem{maincor}[maintheorem]{Corollary}

\newtheorem{prop}[theorem]{Proposition}
\newtheorem{assumption}{Assumption}[section]
\newtheorem{cor}[theorem]{Corollary}
\newtheorem{conj}[theorem]{Conjecture}

\theoremstyle{definition}
\newtheorem{defn}[theorem]{Definition}

\theoremstyle{remark}
\newtheorem{rk}{Remark}

\newenvironment{smatrix}{\left(\begin{smallmatrix}}{\end{smallmatrix}\right)}

\newcommand{\hf}{{\mathbf{f}}}
\newcommand{\hg}{{\mathbf{g}}}
\newcommand{\hh}{{\mathbf h}}

\newcommand{\hr}{{\mathbf r}}
\newcommand{\htt}{{\mathbf t}}

\newcommand{\A}{\mathbb{A}}

\DeclareMathOperator{\R}{\mathbb{R}}
\DeclareMathOperator{\C}{\mathbb{C}}

\DeclareMathOperator{\Q}{\mathbb{Q}}
\DeclareMathOperator{\dep}{dep}
\DeclareMathOperator{\diag}{diag}

\DeclareMathOperator{\FL}{FL}

\DeclareMathOperator{\Frac}{Frac}

\DeclareMathOperator{\GL}{GL}

\DeclareMathOperator{\GSp}{GSp}
\DeclareMathOperator{\Hom}{Hom}

\DeclareMathOperator{\imp}{impr}
\DeclareMathOperator{\Ind}{Ind}
\DeclareMathOperator{\Iw}{Iw}
\DeclareMathOperator{\Kl}{Kl}

\DeclareMathOperator{\ord}{ord}

\DeclareMathOperator{\Si}{Si}
\DeclareMathOperator{\sign}{sign}

\DeclareMathOperator{\Spf}{Spf}

\DeclareMathOperator{\Sym}{Sym}

\providecommand{\ord}{\mathrm{ord}}
\providecommand{\Tr}{\operatorname{Tr}}
\hypersetup{hidelinks,pdftitle={Exceptional zeros for GSp(4) x GL(2) and Yoshida lifts},pdfauthor={Armando Gutierrez Terradillos, Lois Omil-Pazos, Oscar Rivero}}
\begin{document}

\title[Exceptional zeros for $\GSp(4) \times \GL(2)$]{Exceptional zeros for $\GSp(4) \times \GL(2)$ and Yoshida lifts}

\author{Armando Gutierrez Terradillos, Lois Omil-Pazos, and \'Oscar Rivero}

\date{\today}

\address{A. G. T.: Department of Mathematics, Aarhus University, Ny Munkegade 118, DK-8000, Aarhus, Denmark}
\email{armangute@math.au.dk}

\address{L. O.-P.: CITMAga, Santiago de Compostela, Spain}
\email{lois.omil.pazos@usc.es}

\address{O. R.: CITMAga and Departamento de Matem\'aticas, Universidade de Santiago de Compostela, Santiago de Compostela, Spain}
\email{oscar.rivero@usc.es}

\begin{abstract}
The aim of this article is to study an exceptional zero phenomenon for the $p$-adic $L$-functions attached to $\GSp_4\times\GL_2$ in the critical regions (F) and (D). We construct improved $p$-adic $L$-functions by replacing the $p$-depleted Eisenstein family with its $p$-ordinary counterpart. Under a local compatibility hypothesis at the ramified places and classical control for the families, we show that their poles detect exactly the Yoshida lifts whose smaller constituent is dual to the $\GL_2$ representation. The proof compares their residues with distinguished automorphic periods, extending the coherent period formula to the case where one modular form has weight zero. An explicit Euler-factor identity relates each improved function to the usual $p$-adic $L$-function. At these Yoshida points, the vanishing factor cancels the pole, and the resulting value is expressed through the distinguished period and the logarithmic derivative of the ordinary elliptic root. We also discuss the expected factorization along Yoshida families into Hida--Rankin $p$-adic $L$-functions.
\end{abstract}

\maketitle
\enlargethispage{4pt}

\setcounter{tocdepth}{1}
\tableofcontents
\section{Introduction}

The poles of the complex $L$-function for $\GSp_4\times\GL_2$ detect a natural relation between its two inputs: the Siegel representation is a Yoshida lift, and the $\GL_2$ (elliptic) representation is dual to one of its constituents. This relation is also characterized by the nonvanishing of a distinguished period. We seek a $p$-adic counterpart of this criterion in the critical regions (F) and (D). At the relevant weights, an interpolation factor may vanish and cancel a pole. Thus, removing this factor leads to improved $p$-adic $L$-functions whose residues recover the distinguished period.

The two regions lead to different coherent realizations: in (F), a degree-two class on the Siegel threefold pairs with two degree-zero elliptic classes, while in (D) the Siegel class and the second elliptic class both have degree one. In either construction, the first $\GL_2$ factor contributes an Eisenstein series. Before $p$-depletion, its $p$-adic family has a simple pole with constant residue. Taking this residue should then recover the distinguished period. The difficulty is that the resulting coherent pairing has first elliptic weight zero, beyond the range of the positive-weight period formula used in the interpolation theorem.

We establish the period formula in both coherent degrees. Using the automorphic realization of coherent cohomology, we identify the cocycles and coefficient contractions that enter the pairing, and recover the global integral through the Serre trace. The complex and $p$-adic constructions can then be compared as realizations of the same rational pairing. Together with the local test-vector calculations, this gives a common Yoshida criterion for the poles of the two improved functions. Moreover, an explicit Euler-factor identity relates each of them to its ordinary counterpart. We call a Yoshida point {\it matching} if the $\GL_2$ representation is contragredient to the smaller Yoshida constituent, with the corresponding ordinary refinement. At a matching Yoshida point, the ordinary value is governed by the coherent period and the logarithmic derivative of the elliptic ordinary root.

\subsection{The complex pole and the distinguished period}
Let $\Pi$ be a unitary globally generic cuspidal automorphic representation of $\GSp_4(\mathbb A)$, of cohomological weight $(k_1,k_2)$ with $k_1\geq k_2\geq3$, and let $\Sigma$ be a unitary cuspidal representation of $\GL_2(\mathbb A)$. Put
\[
H=\GL_2\times_{\GL_1}\GL_2,
\]
with the equal-determinant embedding $\iota:H\hookrightarrow\GSp_4$. When $\omega_\Pi\omega_\Sigma=1$, the global period is
\[
\mathcal I(\Psi,F)=\int_{Z_G(\mathbb A)H(\mathbb Q)\backslash H(\mathbb A)}\Psi(\iota(h))F(h_2)\,dh,
\]
where $Z_G$ is the common scalar centre. We call the pair distinguished if this integral is nonzero for some $\Psi\in\Pi$ and $F\in\Sigma$.

Write $\Lambda(s,\Pi\times\Sigma)$ for the completed degree-eight Langlands $L$-function in unitary normalization. Its relation to Yoshida lifts follows from the generic spin transfer of Asgari--Shahidi and the Rankin--Selberg pole theorem. Novodvorsky's integral connects the pole with $\mathcal I$: the residue of its Eisenstein series at $s=1$ is constant. For this integral comparison we require equality of the local Novodvorsky and Langlands factors at ramified finite places, as in the assumption \ref{eq:local-compatibility}. This holds, in particular, when each corresponding $\Sigma_v$ is non-supercuspidal \cite{Loeffler:Nov}. We recall the resulting criterion in Section~\ref{SectionSpinStandardRep}.

\begin{thmintro}\label{intro:complex}
Let $\Pi$ and $\Sigma$ be as above, with $\omega_\Pi\omega_\Sigma=1$. Under the local compatibility expressed in \ref{eq:local-compatibility}, the following conditions are equivalent:
\begin{enumerate}[label=\textup{(\roman*)}]
\item $\Lambda(s,\Pi\times\Sigma)$ has a pole at $s=1$;
\item $\Pi$ is a Yoshida lift of $\Sigma^\vee\boxtimes\Sigma'$ for a cuspidal representation $\Sigma'\not\simeq\Sigma^\vee$ with $\omega_{\Sigma'}=\omega_\Sigma^{-1}$;
\item the pair $(\Pi,\Sigma)$ is distinguished.
\end{enumerate}
When these conditions hold, the pole is simple.
\end{thmintro}

This criterion supplies the complex input to our argument. Its proof combines \cite{AsgariShahidi,GT11,JPSS83} with the residue and local-integral results recalled in Section~3. For the $p$-adic application, we evaluate this period on automorphic forms with prescribed compact types at infinity and ordinary vectors at $p$. When the pair is distinguished, the test-vector calculations below show that these choices allow a nonzero period for suitable data at the remaining finite places.

\subsection{The coherent period at weight zero}
The boundary relevant to both regions occurs at elliptic weight \[ c_0=k_1-k_2+2, \] the weight of the smaller Yoshida constituent. Taking the Eisenstein residue replaces the form on the first $\GL_2$ factor by $\mathbf1$, leaving a pairing between the Siegel class and an elliptic cusp class of weight $c_0$.

The positive-weight formula of \cite[Theorem 5.6(b)]{LPSZ1} is stated with first elliptic weight $1+t_1$, where $t_1\geq0$, while our specialization has $t_1=-1$. To reach it, we work with the Dolbeault representatives in the Harris--Su realization \cite{Harris:Partial,Su24} and check that restriction to the $H$-surface, exterior product and algebraic coefficient contraction commute with the automorphic lift. At weight zero, the contractions select Siegel vectors of compact weights $(0,-c_0)$ and $(0,c_0)$, paired respectively with holomorphic and antiholomorphic elliptic representatives. Evaluating the corresponding cochains gives nonzero constants. Then, cuspidal decay allows us to identify the integral with the trace on the toroidal compactification, with the constant form in the first elliptic factor.

The resulting formula is stated in Theorem~\ref{Integral:New:Pairing}. Writing $A_\diamond$ for the proportionality scalar in the weight-zero formula of that theorem, we obtain
\[
P_F(\eta,\mathbf1\boxtimes g)=A_F\mathcal I_F\,e_m,\qquad
P_D(\xi,\mathbf1\boxtimes\upsilon)=A_D\mathcal I_D\,e_m,
\qquad m=k_1+k_2-6.
\]
Here $P_F,P_D$ are rational coherent pairings, and $\mathcal I_F,\mathcal I_D$ are the periods of the corresponding automorphic representatives. The pairings take values in the same one-dimensional character space $E\{m\}$; its basis $e_m$ records the twist relating the arithmetic coefficients to unitary automorphic forms. The nonzero constants $A_F,A_D$ account for the cocycle, contraction, measure and trace normalizations.

On the $p$-adic side, the analytic branching kernels of \cite{LZ-BK,LR3} remain polynomials at weight zero. Lemma~\ref{lem:boundary-branching} therefore identifies their specializations with the algebraic contractions. This allows us to compare the complex period and the $p$-adic residue through the same rational pairing.

\subsection{Improved functions and their poles}
Fix ordinary families $\underline\Pi,\underline\Sigma$ satisfying the classical-control hypotheses of Section~\ref{sec:review}. Their classical weights are denoted $(k_1,k_2,c)$, with the underlying eigenfamily points understood. If $s$ is the unitary variable and $s'$ the motivic variable, then
\[
s'=s+\frac{k_1+k_2+c-4}{2}.
\]
Our constructions lie on the three-variable domains
\[
\mathcal W_F^{\imp}:s'=k_1-1,\qquad
\mathcal W_D^{\imp}:s'=c+k_2-3.
\]
Their Eisenstein weights are
\[
t_F=k_1-k_2-c+2,\qquad t_D=-k_1+k_2+c-2.
\]
At the boundary of interest both Eisenstein characters are trivial and $s=0$ on the two domains. This is the point paired with $s=1$ by the functional equation, which also replaces the representations by their contragredients. For a distinguished pair, the completed function has a pole at both points. In Novodvorsky's integral, the Eisenstein residue at $s=0$ is again constant: with the notation of Proposition~\ref{Period:Equals:Residue}, it is $-c_E\varphi(0)$, instead of $c_E\widehat\varphi(0)$ at $s=1$ \cite[Proposition 2.5]{Moriyama:ZetaIntegrals}. Our comparison passes through this common distinguished period.

We construct $L_p^{(F),\imp}$ and $L_p^{(D),\imp}$ by using the ordinary Eisenstein family before depletion in the two coherent pairings. Its pole supplies the possible pole of each improved function. To compute the residue, we first regularize the Eisenstein family in the analytic coefficient module where the pairing is defined. We can then specialize inside the pairing, replacing the Eisenstein series by a multiple of the constant form. The generators and auxiliary test data are fixed on the family neighbourhood, as explained in Section~\ref{sec:improved}.

Fix a good classical point $P$ of the Siegel family and a crystalline point $Q_0$ of weight $c_0$ of the $\GL_2$ family, and put $\Pi=\underline\Pi(P)$ and $\Sigma=\underline\Sigma(Q_0)$. We assume the unramified ordinary refinements, classical control and \'etale elliptic weight map specified in Section~\ref{subsec:poles-comparison}, with both Eisenstein characters trivial at $Q_0$. We also retain the central-character condition and local compatibility of Theorem~\ref{intro:complex}.

\begin{maintheorem}\label{intro:padic}
For suitable tame test data, the following conditions are equivalent:
\begin{enumerate}[label=\textup{(\roman*)}]
\item $L_p^{(F),\imp}(\Pi\times\underline\Sigma)$ has a pole at $Q_0$;
\item $L_p^{(D),\imp}(\Pi\times\underline\Sigma)$ has a pole at $Q_0$;
\item $\Pi$ is a Yoshida lift with $\Sigma^\vee$ as its constituent of weight $c_0$;
\item the pair $(\Pi,\Sigma)$ is distinguished.
\end{enumerate}
When these conditions hold, both poles are simple. For any choice of tame data, a pole of either improved function implies \textup{(iii)}.
\end{maintheorem}

The equivalence comes from comparing the two residues with the distinguished period. Let $z$ be the logarithmic elliptic weight centred at $Q_0$, so that $\ell(\mathbf t_F)=-z$ and $\ell(\mathbf t_D)=z$. The Kubota--Leopoldt normalization gives $\ell(\kappa)\mathcal E^\flat(\kappa)|_{\kappa=0}=-(1-p^{-1})\mathbf1$, and hence
\[
\operatorname{Res}_{z=0}L_p^{(F),\imp}=(1-p^{-1})B_F(\mathbf1),\qquad
\operatorname{Res}_{z=0}L_p^{(D),\imp}=-(1-p^{-1})B_D(\mathbf1),
\]
where $B_\diamond(\mathbf1)$ is the constant-input pairing defined in \eqref{eq:constant-pairing}. Comparing rational structures as in \cite[Theorem 3.10]{LZ25-poles}, we obtain
\[
\frac{\operatorname{Res}_{z=0}L_p^{(\diamond),\imp}}{\Omega_{\diamond,p}}
=\varepsilon_\diamond(1-p^{-1})\,
\iota_p\!\left(\frac{\nu_K A_\diamond\mathcal I_\diamond}{\Omega_{\diamond,\infty}}\right),
\qquad\varepsilon_F=1,\quad\varepsilon_D=-1.
\]
Here $\nu_K$ is the tame volume multiplier in the definition of the global functions; the constants $A_\diamond$ still refer to the geometric pairings $P_\diamond$. The expression inside $\iota_p$ is algebraic by Proposition~\ref{cor:periodcomparison}: observe that after dividing by their comparison periods, the complex integral and the $p$-adic residue come from the same rational pairing. In particular, the residue is nonzero precisely when the period is nonzero on the chosen vectors.

For the converse, the test-vector calculations show that the distinguished period can be nonzero with our prescribed vectors at $p$ and infinity. Choosing suitable tame data then gives a nonzero residue in both regions. The residue also varies analytically along the Eisenstein divisor; near a matching Yoshida point it identifies the smaller constituent at nearby good classical points, as explained after Theorem~\ref{thm:criterion}.

\subsection{The Euler-factor identities}
Let $\alpha,\beta,\gamma,\delta$ be the motivic Siegel Hecke parameters, ordered by the ordinary refinement, and let $\mathfrak a$ be the ordinary elliptic root. The normalized parameters
\[
\beta^\circ=\beta/p^{k_2-2},\qquad
\gamma^\circ=\gamma/p^{k_1-1},\qquad
\mathfrak a^\vee=\chi_\Sigma(p)^{-1}\mathfrak a
\]
are analytic units on the chosen families, with classical values as displayed. On the trivial Eisenstein-character components above, choose the same generators, measures, Gauss sums and Fourier-related tame data for each ordinary--improved pair. Proposition~\ref{teo:interpolation} then gives the following identities.

\begin{maintheorem}\label{intro:euler-comparison}
With these choices, the ordinary and improved functions satisfy the identities
\begin{equation}\label{intro:improvement}
\begin{split}
L_p^{(F)}&=(1-\gamma^\circ\mathfrak a)L_p^{(F),\imp},\\
L_p^{(D)}&=(1-\beta^\circ/\mathfrak a^\vee)L_p^{(D),\imp}.
\end{split}
\end{equation}
as meromorphic functions on their respective three-variable domains.
\end{maintheorem}

Here the ordinary functions are those of \cite{LZ-BK,LR3}, with the common normalizations of Definition~\ref{def:joint-functions} and Remark~\ref{prop:LR-comparison}. We prove the identities by comparing the local integrals at classical points and then using uniqueness on the three-variable family. The same comparison applies to the ordinary Eisenstein input after regularizing its pole.

At a matching Yoshida point the removed Euler factors vanish. Multiplying their first-order terms by the residues gives the following value outside the classical interpolation range.

\begin{maincor}\label{intro:ordinary-value}
At a matching ordinary Yoshida point $(P_0,Q_0)$ satisfying the hypotheses above, fix $P_0$ and take $z$ to be the logarithmic elliptic weight centred at $Q_0$. For $\diamond\in\{F,D\}$,
\[
L_p^{(\diamond)}(P_0,Q_0)
=-(1-p^{-1})\,\partial_z\log\mathfrak a(Q_0)\,B_\diamond(\mathbf1).
\]
For the tame data of Theorem~\ref{intro:padic}, this value is nonzero if and only if $\partial_z\log\mathfrak a(Q_0)\ne0$.
\end{maincor}

Thus the vanishing Euler factor can cancel the pole and leave a nonzero ordinary value. Corollary~\ref{cor:ordinary-value} proves the formula by differentiating \eqref{intro:improvement} in the elliptic weight, with the Siegel point fixed.

The Yoshida transfer gives a further, endoscopic factorization. For constituents $f,g$ of weights $k_1+k_2-2$ and $k_1-k_2+2$, respectively, its finite motivic form is
\[
L_{\mathrm{tr}}^{\mathrm{mot}}(\Pi\times h,s')=L(f\otimes h,s')L(g\otimes h,s'-(k_2-2)).
\]
Section~\ref{sec:factorization} formulates the corresponding conjectural $p$-adic identity and describes the period and functional-equation comparisons needed to obtain it. Such a factorization would express our residue in terms of a Rankin residue and the factor associated with the larger constituent, complementing the direct proof of the pole criterion through coherent periods. We also recall the complex Saito--Kurokawa factorization, which lies outside our globally generic hypotheses.

The diagram below, following \cite{LR3}, displays the three critical regions for fixed Siegel weights. The two domains meet at the lower Yoshida boundary, of elliptic weight $k_1-k_2+2$. The other constituent, of weight $k_1+k_2-2$, lies on the boundary between (D) and (A) and is not detected by these two weight-zero Eisenstein specializations.

\begin{center}
\tikzset{every picture/.style={line width=0.75pt}} %set default line width to 0.75pt

\begin{tikzpicture}[x=0.75pt,y=0.75pt,yscale=-1,xscale=1]
%uncomment if require: \path (0,347); %set diagram left start at 0, and has height of 347

%Straight Lines [id:da46752639414311403]
\draw [color={rgb, 255:red, 165; green, 165; blue, 165 }  ,draw opacity=1 ][line width=3.75]  [dash pattern={on 9pt off 4.5pt}]  (460,120) -- (270,310) ;
%Straight Lines [id:da7207496603743515]
\draw [color={rgb, 255:red, 165; green, 165; blue, 165 }  ,draw opacity=1 ][line width=3.75]    (230,70) -- (470,310) ;
%Straight Lines [id:da7919728211189669]
\draw    (210,310) -- (488,310) ;
\draw [shift={(490,310)}, rotate = 180] [color={rgb, 255:red, 0; green, 0; blue, 0 }  ][line width=0.75]    (10.93,-3.29) .. controls (6.95,-1.4) and (3.31,-0.3) .. (0,0) .. controls (3.31,0.3) and (6.95,1.4) .. (10.93,3.29)   ;
%Shape: Triangle [id:dp3275490883499559]
\draw  [fill={rgb, 255:red, 241; green, 241; blue, 241 }  ,fill opacity=1 ] (380,220) -- (470,310) -- (290,310) -- cycle ;
%Shape: Diamond [id:dp961347636667798]
\draw  [fill={rgb, 255:red, 241; green, 241; blue, 241 }  ,fill opacity=1 ] (380,140) -- (410,170) -- (380,200) -- (350,170) -- cycle ;
%Straight Lines [id:da2906033262064447]
\draw [fill={rgb, 255:red, 241; green, 241; blue, 241 }  ,fill opacity=1 ]   (460,40) -- (380,120) -- (300,40) ;
%Straight Lines [id:da6593245785109718]
\draw    (210,310) -- (210,42) ;
\draw [shift={(210,40)}, rotate = 90] [color={rgb, 255:red, 0; green, 0; blue, 0 }  ][line width=0.75]    (10.93,-3.29) .. controls (6.95,-1.4) and (3.31,-0.3) .. (0,0) .. controls (3.31,0.3) and (6.95,1.4) .. (10.93,3.29)   ;
%Straight Lines [id:da07020087697054223]
\draw  [dash pattern={on 4.5pt off 4.5pt}]  (380,30) -- (380,310) ;
%Shape: Polygon [id:ds8457648750180513]
\draw  [fill={rgb, 255:red, 241; green, 241; blue, 241 }  ,fill opacity=1 ] (340,180) -- (370,210) -- (280,300) -- (250,270) -- cycle ;
%Straight Lines [id:da8926622329233427]
\draw [fill={rgb, 255:red, 241; green, 241; blue, 241 }  ,fill opacity=1 ]   (280,40) -- (370,130) -- (340,160) -- (220,40) ;

% Text Node
\draw (201,22.4) node [anchor=north west][inner sep=0.75pt]    {$c$};
% Text Node
\draw (481,314.4) node [anchor=north west][inner sep=0.75pt]    {$s$};
% Text Node
\draw (381,52.4) node [anchor=north west][inner sep=0.75pt]    {$( A)$};
% Text Node
\draw (377,162.4) node [anchor=north west][inner sep=0.75pt]    {$( D)$};
% Text Node
\draw (381,282.4) node [anchor=north west][inner sep=0.75pt]    {$( F)$};
% Text Node
\draw (380,315.4) node [anchor=north] [inner sep=0.75pt]    {$\tfrac{1}{2}$};
% Text Node
\draw (292,233.4) node [anchor=north west][inner sep=0.75pt]    {$( E)$};
% Text Node
\draw (281,72.4) node [anchor=north west][inner sep=0.75pt]    {$( B)$};

\fill (370,210) circle (2pt);
\draw (415,207) node [anchor=west] {$Q_0$};
\draw (372,210) -- (414,207);
\end{tikzpicture}
\end{center}
{\small The vertical coordinate is the elliptic weight $c$ and the horizontal coordinate is the unitary variable $s$; the vertical dashed line is $s=1/2$. The solid grey line is $\mathcal W_F^{\imp}$ and the dashed grey line is $\mathcal W_D^{\imp}$; they meet at $Q_0=(0,k_1-k_2+2)$. Regions (B) and (E) are the noncritical regions separating the three critical regions (A), (D) and (F). The left vertical axis is schematic and is not the line $s=0$.}

\subsection{Connection with other works}

The cohomological constructions used here rely on higher Hida and Coleman theory, developed by Pilloni and Boxer--Pilloni \cite{Pil20,boxerpilloni20}. The four-variable work of Graham--Rockwood \cite{GrahamRockwood} and the automorphic constructions of Liu \cite{liu1,liu2} give complementary approaches. We work throughout with the three-variable cohomological constructions. Liu's Yoshida factorization provides a useful point of comparison for Section~\ref{sec:factorization}. In region (D) we use the construction of \cite{LR3}; since its support conditions require a deep auxiliary level, we also make explicit the normalization relating it to the algebraic pairing at level $p^2$ in Proposition~\ref{prop:running-level}.

The closest model for the pole argument is the improved Asai function of Loeffler--Zerbes \cite{LZ25-poles}, where a residue detects a distinguished representation. In the complex setting, residues for $\GSp_4\times\GL_2$ were also studied by F. Lemma \cite{lemma20}. Related questions for $\GSp_6$ appear in Cauchi--Lemma--Rodrigues Jacinto \cite{CLJ}.

\subsection{Organization of the paper}
The argument has three steps: the weight-zero period formula, the specialization of the Eisenstein residue, and the local calculation of the Euler factor. Section~\ref{sec:GL2xGL2} recalls the $\GL_2\times\GL_2$ model, which was the main motivation behind this note. Sections~3 and~4 establish the complex pole criterion and its coherent interpretation. Section~\ref{sec:review} constructs the family pairings. Section~\ref{sec:improved} defines the improved functions, compares their classical specializations with the coherent pairings, and proves the pole criterion and value formula. Section~\ref{sec:factorization} discusses the expected endoscopic factorization.

\subsection{Acknowledgements}
The authors would like to thank Raúl Alonso and David Loeffler for informative conversations and helpful exchanges related to the subject of this article. A.G.T. thanks Antonio Cauchi for several helpful discussions related to this project. The third-named author is grateful to Paul Nelson for funding his research visit to Aarhus University in May 2025, during which several of the ideas in this paper were developed.

During the preparation of this paper, the authors were supported by PID2023-148398NA-I00, funded by MICIU/AEI/10.13039/501100011033 and by ERDF, EU.

\subsection{AI disclosure}
We used the AI assistants ChatGPT and Claude (specifically, the models Astra~6 and Opus~5.5) to help develop and check technical arguments, draft and revise proofs, and improve the exposition. The most substantial mathematical use concerned the weight-zero coherent period formula in Theorem~\ref{Integral:New:Pairing}, the supported pairing and level comparison in Propositions~\ref{prop:supported-pairing} and~\ref{prop:running-level}, and the local calculations and test-vector arguments leading to Proposition~\ref{prop:ptest}. We also used these tools extensively to rewrite parts of the manuscript and make the exposition clearer, particularly in Sections~\ref{subsec:family-duality} and~\ref{subsec:classical-comparisons}, by explaining the motivation for the constructions, reorganizing proofs and reducing repetition.

%%%%%%%%%%%%%%%%%%%%%%%%%%%%%%%%%%%%%%%%%%%%%%%%%%%%%%
% The case of GL2xGL2
%%%%%%%%%%%%%%%%%%%%%%%%%%%%%%%%%%%%%%%%%%%%%%%%%%%%%%

\section{\texorpdfstring{The $\GL_2\times\GL_2$ case}{The GL2 x GL2 case}}\label{sec:GL2xGL2}

We begin with the $\GL_2\times\GL_2$ case, following Hida and Rivero--Rotger \cite{hida2,RR21}. It already displays the relation between poles, distinguished periods and vanishing Euler factors that we will study for $\GSp_4\times\GL_2$. In the complex setting, the Rankin--Selberg $L$-function has a pole precisely when the two cuspidal representations are contragredient. The corresponding $p$-adic statement concerns an improved function of two variables: after fixing the first eigenform, a pole detects a contragredient specialization of the second family. We recall both statements, keeping track of the normalizations and of the Euler factor relating the improved function to the usual three-variable Hida--Rankin function.
\subsection{Notation}
\begin{itemize}
    \item We write $\overset{.}{=}$ to denote equality up to a nonzero constant.
    \item Let $G_1$ and $G_2$ be groups, let $V$ and $W$ be representations of $G_1$, and let $M$ be a representation of $G_2$. We write $V\otimes W$ for the internal tensor product of $V$ and $W$, and $V\boxtimes M$ for the external tensor product, viewed as a representation of $G_1\times G_2$.
    \item Let $E$ be a global field. For each place $v$ of $E$, we denote by $E_v$ the corresponding completion. If $v$ is finite, we fix a uniformizer $\varpi_v$ and denote by $q_v$ the cardinality of the residue field.
    \item The adelic absolute value is normalized as $|\cdot|_{\A}=|\cdot||\cdot|_f$, with the local absolute values satisfying $|\varpi_v|_v=q_v^{-1}$ at every finite place $v$.
\end{itemize}
\subsection{The complex case and normalizations}

Let $f$ and $g$ be normalized cuspidal newforms of weights $k,\ell\geq2$ and nebentypes $\chi_f,\chi_g$. We write $\Sigma_f,\Sigma_g$ for their unitary automorphic representations. Throughout this section, $s$ denotes the unitary automorphic variable and $j$ the motivic variable, related by
\[
L(\Sigma_f\times\Sigma_g,s)=L(f\otimes g,j),
\qquad j=s+\frac{k+\ell-2}{2}.
\]
Here both sides denote primitive finite $L$-functions. Write $f^*$ for the normalized newform associated with $\Sigma_f^\vee$. Thus $f^*$ is the primitive form associated with the twist $f\otimes\chi_f^{-1}$; under a complex embedding, its Fourier coefficients are the conjugates of those of $f$. In particular, $f^*$ has weight $k$.

\begin{theorem}\label{thm:GL2xGL2}\label{intro:thm:GL2xGL2}
Let $\Sigma_1,\Sigma_2$ be unitary cuspidal automorphic representations of $\GL_2(\mathbb A)$. The following conditions are equivalent:
\begin{enumerate}
\item[(a)] $L(\Sigma_1\times\Sigma_2,s)$ has a pole at $s=1$;
\item[(b)] $\Sigma_2\simeq\Sigma_1^\vee$;
\item[(c)] $\omega_{\Sigma_1}\omega_{\Sigma_2}=1$ and the diagonal period
\[
\int_{Z(\mathbb A)\GL_2(\mathbb Q)\backslash\GL_2(\mathbb A)}
\varphi_1(h)\varphi_2(h)\,dh
\]
is nonzero for some $\varphi_i\in\Sigma_i$.
\end{enumerate}
When these conditions hold, the pole of $L(\Sigma_1\times\Sigma_2,s)$ is simple.
\end{theorem}

\begin{proof}
The equivalence of (a) and (b) is the Rankin--Selberg pole theorem; see \cite{JPSS83}. The diagonal period converges for cusp forms and gives an invariant pairing. Orthogonality of cuspidal representations shows that it vanishes unless the representations are contragredient. Conversely, in the contragredient case it recovers the Petersson pairing and is nonzero for suitable vectors.
\end{proof}

For example, in the self-dual situation Artin formalism gives
\begin{equation}\label{eq:rankin-complex}
L(f\otimes f^*,j)
=L(\operatorname{ad}^0\Sigma_f,j-k+1)\,\zeta(j-k+1).
\end{equation}
The degree-three factor is the symmetric-square factor with the nebentype twist $\chi_f^{-1}$. In particular, the complex pole occurs at $j=k$, corresponding to $s=1$. This convention will be useful when comparing the complex and $p$-adic statements below.

\subsection{Weight space and ordinary families}

For the Hida-theoretic results recalled here, assume $p\geq5$. Fix embeddings $\overline{\mathbb Q}\hookrightarrow\mathbb C$ and $\overline{\mathbb Q}\hookrightarrow\overline{\mathbb Q}_p$, and a finite extension $K/\mathbb Q_p$ containing the coefficients under consideration. Put
\[
\Lambda_{\mathcal O}=\mathcal O_K[\![\mathbb Z_p^\times]\!],
\qquad \mathcal W=(\Spf\Lambda_{\mathcal O})^{\mathrm{rig}}.
\]
An arithmetic character has the form $z\mapsto z^n\varepsilon(z)$, with $\varepsilon$ of finite order on $\mathbb Z_p^\times$. The decomposition $z=\omega(z)\langle z\rangle$ separates the character on $\mu_{p-1}$ from its restriction to $1+p\mathbb Z_p$. On each component, a parameter is $T=[1+p]-1$, which specializes to $(1+p)^n\varepsilon(1+p)-1$.

We use cohomological weights for the weight map of a Hida family: the character $z\mapsto z^{k-2}\varepsilon(z)$ corresponds to modular weight $k$. A primitive cuspidal Hida family $\hf$ is a normalized ordinary eigenform \[ \hf=\sum_{n\geq1}a_n(\hf)q^n \] over a normal branch finite flat over a component of $\Lambda_{\mathcal O}$, with the classical specializations given by ordinary control. We work with primitive branches and the ordinary projectors of \cite{hida2,KLZ}. Let $U_{\hf}$ be a sufficiently small affinoid in the rigid space of the branch, and put $\mathcal O_{\hf}=\mathcal O(U_{\hf})$. Denote its weight character by $\kappa_{\hf}$ and its \emph{full modular weight character} by \[ \mathbf k_{\hf}(z)=z^2\kappa_{\hf}(z). \]
We use analogous notation for a second, independently chosen family $\hg$.

A classical point $x\in U_{\hf}$ is called crystalline if $\hf_x$ is the ordinary $p$-stabilization of a newform $\hf_x^\circ$ of level prime to $p$. At such a point of weight $k\geq2$, write \[ \alpha_{\hf_x}=a_p(\hf_x),\qquad \beta_{\hf_x}=\frac{\chi_f(p)p^{k-1}}{\alpha_{\hf_x}}. \] Here $\chi_f$ denotes the prime-to-$p$ nebentype, fixed on the chosen component. The function $\alpha_{\hf}=a_p(\hf)$ is an analytic unit. The notation $\beta_{\hf_x}$ is used only at classical points: the factor $p^{k-1}$ does not vary analytically with $k$.

Throughout, we specialize only at crystalline points of weight at least two. Around each such point, we work in an ordinary neighbourhood as in \cite[\S7.7]{KLZ}: the weight gives a local coordinate, the classical specialization maps are isomorphisms by ordinary control, and Hida's normalized eigenprojector has no poles. We keep the primitive newform data and the maps used to change level fixed throughout. The improvement and the ordinary function use these same maps. At a matching specialization the primitive tame levels coincide, and the residue calculation below uses the normalized newform pairing at that level.

\subsection{The three-variable Hida--Rankin function}

Let $U_{\mathbf j}\subset\mathcal W$ be an affinoid in cyclotomic weight space. Write \[ \begin{split} \mathcal O_{\hf\hg}&=\mathcal O_{\hf}\widehat\otimes_K\mathcal O_{\hg},\\ \mathcal O_{\hf\hg\mathbf j}&=\mathcal O_{\hf\hg}\widehat\otimes_K\mathcal O(U_{\mathbf j}). \end{split} \]
We work componentwise when taking fraction fields, enlarging $K$ to contain the required roots of unity at the tame levels. The cyclotomic arithmetic character $z\mapsto z^j$ is indexed by the motivic integer $j$. Let $I_{\hf}$ be the congruence ideal of the first family, extended to either product algebra when necessary. Our convention is that $I_{\hf}$ is an integral ideal, so the normalized ordinary eigenfunctional takes values in the fractional module $I_{\hf}^{-1}$. On our smooth local branches, this fractional module is the one denoted $I_{\mathbf a}$ in \cite[\S7.7]{KLZ}.

The interpolation formula uses the Rankin series of \cite[\S2.7]{KLZ}. To distinguish it from the primitive tensor-product function above, our convention is
\[
L(f,g,j)=L^{(N_fN_g)}(\chi_f\chi_g,2j+2-k-\ell)
\sum_{n\geq1}\frac{a_n(f)a_n(g)}{n^j},
\]
initially in its half-plane of convergence and then by meromorphic continuation. The superscript on the Dirichlet $L$-function removes the factors at primes dividing $N_fN_g$. %At primes dividing the levels this convention need not agree with the primitive tensor-product $L$-function.

\begin{theorem}[Hida]\label{thm:Hida-3var}
There exists a three-variable function
\[
L_p(\hf,\hg,\mathbf j)\in I_{\hf}^{-1}\mathcal O_{\hf\hg\mathbf j}.
\]
For crystalline points $x,y$ of full weights $k,\ell\geq2$, and integers $\ell\leq j<k$, it satisfies
\[
L_p(\hf,\hg,\mathbf j)(x,y,j)
=\frac{\mathcal E(x,y,j)}{\mathcal E_0(x)\mathcal E_1(x)}
\frac{\Gamma(j)\Gamma(j-\ell+1)L(\hf_x^\circ,\hg_y^\circ,j)}
{\pi^{2j-\ell+1}(-i)^{k-\ell}2^{2j+k-\ell}
\langle\hf_x^\circ,\hf_x^\circ\rangle},
\]
with the Petersson normalization of \cite{KLZ}, where
\begin{equation}\label{Euler:Factor:Triple}
\begin{split}
\mathcal E_0(x)&=1-\frac{\beta_{\hf_x}}{\alpha_{\hf_x}},\qquad
\mathcal E_1(x)=1-\frac{\beta_{\hf_x}}{p\alpha_{\hf_x}},\\
\mathcal E(x,y,j)&=
\left(1-\frac{p^{j-1}}{\alpha_{\hf_x}\alpha_{\hg_y}}\right)
\left(1-\frac{p^{j-1}}{\alpha_{\hf_x}\beta_{\hg_y}}\right)\left(1-\frac{\beta_{\hf_x}\alpha_{\hg_y}}{p^j}\right)
\left(1-\frac{\beta_{\hf_x}\beta_{\hg_y}}{p^j}\right).
\end{split}
\end{equation}
\end{theorem}

This is \cite[Theorem 5.1d]{hida2}, with the normalization and interpolation formula of \cite[Theorems 2.7.4 and 7.7.2]{KLZ}. Weight notation $(k,\ell,j)$ will occasionally abbreviate $(x,y,j)$, but the eigenforms, and not just their weights, are part of the specialization data.

At a crystalline point the second Euler factor can be written as
\[
1-\frac{p^{j-1}}{\alpha_{\hf_x}\beta_{\hg_y}}
=1-\frac{\chi_g^{-1}(p)\alpha_{\hg_y}}{\alpha_{\hf_x}}p^{j-\ell}.
\]
It therefore varies analytically on the two-variable surface $\mathbf j=\mathbf k_{\hg}$. This is an equality of full characters, including their finite-order parts.

\subsection{The improved function and its poles}

We continue with two arbitrary primitive ordinary families $\hf,\hg$ as above. On the subspace $\mathbf j=\mathbf k_{\hg}$ put
\[
e(x,y)=1-\frac{\chi_g^{-1}(p)\alpha_{\hg_y}}{\alpha_{\hf_x}}
\in\mathcal O_{\hf\hg}.
\]
The improved construction uses the ordinary Eisenstein series in Hida's convolution in place of its $p$-depletion, thereby removing the factor $e$. Its only possible Eisenstein pole occurs when the Eisenstein character is trivial. This requires the two full $p$-adic weight characters to agree and the fixed tame Eisenstein character to be trivial. At crystalline points, the equality of weight characters is simply $k=\ell$. We work locally near this locus.

\begin{theorem}\label{prop:rankin-improvement}
There exists a meromorphic two-variable function
\[
L_p^{\imp}(\hf,\hg)\in\Frac(\mathcal O_{\hf\hg})
\]
such that
\begin{equation}\label{eq:rankin-improvement}
\left.L_p(\hf,\hg,\mathbf j)\right|_{\mathbf j=\mathbf k_{\hg}}
=e\,L_p^{\imp}(\hf,\hg).
\end{equation}
Away from the locus where the Eisenstein character is trivial, the function belongs to $I_{\hf}^{-1}\mathcal O_{\hf\hg}$. Along that locus it has at most a simple pole, with the same congruence-ideal bound: multiplying by any local equation of the locus gives an element of $I_{\hf}^{-1}\mathcal O_{\hf\hg}$. At crystalline points of weights $k>\ell\geq2$, one has
\[
L_p^{\imp}(\hf,\hg)(x,y)
=\frac{\mathcal E^{\imp}(x,y)}{\mathcal E_0(x)\mathcal E_1(x)}
\frac{\Gamma(\ell)L(\hf_x^\circ,\hg_y^\circ,\ell)}
{\pi^{\ell+1}(-i)^{k-\ell}2^{k+\ell}
\langle\hf_x^\circ,\hf_x^\circ\rangle},
\]
where the precise improved Euler factor is
\begin{equation}\label{eq:rankin-improved-euler}
\mathcal E^{\imp}(x,y)=
\left(1-\frac{p^{\ell-1}}{\alpha_{\hf_x}\alpha_{\hg_y}}\right)
\left(1-\frac{\beta_{\hf_x}\alpha_{\hg_y}}{p^\ell}\right)
\left(1-\frac{\beta_{\hf_x}\beta_{\hg_y}}{p^\ell}\right).
\end{equation}
\end{theorem}

\begin{proof}
We use Hida's smoothed Eisenstein convolution for two independent primitive families \cite[Theorem 9.1, pp.~80--81]{hida2}. Hida writes the second input as $g_Q^\rho$, so his second family is $\hg^*$ in our conventions; its ordinary root is $\chi_g^{-1}(p)\alpha_{\hg}$. The eigenfunctional and periods are those of Theorem~\ref{thm:Hida-3var}.

To bound the order of the pole, we write the smoothing factor in local weight coordinates. Choose a smoothing integer $b$ prime to the levels and to $p$, with $u=\langle b\rangle$ and $\log_p(u)\ne0$. On the trivial Eisenstein-character component its smoothing factor is
\[
S_b=1-\frac{\kappa_{f,x}(u)}{\kappa_{g,y}(u)}.
\]
Writing $X=\kappa_{f,x}(u)-1$ and $Y=\kappa_{g,y}(u)-1$, we have $S_b=-(X-Y)/(1+Y)$. On a sufficiently small neighbourhood of the trivial Eisenstein character, this is a unit times a local equation of the weight diagonal. Dividing the smoothed convolution by $S_b$ therefore gives at most a simple pole, with the congruence-ideal bound of Hida's theorem. Away from that locus one chooses a smoothing factor which is a unit. The resulting functions agree on overlaps by their classical interpolation formulas.

Replacing the depleted Eisenstein input by the ordinary one removes the factor $e$. The remaining local factors are \eqref{eq:rankin-improved-euler}, and substitution of $j=\ell$ gives the stated interpolation formula. Since the two constructions use the same eigenfunctional and periods, their comparison gives \eqref{eq:rankin-improvement}.

\end{proof}

We now fix a crystalline point $x_0\in U_{\hf}$ of weight $k\geq2$, and write $f=\hf_{x_0}^{\circ}$. All its data, including the ordinary refinement, remain fixed. The function
\[
L_p^{\imp}(f,\hg):=L_p^{\imp}(\hf,\hg)(x_0,-)
\]
then has just one variable, the point of $U_{\hg}$. 
\begin{theorem}\label{prop:interpolation-gl2xgl2}
On the ordinary classical neighbourhoods specified above, $L_p^{\imp}(f,\hg)$ has a pole at a crystalline point $y_0$ if and only if
\[
\hg_{y_0}^{\circ}=f^*,
\qquad\text{equivalently}\qquad
\Sigma_{\hg_{y_0}^{\circ}}\simeq\Sigma_f^\vee.
\]
Such a point has weight $k$, and the pole is simple.
\end{theorem}

\begin{proof}
The polar part of the Eisenstein family is a nonzero multiple of the constant form $\mathbf1$. Let $f_\alpha$ be the fixed ordinary stabilization, and let $h$ be the ordinary stabilization of $(\hg_{y_0}^{\circ})^*$. If their primitive tame levels are $N,J$, put $L=\operatorname{lcm}(N,J)$. The equal-weight calculation in \cite[Theorem 9.1 and its proof, pp.~80--81]{hida2} expresses the residue, up to a nonzero scalar, as
\[
\ell_f\!\left(\operatorname{Tr}_{L/N}e^{\ord}(h|[L/J])\right),
\]
where $\ell_f(f_\alpha)=1$ and the trace has Hida's normalization. All these maps commute with the Hecke operators away from $Lp$. Distinct primitive eigensystems are separated by those operators, so this expression vanishes unless the two primitive forms agree. In the latter case $N=J=L$, and uniqueness of the ordinary root in weight at least two gives $h=f_\alpha$. The expression is then nonzero.

Ordinary control identifies the specializations with classical eigenspaces. The eigenfunctional is regular on our neighbourhood. The weight map is \'etale, so the possible Eisenstein pole remains simple on the fixed-$f$ fibre. Its residue vanishes exactly in the nonmatching case, when the singularity is removable.

\end{proof}

Together with Theorem~\ref{thm:GL2xGL2}, this gives the desired analogy: the one-variable improved function has a pole precisely when the complex Rankin function has its pole at $s=1$, or equivalently when the diagonal automorphic period is nonzero. In both settings, taking the Eisenstein residue leaves the classical pairing with $\mathbf1$. This is the model for the comparison in Section~\ref{sec:improved}.

In particular, if $\hg=\hf^*$, the twist identifies the two parameter spaces and
\[
L_p(\hf,\hf^*,\mathbf j)(x,y,k_y)
=\left(1-\frac{\alpha_{\hf_y}}{\alpha_{\hf_x}}\right)
L_p^{\imp}(\hf,\hf^*)(x,y).
\]
After fixing $x=x_0$, the pole occurs at $y=x_0$. This is the diagonal formulation of \cite[Theorem 5.1d']{hida2} and \cite[Theorem 2.4]{RR21}. The cancellation can be made explicit for the original two families. At a matching point, choose a local coordinate $z$ vanishing at $y_0$ and put $A(y)=\chi_g^{-1}(p)\alpha_{\hg_y}$. If $L_p^{\imp}(f,\hg)=R/z+O(1)$, then
\begin{equation}\label{eq:rankin-limit}
L_p(\hf,\hg,\mathbf j)(x_0,y_0,k)
=-\frac{1}{\alpha_{\hf_{x_0}}}
\left.\frac{dA}{dz}\right|_{y_0}R.
\end{equation}
The derivative and Laurent coefficient use the same coordinate. In particular, cancellation of the pole by the Euler factor does not by itself force the ordinary function to vanish. For $p>2$, under the hypotheses of Hida's adjoint formula \cite{Hida:Adjoint,Hida:AdjointConstancy}, the algebraic weight coordinate $k$ gives $\mathcal L(\operatorname{ad}^0 V_f)=-2\,\partial_k\log\alpha_{\hf}|_f$.

%%%%%%%%%%%%%%%%%%%%%%%%%%%%%%%%%%%%%%%%%%%%%%%%%%%%%%
% Preliminaries on the Gsp4
%%%%%%%%%%%%%%%%%%%%%%%%%%%%%%%%%%%%%%%%%%%%%%%%%%%%%%

\section{\texorpdfstring{The degree-eight $L$-function and its periods}{The degree-eight L-function and its periods}}
In this section, we recall the basic definitions and results from automorphic theory needed throughout the paper. The main result of this section is a known characterization of the poles of the degree-eight $L$-function of $\GSp_4\times\GL_2$ in terms of the nonvanishing of an automorphic period and the occurrence of the automorphic representation as a Yoshida lift. This result inspired the main theorem of our paper, which establishes an analogous characterization for improved $p$-adic $L$-functions specialized at classical weights.
\subsection{General notations and definitions}
\subsubsection{The algebraic groups}
Let $J_2$ be the $2\times 2$ antidiagonal matrix with all entries $1$ and let $J = \left(\begin{smallmatrix}&J_2\\-J_2& \end{smallmatrix}\right)$. We denote by $G$ the group scheme over $\mathbb{Z}$ whose $R$-points, for any commutative $\mathbb{Z}$-algebra $R$, are given by
\[G(R) := \GSp_4(R) = \{(g,m_g)\in \GL_4(R)\times \GL_1(R),\;s.t.\;{}^tgJg = m_gJ\}.\]
The projection \[\nu: G\to \GL_1,\;\;\;g\mapsto m_g,\]
is the similitude character. Let $H := \GL_2 \times_{\GL_1} \GL_2$, where both maps $\GL_2\to \GL_1$ are given by the determinant. We embed $H$ into $G$ via the map
\[ \iota \, : \, \left[
\begin{pmatrix}a&b\\c&d\end{pmatrix}, \begin{pmatrix}a'&b'\\c'&d'\end{pmatrix} \right]
\mapsto \begin{pmatrix}a&&&b\\&a'&b'\\&c'&d'\\c&&&d\end{pmatrix}.
\]
We sometimes write $H_i$ for the $i$-th $\GL_2$ factor of $H$.
\subsubsection{Parabolic subgroups}\label{subsec:parabolic}
Let $T$ be the diagonal torus of $G$, which is contained in $H$ and is a maximal torus in either $H$ or $G$. We write $B_G = N_GT$ for the upper-triangular Borel subgroup of $G$, where $N_G$ is the maximal unipotent subgroup of $G$. For the rest of the paper, we fix a generic character denoted by
\[\psi_{G}:[N_G]\to \C^{\times}.\]
Moreover, we denote by $P_{\Si}$ and $P_{\Kl}$ the standard Siegel and Klingen parabolics containing $B_G$, so that 
\[ P_{\Si} =
\begin{smatrix}\star & \star & \star & \star \\ \star & \star & \star & \star\\
&& \star & \star\\ && \star & \star \end{smatrix}, \qquad
P_{\Kl} =
\begin{smatrix}\star & \star & \star & \star \\  & \star & \star & \star\\
&\star & \star & \star\\ &&& \star \end{smatrix}.
\]
Let $M_{\Si}$ and $M_{\Kl}$ denote the Levi factors of the Siegel and Klingen parabolic subgroups, respectively. Both are isomorphic to $\GL_2\times \GL_1$. 
We denote by $B_{\GL_2} = T_{\GL_2}N_{\GL_2}$ the upper triangular Borel subgroup of $\GL_2$, and we fix a generic character of the form 
\[\psi_{\GL_2}:[N_{\GL_2}]\to \C^{\times}.\]
We denote $B_H := B_{\GL_2}\times_{\GL_1} B_{\GL_2}$ for the upper-triangular Borel of $H$ and write its Levi decomposition as $B_{H} = T_{H}N_{H}$, where $T_{H}$ is the maximal torus and $N_{H}$ its unipotent radical.  
\subsubsection{Flag varieties}
We write $\FL_G$ for the Siegel flag variety $P_{\Si} \backslash G$, with its natural right $G$-action. There are four orbits for the Borel $B_G$ acting on $\FL_G$, the \emph{Bruhat cells}, represented by a subset of the Weyl group of $G$, the \emph{Kostant representatives}, which are the smallest-length representatives of the quotient $W_M \backslash W_G$. We denote these by $w_0, \dots, w_3$. Note that the cell $C_{w_i} = P_{\Si} \backslash P_{\Si} w_i B_G \subset \FL_G$ has dimension $\ell(w_i) = i$.

Analogously, for the $H$-flag variety $\FL_H = B_H \backslash H$, we have 4 Kostant representatives $w_{00} = \mathrm{id}$, $w_{10} = \left(\begin{pmatrix}0&1\\-1&0\end{pmatrix}, \mathrm{id}\right)$, similarly $w_{01}$, $w_{11}$ (with the cell $C_{w_{ij}}$ having dimension $i + j$). (This is the whole of the Weyl group of $H$, since the Levi subgroup $M_H=T$ has trivial derived group.)
\subsubsection{Algebraic representations}\label{Algebraic:Representations}
Since the maximal split tori of $G$, $H$ and $M_{\Si}$ coincide, the representation theory of these groups is controlled by the same combinatorial object, the character group $X^{\bullet}(T)$. We identify it with the set of triples $(a,b;c) \in \mathbb Z^3$ such that $a+b \equiv c \pmod{2}$, by defining $\lambda(a,b;c)$ as the unique character of $T$ such that \[ \begin{smatrix} st_1 \\ & st_2 \\ && st_2^{-1} \\ &&&st_1^{-1} \end{smatrix} \mapsto t_1^{a}t_2^{b} s^c. \]
\begin{prop}
We have the following classification of the algebraic representations in terms of their highest weights:
\begin{itemize}
\item The weights $\lambda(a,b;c)$ which are dominant for $M_{\mathrm{Si}}$ are precisely those satisfying $a \geq b$. For such triple, we denote by $W_{G}(a,b;c)$ the irreducible representation of $M_{\Si}$ with highest weight $\lambda(a,b;c)$. Since $M_{\Si}\simeq \GL_2\times \GL_1$, we have
\[W_G(a,b;c)\simeq (\mathrm{Sym}^{a-b}\otimes \mathrm{det}^b)\boxtimes \nu^{(c-a-b)/2}.\]
\item If, in addition, $b\geq 0$ then $\lambda(a,b;c)$ is dominant for $G$. We denote by $V_G(a,b;c)$ the irreducible representation with highest weight $\lambda(a,b;c)$.
\item If $a,b\geq 0$, then $\lambda(a,b;c)$ is a highest weight for $H$. We write $V_H(a,b;c)$ for the corresponding irreducible representation of $H$.
\end{itemize}
\end{prop}
\begin{lemma}\label{Restriction:WG}
The restriction of $W_G(a,b;c)$ to $B_H$ is
\[
W_G(a,b;c)|_{B_H}=\bigoplus_{j=0}^{a-b}\lambda(a-j,b+j;c).
\]
Indeed, the monomials of $\Sym^{a-b}$ give the displayed torus weights, and $N_H\subset N_{\Si}$ acts trivially on the inflated Levi representation.
\end{lemma}
\subsubsection{Representations of the maximal compact subgroup}
Let
\[
G(\R)^+=\{g\in G(\R):\nu(g)>0\}
\]
be the identity component of $G(\R)$. Denote by $K_{G,\infty}=\R^{\times}\cdot\; U(2)$ the maximal compact
modulo centre subgroup of $G(\R)^+$. Let
\[
T_U=\{\operatorname{diag}(z_1,z_2)\in U(2): |z_1|=|z_2|=1\}
\]
be the diagonal maximal torus of $U(2)$. We identify its character
lattice with $\mathbb{Z}^2$ by associating $(A,B)$ with the character
$\operatorname{diag}(z_1,z_2)\mapsto z_1^A z_2^B$, and choose $(1,-1)$
as the positive root. For integers $A\ge B$, the irreducible
$U(2)$-representation of highest weight $(A,B)$ is
\[
V_{A,B}=\operatorname{Sym}^{A-B}(\C^2)\otimes\mathrm{det}^B.
\]
Its weights are $(A-j,B+j)$, for $0\le j\le A-B$, each occurring
with multiplicity one.
\begin{rk}
Since
\[
K_{G,\infty}\simeq
(\R^{\times}\times U(2))/\langle(-1,-I_2)\rangle,
\]
an irreducible finite-dimensional complex representation of
$K_{G,\infty}$ is determined by a representation $V_{A,B}$ of $U(2)$
and a character $\chi:\R^{\times}\to\C^{\times}$ satisfying
\[
\chi(-1)=(-1)^{A+B}.
\]
\end{rk}
\begin{defn}
If $\chi(z)=z^C$ for some $C\in \mathbb{Z}$, necessarily satisfying
$C\equiv A+B\pmod{2}$, we say that the resulting representation of $K_{G,\infty}$ has highest weight $(A,B;C)$.
\end{defn}
\begin{rk}\label{Representation:Siegel:Compact}
By \cite[\S 6.3]{LPSZ1}, there is an isomorphism
\[
M_{\mathrm{Si},\C}\simeq K_{G,\infty,\C}.
\]
Thus every algebraic representation of $M_{\mathrm{Si}}$ can be
regarded as a complex representation of $K_{G,\infty}$, and we use
the same notation for both realizations. With the weight conventions
of \cite[Remark 2.4]{Horawa:Prasanna}, the representation
$W_G(A,B;C)$ of $M_{\mathrm{Si}}$ corresponds to the
$K_{G,\infty}$-representation of highest weight $(-B,-A;C)$. Moreover, under this identification, a vector of weight $(A,B;C)$ becomes a vector of weight $(-A,-B;C)$.
\end{rk}
\subsubsection{Discrete series representations}\label{Discrete:Series:Subsection}
We identify a maximal compact subgroup of $\mathrm{Sp}_4(\R)$ with $U(2)$ and use the two diagonal characters $e_1,e_2$ of its compact torus. For $\lambda_1>\lambda_2>0$, the discrete-series packet consists of the generic representations
\[
D(\lambda_1,-\lambda_2),\qquad D(\lambda_2,-\lambda_1),
\]
and the holomorphic and antiholomorphic representations
\[
H(\lambda_1,\lambda_2),\qquad H(-\lambda_2,-\lambda_1).
\]
The pairs displayed are their Harish--Chandra parameters; see \cite[Theorem 9.20]{Knapp:Representations} and \cite[\S1.2]{Moriyama:ZetaIntegrals}.

For a central character $\omega:\R^\times\to\C^\times$ with $\omega(-1)=(-1)^{\lambda_1+\lambda_2+1}$, the two essentially discrete-series representations of $\GSp_4(\R)$ are characterized by
\begin{align*}
\Pi^W(\lambda_1,\lambda_2)_\omega|_{\mathrm{Sp}_4(\R)}
 &=D(\lambda_1,-\lambda_2)\oplus D(\lambda_2,-\lambda_1),\\
\Pi^H(\lambda_1,\lambda_2)_\omega|_{\mathrm{Sp}_4(\R)}
 &=H(\lambda_1,\lambda_2)\oplus H(-\lambda_2,-\lambda_1).
\end{align*}
Thus $\Pi^W$ is generic and $\Pi^H$ is nongeneric. Our cohomological weights satisfy $k_1\geq k_2\geq3$, so the archimedean representations used below lie in the discrete-series range.

\subsubsection{Hecke parameters}

Let $\Pi$ be a unitary globally generic cuspidal automorphic representation of $G(\mathbb A)$, unramified at $p$, with $\Pi_{\infty}$ of weights $(k_1,k_2)$, with $k_1 \geq k_2 \geq 3$. Let us put $ r_i := k_i - 3$, for $i=1,2$. Let $\{\alpha,\beta,\gamma,\delta\}$ be the Hecke parameters of $\Pi_p$ in the sense of \cite[\S5.1.5]{Pil20} or \cite[Theorem 10.1.3]{LSZ22}, where we borrow our conventions. More precisely, they are the reciprocal of the roots of a precise Hecke polynomial at $p$, $P_p(X)$, defined in Theorem 10.1.3 (3) of loc.\,cit.
This means that we can normalize $(\alpha,\beta,\gamma,\delta)$ in such a way that all four are algebraic integers of complex absolute value $p^{(k_1+k_2-3)/2}$, and ordered such that $\alpha \delta = \beta \gamma = p^{(k_1+k_2-3)} \chi_{\Pi}(p)$ and $0 \leq v_p(\alpha) \leq \ldots \leq v_p(\delta)$. 
Here $\chi_\Pi$ denotes the Dirichlet character associated with the central character of $\Pi$. For general type, purity follows from the Ramanujan theorem for non-CAP Siegel cusp forms \cite{Weissauer:Endoscopy}. In the Yoshida case the same purity follows from the two elliptic constituents: the roots for the form of weight $k_1+k_2-2$, and $p^{k_2-2}$ times the roots for the form of weight $k_1-k_2+2$, all have this absolute value by \cite{Del74}.

Let $\Iw_G(p)$ denote the Iwahori subgroup of $G$ associated with $B_G$. We consider the following operators in the Hecke algebra of level $\Iw_G(p)$, acting on the cohomology of any of the sheaves introduced above:

\begin{itemize}
\item The Siegel operator $\mathcal U_{\Si} = [\diag(p,p,1,1)]$, as well as its dual $\mathcal U_{\Si}' = [\diag(1,1,p,p)]$.
\item The Klingen operator $\mathcal U_{\Kl} = p^{-r_2} \cdot [\diag(p^2,p,p,1)]$, as well as its dual $\mathcal U_{\Kl}' = p^{-r_2} \cdot [\diag(1,p,p,p^2)]$.
\item The Borel operator $\mathcal U_B = \mathcal U_{\Si} \cdot \mathcal U_{\Kl}$, as well as its dual $\mathcal U_B' = \mathcal U_{\Si}' \cdot \mathcal U_{\Kl}'$.
\end{itemize}

A Hecke eigenform is called \emph{Siegel} (resp.~\emph{Klingen}) ordinary if the corresponding eigenvalue is a $p$-adic unit. It is called \emph{Borel} ordinary if it is both Siegel and Klingen ordinary. If $\Pi_p$ is Borel-ordinary, it is both Klingen-ordinary and Siegel-ordinary. In this case, we may and will assume that $(\alpha,\beta,\gamma,\delta)$ have $p$-adic valuations $(0,\,k_2-2,\,k_1-1,\,k_1+k_2-3)$, respectively.
\subsection{Automorphic representations and $L$-functions}
In this subsection, we introduce the degree-eight $L$-function of $\GSp_4\times \GL_2$ and recall its integral representation due to Novodvorsky. We then relate the pole of this $L$-function to an integral of a cusp form on $\GSp_4\times \GL_2$ over $H$, and recall that this integral does not vanish if and only if $\Pi$ is a Yoshida lift. These results are summarized in the main theorem of this section, which is already known and serves as the complex counterpart of the main theorem of this paper. Indeed, this classical result inspired our main theorem, in which we seek to establish an analogous characterization for improved $p$-adic $L$-functions specialized at classical weights.
\subsubsection{Novodvorsky's integral representation of the degree-eight $L$-function of $\GSp_4\times \GL_2$}\label{SectionSpinStandardRep}
From now on, let $\Pi$ be a unitary, globally generic, cuspidal automorphic representation whose archimedean component $\Pi_\infty$ is the generic discrete series described before, and let $\Sigma$ be a unitary cuspidal automorphic representation of $\GL_2(\A)$. Write $\omega=\omega_\Pi\omega_\Sigma$. For $\varphi\in\mathcal S(\A^2)$, set
\begin{equation}\label{Section:f:Archimedean}
f_{\varphi,\omega}(g,s):=|\mathrm{det} (g)|^s\int_{\A^\times}\varphi((0,t)g)\omega(t)|t|^{2s}\,d^\times t.
\end{equation}
This is a section of the normalized induction
\[
\Ind_{B_{\GL_2}(\A)}^{\GL_2(\A)}
\bigl(|\cdot|^{s-1/2}\boxtimes\omega^{-1}|\cdot|^{1/2-s}\bigr).
\]
Its Eisenstein series
\[
E(g,s,\varphi)=\sum_{\gamma\in B_{\GL_2}(\Q)\backslash\GL_2(\Q)}f_{\varphi,\omega}(\gamma g,s)
\]
converges for $\Re(s)$ sufficiently large and has meromorphic continuation. Novodvorsky's global zeta integral is defined as follows:
\[
Z(s,\Psi,F,\varphi):=\int_{Z_G(\A)H(\Q)\backslash H(\A)}
\Psi(\iota(h))F(h_2)E(h_1,s,\varphi)\,dh.
\]
When $\omega = 1$, we define the period
\[
\mathcal I(\Psi,F):=\int_{Z_G(\A)H(\Q)\backslash H(\A)}\Psi(\iota(h))F(h_2)\,dh.
\]
\begin{prop}\label{Period:Equals:Residue}
Assume $\omega=1$. The only possible poles of $Z(s,\Psi,F,\varphi)$ are simple poles at $0,1$, and
\[
\operatorname{Res}_{s=1}Z(s,\Psi,F,\varphi)
\overset{.}{=}\widehat\varphi(0)\mathcal I(\Psi,F),\qquad
\widehat\varphi(0)=\int_{\A^2}\varphi(x,y)\,dx\,dy.
\]
\end{prop}
\begin{proof}
This is \cite[Propositions 2.5 and 3.1]{Moriyama:ZetaIntegrals}.
\end{proof}
We denote by $c_E$ the nonzero Eisenstein residue constant for the chosen measures, so that $\operatorname{Res}_{s=1}Z(s,\Psi,F,\varphi)=c_E\widehat\varphi(0)\mathcal I(\Psi,F)$.
Fix compatible Whittaker models for $\Pi$ and $\Sigma$. For factorizable data, \cite[\S3.2]{Moriyama:ZetaIntegrals} and \cite{Shapiro:Soudry} show
\begin{equation}\label{eq:global-unfolding}
Z(s,\Psi,F,\varphi)=\prod_vZ_v(s,W_{\Psi_v},W_{F_v},\varphi_v),
\end{equation}
where
\[
Z_v(s,W,W',\varphi_v)=\int_{Z_G(\Q_v)N_H(\Q_v)\backslash H(\Q_v)}
W(\iota(h))W'(h_2)f_{\varphi_v,\omega_v}(h_1,s)\,dh.
\]

\begin{prop}\label{GCD:Property:Zeta}
Let $v$ be a finite place. As the local data $(W,W',\varphi_v)$ vary, the integrals $Z_v(s,W,W',\varphi_v)$ generate the fractional ideal
\[
L^{\mathrm{Nov}}(s,\Pi_v\times\Sigma_v)
\C[q_v^s,q_v^{-s}],
\]
where $L^{\mathrm{Nov}}(s,\Pi_v\times\Sigma_v)
=P_v(q_v^{-s})^{-1}$, with $P_v(0)=1$.
Moreover, for any $s_0\in\C$, there exist local data $(W,W',\varphi_v)$ such that
\[
\left.
\frac{Z_v(s_0,W,W',\varphi_v)}
{L^{\mathrm{Nov}}(s_0,\Pi_v\times\Sigma_v)}
\right|_{s=s_0}\neq 0.
\]
\end{prop}
\begin{proof}
The first assertion is due to \cite{Soudry}. For the second, observe that the quotients in the statement generate the unit ideal in $\C[q_v^s,q_v^{-s}]$. Hence, they cannot all vanish at $s=s_0$.
\end{proof}
\begin{defn}\label{Definition:L:Function}
Let $\mathrm{spin}:\GSp_4(\C)\hookrightarrow\GL_4(\C)$ denote the spin representation, and let $\mathrm{std}_2$ be the standard representation of $\GL_2(\C)$. For each place $v$, let $L(s,\Pi_v\times\Sigma_v)$ denote the local Langlands $L$-factor associated with $\mathrm{spin}\boxtimes\mathrm{std}_2$. We define the completed degree-eight $L$-function by
\[
\Lambda(s,\Pi\times\Sigma)
=\prod_v L(s,\Pi_v\times\Sigma_v).
\]
\end{defn}

\begin{assumption}\label{eq:local-compatibility}
Throughout the remainder of the paper, we assume that
\[
L^{\mathrm{Nov}}(s,\Pi_v\times\Sigma_v)
= L(s,\Pi_v\times\Sigma_v)
\]
at every finite place $v$.
\end{assumption}
\begin{rk}\label{Remark:L:functions:coincide}
By \cite[\S3.4]{Moriyama:ZetaIntegrals},
\cite[Proposition 5.1]{Moriyama:ZetaIntegrals} and
\cite{Loeffler:Nov}, Assumption \ref{eq:local-compatibility}
is known to hold at a finite place $v$ whenever $\Sigma_v$ is non-supercuspidal.
\end{rk}
\begin{cor}\label{Global:Zeta:Integral:Result}
Under Assumption \ref{eq:local-compatibility}, taking factorizable $\Psi$, $F$ and $\varphi$, we have that
\[
Z(s,\Psi,F,\varphi)=\Lambda(s,\Pi\times\Sigma)\,\prod_{v\in S}Q_v(s),
\qquad Q_v(s)=\frac{Z_v(s,W_{\Psi_v},W_{F_v},\varphi_v)}{L(s,\Pi_v\times\Sigma_v)},
\]
where $S$ contains the infinite place, the ramified finite places and every finite place with nonspherical data.
\end{cor}
\begin{proof}
Combine \eqref{eq:global-unfolding} with Proposition \ref{GCD:Property:Zeta} and Remark \ref{Remark:L:functions:coincide}.
\end{proof}

\subsubsection{Yoshida lifts and distinguished periods}
Let $V$ be a four dimensional totally split quadratic space. Then
\[
\mathrm{GSO}(V)\simeq
(\GL_2\times\GL_2)/\{(zI_2,z^{-1}I_2)\}.
\]
Then, given two automorphic representations $\Sigma_1$ and $\Sigma_2$ of $\GL_2$, their external tensor product $\Sigma_1\boxtimes\Sigma_2$ defines an automorphic representation of $\mathrm{GSO}(V)$ if and only if $\omega_{\Sigma_1}=\omega_{\Sigma_2}$. In this case, we denote the theta lift of $\Sigma_1\boxtimes \Sigma_2$ to $\GSp_4$ by $\mathcal Y(\Sigma_1\boxtimes\Sigma_2)$. When nonzero and cuspidal, this lift has central character $\omega_{\Sigma_1}=\omega_{\Sigma_2}$, and its functorial lift to $\GL_4$ via the spin representation is the isobaric sum $\Sigma_1\boxplus\Sigma_2$.
\begin{lemma}\label{Poles:If:Yoshida}
The function $\Lambda(s,\Pi\times\Sigma)$ has a pole at $s=1$ if and only if
\[
\Pi\simeq\mathcal Y(\Sigma^\vee\boxtimes\Sigma')
\]
for some unitary cuspidal automorphic representation $\Sigma'\not\simeq \Sigma^{\vee}$ of $\GL_2$ satisfying $\omega_{\Sigma'}=\omega_\Sigma^{-1}$. Moreover, whenever such a pole exists, it is simple.
\end{lemma}
\begin{proof}
By \cite[Theorem 2.4]{AsgariShahidi}, the functorial lift of $\Pi$ to $\GL_4$ is either cuspidal or the isobaric sum $\Sigma_1\boxplus\Sigma_2$ of two distinct unitary cuspidal representations of $\GL_2$. Moreover, the latter case occurs if and only if $\Pi\simeq\mathcal Y(\Sigma_1\boxtimes\Sigma_2)$. In this case, the computation of the local Langlands parameters of the Yoshida lift in \cite{GT11} gives
\[
\Lambda(s,\Pi\times\Sigma)
=\Lambda(\Sigma_1\times\Sigma,s)
\Lambda(\Sigma_2\times\Sigma,s).
\]
By \cite{JPSS83}, each factor $\Lambda(\Sigma_i\times\Sigma,s)$ has a pole at $s=1$ if and only if $\Sigma_i\simeq\Sigma^\vee$, in which case the pole is simple. Since $\Sigma_1$ and $\Sigma_2$ are distinct, at most one factor can have a pole, while the other is nonzero at $s=1$ by \cite{Sha81}.
\end{proof}

\begin{theorem}\label{thm:complex-period}
Assume $\omega_\Pi\omega_\Sigma=1$. The following are equivalent:
\begin{enumerate}
\item $\Lambda(s,\Pi\times\Sigma)$ has a simple pole at $1$,
\item $\Pi\simeq\mathcal Y(\Sigma^\vee\boxtimes\Sigma')$ for a cuspidal $\Sigma'\not\simeq\Sigma^\vee$ with $\omega_{\Sigma'}=\omega_\Sigma^{-1}$,
\item the period $\mathcal I$ is not identically zero.
\end{enumerate}
\end{theorem}
\begin{proof}
The equivalence of (1) and (2) follows from Lemma~\ref{Poles:If:Yoshida}. For (1) and (3), recall that Corollary \ref{Global:Zeta:Integral:Result} shows
\[
Z(s,\Psi,F,\varphi)=\Lambda(s,\Pi\times\Sigma)\prod_{v\in S}Q_v(s),
\]
At finite places, $Q_v(1)$ are holomorphic at $s=1$ and can be chosen nonzero by Proposition \ref{GCD:Property:Zeta} and Assumption \ref{eq:local-compatibility}. At infinity, holomorphy near $s=1$ follows from the tempered Whittaker bounds of \cite{Wal92}, while nonvanishing for suitable data follows from \cite[Proposition 5.1]{Moriyama:ZetaIntegrals}. Consequently, $\Lambda(s,\Pi\times\Sigma)$ has a pole at $s=1$ if and only if some global zeta integral has a nonzero residue there. By Proposition \ref{Period:Equals:Residue}, this is equivalent to the nonvanishing of $\mathcal I$.
\end{proof}

\section{Automorphic periods and cohomology}\label{subsec:coherent-periods}
In this section, we use Serre duality to construct a pairing between the cohomology of the Shimura varieties attached to $G$ and $H$. We then express this pairing as an integral over $[H]$ involving an automorphic form on $G$ and a pair of modular forms. A similar pairing was constructed in \cite{LPSZ1} for modular forms of strictly positive weight. Here, we introduce a different pairing that allows us to consider modular forms of weight $0$. This construction will play a key role in the next section, where it provides a link between the $p$-adic pairing of \cite{LPSZ1} and the complex geometry in the weight $0$ setting, which is essential for establishing our main result.
\subsection{Shimura varieties and their cohomology}\label{Section:Shimura:Varieties}
In this subsection, we recall the basic theory of Shimura varieties and introduce automorphic line bundles, constructed from representations of maximal compact (and equivalently Levi) subgroups. These bundles encode arithmetic information and play a central role in the construction of the $p$-adic $L$-function in the next section. We conclude by recalling Serre duality for Shimura varieties.
\subsubsection{The Shimura varieties}\label{subsubsec:shimura-varieties}
Let $K_f\subset G(\A_f)$ be a neat open compact subgroup, and let $Y_{G,\Q}$ denote the canonical model over $\Q$ of the Siegel threefold of level $K_f$. This is a smooth quasi-projective variety whose complex points admit the analytic uniformization
\[
Y_G(\C)\simeq G(\Q)_+\backslash \left(\mathcal{H}_2\times G(\A_f)/K_f\right),
\]
where $\mathcal{H}_2$ is the Siegel upper half-space of genus $2$, and $G(\Q)_+$ denotes the subgroup of elements with positive similitude. Fix a toroidal compactification $X_{G,\Q}$ of $Y_{G,\Q}$, with boundary divisor $D_G$.

Let $Y_{H,\Q}$ denote the canonical model over $\Q$ of the Shimura variety associated with $H$, of level $K_{H,f}=K_f\cap H(\A_f)$. In the particular case where $K_f$ is the principal congruence subgroup of level $N\geq 3$, namely
\[
K_f\simeq\ker\left(\GSp_4(\widehat{\mathbb Z})\longrightarrow\GSp_4(\mathbb Z/N\mathbb Z)\right),
\]
the Shimura variety $Y_{H,\Q}$ admits a description in terms of modular curves. More precisely, letting $Y(N)$ denote the full-level modular curve over $\Q$, and $X(N)$ its smooth projective compactification, we have
\begin{equation}\label{Union:Modular:Curves}
Y_{H,\Q}\simeq Y(N)\times_{\mu_N}Y(N),
\qquad
X_{H,\Q}\simeq X(N)\times_{\mu_N}X(N),
\end{equation}
where the latter isomorphism is compatible with the compactifications. We use principal congruence level for the product description \eqref{Union:Modular:Curves}; at other levels, the required maps are obtained through a common finite cover. We write
\[
D_{X(N)}:=(X(N)\setminus Y(N))_{\mathrm{red}}.
\]
The inclusion $\iota:H\hookrightarrow G$ induces a natural morphism $\iota:Y_{H,\Q}\to Y_{G,\Q}$. We choose toroidal compactifications $X_{H,\Q}$ and $X_{G,\Q}$ compatibly, so that $\iota$ extends to a morphism
\[
\iota:X_{H,\Q}\longrightarrow X_{G,\Q},
\]
as explained in \cite[\S 2.3.4]{LPSZ1}.
\subsubsection{Archimedean components and the cohomology of Shimura varieties}\label{Archimedean:Component:Siegel}
Let $\Pi$ be a unitary globally generic automorphic representation of $G(\mathbb A)$, and fix integers $r_1\geq r_2\geq0$. The cohomological discrete-series packet with derived infinitesimal character $(r_1+2,r_2+1)$ consists of
\[
\Pi^W(r_1+2,r_2+1)\quad\text{and}\quad\Pi^H(r_1+2,r_2+1),
\]
normalized to have central character $\operatorname{sign}^{r_1+r_2}$ on $\mathbb R^\times$; see \cite{Vogan:Zuckerman} and \cite[\S5.1]{LPSZ1}. The coefficient on $\operatorname{Sp}_4$ has highest weight $(r_1,r_2)$. An algebraic extension to $G$ is $V_G(r_1,r_2;m)$, where $m=r_1+r_2$; its real normalization $V_G(r_1,r_2;m)\otimes|\nu|^{-m/2}$ has central character $\operatorname{sign}^m$. This real normalization need not be algebraic. We use below the coherent bundles $[L_1]$ and $[L_2]$ of \cite[\S5.2]{LPSZ1}, for which the finite automorphic normalization is $\Pi_f\otimes|\nu|_f^{-m/2}$. Tensoring these bundles by $C_{u_\diamond}$ gives the normalization in \eqref{Normalization:Pi}.

We assume from now on that $\Pi_\infty\simeq\Pi^W(r_1+2,r_2+1)$. Its restriction to $G(\mathbb R)_+$ decomposes as
\[
\Pi^{W}(r_1+2,r_2+1)|_{G(\R)_{+}}\simeq \Pi_1\oplus\Pi_2,
\]
where $\Pi_1$ and $\Pi_2$ have minimal $U(2)$-types of highest
weights $(r_1+3,-r_2-1)$ and $(r_2+1,-r_1-3)$, respectively.
\begin{defn}\label{Minimal:U:Types}
We denote by
\[
\tau_1:=\Sym^{r_1+r_2+4}U\otimes \mathrm{det}^{-r_2-1},\qquad \tau_2:=\Sym^{r_1+r_2+4}U\otimes \mathrm{det}^{-r_1-3},
\]
the minimal $U(2)$-types of $\Pi_1$ and $\Pi_2$ respectively, where $U$ is the standard representation of $U(2)$.
\end{defn}

\subsubsection{Automorphic vector bundles}
Following \cite[\S 2.4.2]{LPSZ1}, an algebraic representation $V$ of $M_{\mathrm{Si}}$ (resp. $T$) gives rise to an automorphic vector bundle on $Y_{G,\Q}$ (resp. $Y_{H,\Q}$), whose canonical extension to $X_{G,\Q}$ (resp. $X_{H,\Q}$) we denote by $[V]$. These constructions are compatible with duals and tensor products, in the sense that
\[
[V_1]^\vee \simeq [V_1^\vee], \qquad
[V_1]\otimes [V_2]\simeq [V_1\otimes V_2].
\]
Moreover, any morphism $V_1\to V_2$ of $P_{\mathrm{Si}}$-representations (resp. $B_H$-representations) induces a morphism of the corresponding automorphic vector bundles $[V_1]\to [V_2]$. For $J\in\{G,H\}$, we denote by $[V](-D_J)$ the subcanonical extension of $[V]$. These constructions are also compatible with the embedding $\iota:X_{H,\Q}\to X_{G,\Q}$. More precisely, we have a natural identification and a morphism
\begin{equation}\label{Pullback:Vector:Bundles}
\iota^*[V]\simeq [V|_{B_H}]\qquad \iota^*\bigl([V](-D_G)\bigr)
\longrightarrow [V|_{B_H}](-D_H).\end{equation}
\begin{defn}
We now introduce the automorphic vector bundles on $X_G$ and $X_H$ that will play a central role throughout this paper.
\begin{itemize}
\item For even $u\in\mathbb{Z}$, denote by $C_u$ the character bundle $[W_G(0,0;u)]$ on $X_G$, or $[\lambda(0,0;u)]$ on $X_H$.
\item Let $r_1,r_2\in\mathbb{Z}$, and set $m:=r_1+r_2$. We define the representations
\[
L_1:=W_G(r_1+3,1-r_2;m),\qquad
L_2:=W_G(r_2+2,-r_1;m),
\]
and denote by $[L_1]$ and $[L_2]$ the corresponding automorphic vector bundles on $X_{G,\Q}$. 
\item For integers $a,b,c$ satisfying $a+b\equiv c\pmod{2}$, we define the automorphic vector bundles on $X_{H,\Q}$
\[
\omega_H(a,b;c):=[\lambda(a,b;c)],
\qquad
\omega(a,b):=\omega_H(a,b;a+b-4).
\]
\end{itemize}
\end{defn}
\begin{rk}\label{Remark:Trivialization}
The bundles $C_u$ are isomorphic to the trivial line bundle on
$X_G$ and $X_H$, respectively. However, these identifications are
not, in general, Hecke-equivariant. More precisely, the induced
Hecke action is twisted by the character
\begin{equation}\label{eq:character-line-action}
    \chi_u(g)=|\nu(g)|_f^{-u/2}.
\end{equation}
We denote the corresponding one-dimensional Hecke module by
$E\{u/2\}$, and fix a basis vector $e_{u/2}$, normalized so that
$e_s\otimes e_t=e_{s+t}$.
\end{rk}
\subsubsection{Serre duality}\label{Subsection:Serre:Duality}
First let us recall Serre duality for the varieties $X_{H}$ and $X_{G}$ and automorphic line bundles. Let $E$ be a number field and let $X$ denote either $X_{G,E}$ or $X_{H,E}$. We write $D_X$ for the boundary divisor, and set $d_X := \mathrm{dim}(X)$. Let $V$ be an algebraic representation of $M_{\Si}$ in the $G$-case or of $T$ in the $H$-case. Let $K_X$ denote the canonical bundle of $X$. Serre duality is induced by the natural contraction map
\begin{equation}\label{contraction:Map}
[V]\otimes \left([V]^{\vee}\otimes K_X\otimes C_{u}\right)
\to K_X\otimes C_{u},
\end{equation}
which, together with the cup product and the algebraic trace map $\mathrm{Tr}_{X}$, gives a perfect pairing
\begin{equation}\label{Serre:Duality:Map}
\langle \cdot,\cdot\rangle_{X}:H^{d_X-i}(X,[V])\times
H^{i}(X,[V]^{\vee}\otimes K_X\otimes C_{u})
\to E\{u/2\},\quad (a,b)\mapsto \mathrm{Tr}_X(a\cup b).
\end{equation}
Here $E\{u/2\}$ is the field $E$ equipped with the action defined in Remark \ref{Remark:Trivialization}. We will say that the vector bundle $\left([V]^{\vee}\otimes K_X\otimes C_u\right)$ is Serre-dual to $[V]$.
Kodaira--Spencer \cite[\S 2.4.2]{LPSZ1} identifies the canonical bundles as
\begin{equation}\label{Kodaira:Spencer}
K_{X_G}=[W_G(3,3;0)](-D_G),\qquad K_{X_{H}}=\omega_H(2,2;0)(-D_H).
\end{equation}
\subsection{A pairing between the cohomology of $X_{H}$ and $X_{G}$}\label{Section:Identification:Modular:Forms}
In the preceding subsection, we described Serre duality separately for the Shimura varieties $X_G$ and $X_H$. We now explain how the map $\iota$ induces pullback maps on cohomology, allowing us to construct pairings between cohomology groups of these two varieties. Our construction is similar to that of \cite{LPSZ1}, but allows us to consider modular forms of the form $1\boxtimes g$, which lie outside the range considered in \textit{loc.\ cit.} It is not, however, an extension of their pairing, since the cohomology groups involved are different. Rather, it allows us to interpret the specialization at weights $(0,\bullet)$ of the $p$-adic pairing introduced in \S\ref{subsec:family-duality} as a pairing on the cohomology of Shimura varieties. As we will see in the next section, this cohomological pairing admits an expression as an integral of automorphic forms. Thus, our construction provides a link between the specialized $p$-adic pairing and an automorphic integral.
\subsubsection{The definition of the pairing $P_{\diamond}$}
\begin{defn}\label{Definition:Vector:Bundles}
	Let $r_1\geq r_2\geq 0$ be integers, and set
	$c_0:=r_1-r_2+2$.
	We define the following automorphic vector bundles:
	\begin{enumerate}
	\item Over $X_G$, we set $u_F=2r_2+2$, $u_D=2r_1+4$.
    We define the representations and their associated vector bundles by
    \[
    V_F=L_2\otimes \lambda(0,0,u_F),\quad V_D=L_1\otimes \lambda(0,0,u_D),\quad
    E_F=[V_F],\quad E_D=[V_D].
    \]
		
		\item Over $X_H$, choose
        $(a_F,b_F),(a_D,b_D)\in\mathbb Z^2$ satisfying
        \[
        a_F+b_F=c_0,\quad -r_2\leq a_F\leq r_1+2,\quad
        b_D-a_D=c_0,\quad -r_1-1\leq a_D\leq r_2+1,
        \]
        and define
        \[
        H_F=\omega(a_F,b_F),\qquad
        H_D=\omega(a_D,2-b_D).
        \]
        We also denote by $W_{F}$ and $W_D$ the representations underlying these two automorphic vector bundles.
	\end{enumerate}
\end{defn}
\begin{rk}
The twists by $\lambda(0,0,u_\diamond)$ put both pairings in the common target $E\{m\}$. For positive first elliptic weight, they agree with the scalar pairings of \cite{LPSZ1,LR2,LR3} after the tensor transport of Theorem~\ref{thm:scalar-period-comparison}. The same algebraic contractions remain defined at weight zero, and Lemma~\ref{lem:boundary-branching} identifies their analytic specializations.
\end{rk}
From now on, let $\diamond\in\{F,D\}$, and set $i_F:=2$ and $i_D:=1$. Our aim is to construct pairings of the form
\[
P_{\diamond}:H^{i_{\diamond}}(X_G,E_{\diamond}(-D_G))\times
H^{2-i_{\diamond}}(X_H,H_{\diamond})
\longrightarrow E\{m\},
\]
for $H_{\diamond}$ with $a_{\diamond}$ and $b_{\diamond}$ in the regions stated in Definition \ref{Definition:Vector:Bundles}. We say that we are in region $(\diamond)$ when working with the vector bundles $E_{\diamond}$ and $H_{\diamond}$. By \eqref{Pullback:Vector:Bundles}, there is a pullback map
\[
\iota^*:H^{i_\diamond}(X_G,E_\diamond(-D_G))
\to
H^{i_\diamond}(X_H,[V_{\diamond}|_{B_H}]
(-D_H)).
\]
Set
\[
(A_{\diamond},B_{\diamond};C_{\diamond})=
\begin{cases}
(r_2+2,-r_1;m+u_F), & \diamond=F,\\
(r_1+3,1-r_2;m+u_D), & \diamond=D,
\end{cases}
\]
the highest weight of $V_{\diamond}$. Then, by Lemma \ref{Restriction:WG}, we have 
\begin{equation}\label{RestrictionLi}[V_{\diamond}|_{B_H}]\simeq \bigoplus_{j = 0}^{m+2}\omega_{H}(A_{\diamond}-j,B_{\diamond}+j;C_{\diamond}),
\end{equation}
where we recall that $m = r_1+r_2$. Hence, projection onto the $j$-th constituent induces
\begin{equation}\label{First:Projection}
p_j^*:H^{i_\diamond}(X_G,E_\diamond(-D_G))
\to
H^{i_\diamond}(X_H,
\omega_{H}(A-j,B+j;C)(-D_H)).
\end{equation}
Our goal is now to pair the image of \eqref{First:Projection}
with classes in $H^{2-i_\diamond}(X_H,H_\diamond)$ via Serre duality.
This requires choosing a projection $p_{j_{\diamond}}^*$ to certain 
$\omega_H(A-j_{\diamond},B+j_{\diamond};C)(-D_H)$ that is Serre dual to
$H_\diamond$, i.e.
\begin{equation}\label{Condition:Serre:Dual}
\omega_H(A-j_{\diamond},B+j_{\diamond};C)(-D_H)
\simeq H_\diamond^\vee\otimes K_{X_H}\otimes C_{m'_{\diamond}},
\end{equation}
for some $m'_\diamond\in 2\mathbb{Z}$. Such a choice determines a contraction map
\begin{equation}\label{Contraction:After:Projection}
p_{j_{\diamond}}^*E_{\diamond}(-D_G)\otimes H_{\diamond}
\longrightarrow K_{X_H}\otimes C_{m'_{\diamond}},
\end{equation}
which induces a pairing on cohomology. A direct computation using \eqref{Kodaira:Spencer} and \eqref{Condition:Serre:Dual} shows that, for fixed $r_1$ and $r_2$, the index $j_{\diamond}$ and the integer $m'_\diamond$ are uniquely determined. Their values, together with the corresponding vector bundles, are given in the following table:
\begin{equation}\label{Table:Images:p}
\begin{array}{c|c|c|c|c}
& j_{\diamond} & m'_{\diamond}
& \text{Selected summand of $\iota^*E_{\diamond}$}
& \text{Vector bundle on $X_H$}\\ \hline
F & r_2+a_F & 2m
& S_F:=\omega_H(2-a_F,2-b_F;m+u_F)
& H_F=\omega(a_F,b_F)\\
D & r_1+1+a_D & 2m
& S_D:=\omega_H(2-a_D,b_D;m+u_D)
& H_D=\omega(a_D,2-b_D)
\end{array}
\end{equation}
We denote by $M_\diamond$ the one-dimensional $B_H$-representation underlying $S_\diamond$, i.e. 
\[M_{F} := \lambda(2-a_F,2-b_F;m+u_F),\qquad M_{D} := \lambda(2-a_D,b_D;m+u_D).\]
We will also denote by $p_\diamond:V_\diamond|_{B_H}\to M_{\diamond}$ the projection corresponding to $E_{\diamond}\to S_{\diamond}$ at the level of vector spaces and by $p_\diamond^*:=p_{j_\diamond}^*$ the induced map on cohomology. 
\begin{defn}\label{Pairing:With:Projection}
We define the pairing
\begin{align*}
P_{\diamond}: H^{i_\diamond}(X_G,E_\diamond(-D_G))
\times H^{2-i_\diamond}(X_H,H_\diamond)
\to E\{m'_{\diamond}/2\},\;\;\;
(a,b)\mapsto \langle p_{j_{\diamond}}^*(a),b\rangle_{H},
\end{align*}
where $j_{\diamond}$ and $m'_{\diamond}$ are defined as in Table \eqref{Table:Images:p} and $\langle\cdot,\cdot\rangle$ is Serre duality for $X_{H}$. 
\end{defn}
\begin{rk}
Note that, at the specializations $a_F=a_D=0$ and
$b_F=b_D=c_0$, the corresponding indices $j_F$ and $j_D$
are $r_2$ and $r_1+1$, respectively. Both lie in $[0,m+2]$
and are therefore valid indices in \eqref{RestrictionLi}.
Consequently, the pairings $P_\diamond$ remain well-defined
at these specializations, even though they lie outside the
range considered in \cite{LPSZ1}.
\end{rk}
\begin{rk}
Alternatively, these pairings can be constructed using the Gysin pushforward
\[
\iota_*:H^q(X_H,[V|_{B_H}\otimes\lambda(-1,-1;0)])
\longrightarrow H^{q+1}(X_G,[V]),
\]
followed by Serre duality on $X_G$. As explained in \cite[\S 2.5]{LPSZ1}, the Gysin map is adjoint to the pullback under Serre duality, while the inclusions of the relevant summands are dual to the projections $p_{\diamond}^*$. Consequently, both constructions yield the same pairings $P_{\diamond}$.
\end{rk}
\subsubsection{Cohomology and automorphic representations}\label{subsec:scalar-normalization}\label{Remark:Dimension:One}
Let $X$ be a compactified modular curve, with Hodge bundle $\omega_X$
and cuspidal divisor $D_X$. Recall the standard identifications:
\begin{equation}\label{Modular:Forms:Classes:Cohom}
	M_k\simeq H^0(X,\omega_X^k),\quad S_k\simeq H^0(X,\omega_X^k(-D_X)),\quad S_k^\vee\simeq H^1(X,\omega_X^{2-k}).
\end{equation}
\begin{rk}
As explained in \eqref{Union:Modular:Curves}, at suitable levels, $X_{H}$ is isomorphic to the product of two compactified modular curves $X_1$ and $X_2$. In this case, we can interpret certain cohomology classes in terms of modular forms. In region $(F)$, we have that the pullback of $H_{F}$ under the isomorphism \eqref{Union:Modular:Curves} is isomorphic to $\omega_{X_1}^{a_F}\otimes \omega_{X_2}^{b_F}\otimes [\nu^{-2}]$. Then, the Künneth formula and the isomorphisms \eqref{Modular:Forms:Classes:Cohom} gives an isomorphism of Hecke modules of the form
\begin{equation}\label{Isom:H0:F}
H^0(X_H,H_F) \simeq (\mathcal{M}_{a_F}\boxtimes \mathcal{M}_{b_F})\otimes |\nu|_{\A}^{2},
\end{equation}
where $\mathcal{M}_{a_F}$ denotes modular forms of weight $a_F$. In region $(D)$, the Künneth formula gives
\begin{equation}\label{Kunneth:H1}H^1(X_{H},H_{D})\simeq H^0(X_1,\omega_{X_1}^{a_D})\otimes H^1(X_2,\omega_{X_2}^{2-b_D})\oplus H^1(X_1,\omega_{X_1}^{a_D})\otimes H^0(X_2,\omega_{X_2}^{2-b_D}).\end{equation}
By Serre duality, the first summand is naturally identified as a Hecke module with
$(\mathcal{M}_{a_D}\boxtimes S_{b_D}^{\vee})\otimes |\nu|_{\A}$, where $S_{b_D}^{\vee}$ denotes the dual of the space of cusp forms of weight $b_D$. 
\end{rk}
\begin{defn}
    In region $(F)$, let $f$ and $g$ be modular forms of weights
    $a_F$ and $b_F$, respectively considered as (possibly non unitary) functions on $[\GL_2]$. We denote by
    $f\boxtimes g\in H^0(X_H,H_F)$ the class corresponding to
    $(f,g)$ under the isomorphism \eqref{Isom:H0:F}. In region $(D)$, let $f$ be a modular form of weight $a_D$ and $g\in S_{b_D}^{\vee}$ a Serre-dual cohomology class, represented by the adelization of the complex conjugate of a cusp form of weight $b_D$. We denote by $f\boxtimes g\in H^1(X_H,H_D)$ the corresponding cohomology class.
\end{defn}
Let $\Pi$ be a unitary, globally generic automorphic representation of $G$. We define its arithmetic normalization by
\begin{equation}\label{Normalization:Pi}
\Pi_{\diamond}^{\mathrm{arith}}
:=\Pi\otimes|\nu|^{-\frac{m+u_\diamond}{2}},
\end{equation}
where we recall that $m=r_1+r_2$, $u_F=2r_2+2$, and $u_D=2r_1+4$.
\begin{theorem}[{\cite[Theorem 5.2]{LPSZ1}}]
\label{thm:Yoshida:Cohomology}
Assume that $\Pi$ is globally generic, cuspidal and non-CAP, of general or Yoshida type. At the fixed level $K_f$, the $\Pi_{\diamond,f}^{\mathrm{arith}}$-isotypical contribution to each of the Hecke modules
\[
H^{i_\diamond}(X_{G,\Q},E_\diamond(-D_G))\otimes\C
\quad\text{and}\quad
H^{i_\diamond}(X_{G,\Q},E_\diamond)\otimes\C
\]
is a direct summand isomorphic to $(\Pi_{\diamond,f}^{\mathrm{arith}})^{K_f}$. Here the acting algebra is the Hecke algebra $\mathcal H_{K_f}$ of $(G(\A_f),K_f)$; whenever $(\Pi_{\diamond,f}^{\mathrm{arith}})^{K_f}\ne0$, its multiplicity is one.
\end{theorem}
\begin{defn}
We denote the $(\Pi_{\diamond,f}^{\mathrm{arith}})^{K_f}$-isotypical summand of $H^{i_\diamond}(X_{G,\mathbb{Q}},E_\diamond(-D_G))\otimes\C$ by $H^{i_{\diamond}}(\Pi_f)_{\C}$.
\end{defn}
\subsection{Automorphic realization of $P_{\diamond}$}
In this subsection, following \cite{Harris:Partial}, we explain how to realize the coherent cohomology groups of Shimura varieties first as vector valued differential forms and subsequently as vector-valued automorphic forms. This interpretation will later allow us to express $P_{\diamond}$ in terms of integrals of automorphic forms.
\subsubsection{Coherent cohomology, differential forms and vector valued modular forms}\label{subsec:weight-zero}
Let $J\in\{G,H\}$, and write \[K_{G,\infty} = \R^{\times}\cdot\;U(2),\qquad K_{H,\infty} = \R^{\times}\cdot\,(SO(2)\times SO(2)).\] 
Let $Y_J\subset X_J$ denote the corresponding open Shimura variety of complex dimension $n$, with boundary divisor $D_J$. We denote by $\mathfrak p_J^+$ and $\mathfrak p_J^-$ the holomorphic and antiholomorphic tangent spaces at the standard base point $o$, respectively. Fix bases $\{\xi_1^+,\ldots,\xi_n^+\}$ and $\{\xi_1^-,\ldots,\xi_n^-\}$ of $\mathfrak p_J^+$ and $\mathfrak p_J^-$, corresponding to $\{\partial/\partial z_i\}_{i=1}^n$ and $\{\partial/\partial\overline z_i\}_{i=1}^n$, respectively, where $z_1,\ldots,z_n$ are local holomorphic coordinates at $o$. Recall that then $K_{Y_J} \simeq [\bigwedge^{n}(\mathfrak{p}^+_J)^{\vee}]$. For an automorphic vector bundle $[\mathcal E]$
on $X_J$, let
\[
\mathcal A_\star^{0,\bullet}(Y_J,[\mathcal E])
\]
denote the Dolbeault complex of $[\mathcal E]$-valued differential
forms satisfying the growth conditions of \cite{Harris:Partial},
where $\star=\mathrm{rd}$ (rapidly decreasing) or
$\star=\mathrm{si}$ (slowly increasing). These complexes compute
coherent cohomology:
\[
H^\bullet(X_J,[\mathcal E](-D_J))
\simeq H^\bullet_{\bar\partial}
(\mathcal A_{\mathrm{rd}}^{0,\bullet}(Y_J,[\mathcal E])),
\qquad
H^\bullet(X_J,[\mathcal E])
\simeq H^\bullet_{\bar\partial}
(\mathcal A_{\mathrm{si}}^{0,\bullet}(Y_J,[\mathcal E])).
\]
\begin{defn}\label{Definition:Harris:Cocycles}
Let $\mathcal{V}_{\mathcal{E}}^q :=
\mathcal E \otimes
\bigwedge\nolimits^q(\mathfrak p_J^-)^\vee$,
and let $\sigma_q:K_{J,\infty}\to\mathrm{GL}(\mathcal{V}_{\mathcal{E}}^q)$
denote the natural representation induced by the
$K_{J,\infty}$-actions on $\mathcal E$ and
$\mathfrak p_J^-$ by Remark \ref{Representation:Siegel:Compact}. We define
\[
C_{J,q,\mathcal{E}}
:=
\left\{
f\in C^\infty\left(
J(\Q)\backslash J(\A)/K_f,\mathcal{V}_{\mathcal{E}}^q
\right)
\;\middle|\;
f(gk)=\sigma_q(k)^{-1}f(g),
\quad
\forall g\in J(\A),\ k\in K_{J,\infty}
\right\}.
\]
\end{defn}
\begin{prop}[{\cite[(3.3.1)]{Harris:Partial}}]
There is a canonical isomorphism
\[
\underline{\mathrm{Lift}}_J:
\mathcal{A}^{0,q}(Y_J,[\mathcal E])
\xrightarrow{\sim} C_{J,q,\mathcal E},
\]
obtained by expressing differential forms in the frames induced by the action of $J(\A)$ on the fibre at a fixed base point.
\end{prop}
Let $[\mathcal E'_{m'}]=[\mathcal E]^\vee\otimes K_{X_J}(D_J)\otimes C_{m'}$ for some $m'\in 2\mathbb{Z}$ and consider
\[
\alpha\in\mathcal A_{\mathrm{rd}}^{0,n-q}(Y_J,[\mathcal E'_{m'}]),
\qquad
\beta\in\mathcal A_{\mathrm{si}}^{0,q}(Y_J,[\mathcal E]).
\]
Write $f=\underline{\mathrm{Lift}}_J(\alpha)$ and $g=\underline{\mathrm{Lift}}_J(\beta)$.
On the compactification the coefficient contraction is
\begin{equation}\label{Pairing:Twisted}
[\mathcal E'_{m'}](-D_J)\otimes[\mathcal E]\longrightarrow K_{X_J}\otimes C_{m'}.
\end{equation}
Indeed, $[\mathcal E'_{m'}](-D_J)=[\mathcal E]^\vee\otimes K_{X_J}\otimes C_{m'}$. The rapidly decreasing representative supplies this subcanonical factor, while the slowly increasing representative belongs to the canonical extension. On $Y_J$ the boundary twist disappears, and contraction followed by exterior product induces the bilinear map

\[
C_{J,n-q,\mathcal E'_{m'}}\otimes C_{J,q,\mathcal E}
\longrightarrow C_{J,n,\bigwedge^{n}(\mathfrak{p}^+_J)^{\vee}\otimes\lambda(0,0;m')},
\qquad
f\otimes g\longmapsto [f\wedge g].
\]
This map admits a very precise description. Consider the $K_{J,\infty}$-equivariant map of representations
\[
\begin{aligned}
m:
\left(\mathcal E'_{m'}\otimes\bigwedge^{n-q}(\mathfrak p_J^-)^\vee\right)
\otimes
\left(\mathcal E\otimes
\bigwedge^q(\mathfrak p_J^-)^\vee\right)
\to
\bigwedge^{2n}(\mathfrak p_J^+\oplus\mathfrak p_J^-)^{\vee}\otimes \lambda(0,0,m'),
\end{aligned}
\]
using $\mathcal E'_{m'}=\mathcal E^\vee\otimes\bigwedge^n(\mathfrak p_J^+)^\vee\otimes\lambda(0,0;m')$ in the compact realization. We contract the $\mathcal E$-factors and place the holomorphic exterior factor first, followed by the antiholomorphic factors in the order $f,g$. Then
\begin{equation}\label{Wedge:Contraction}
[f\wedge g](h):=m\bigl(f(h)\otimes g(h)\bigr),
\quad h\in J(\mathbb A),
\end{equation}  
In particular, evaluation on any element of
$v\in \bigwedge^{2n}(\mathfrak p_J^+\oplus\mathfrak p_J^-)$ gives
\[
[f\wedge g](h)(v)\in\lambda(0,0;m').
\]
\begin{rk}
    Since $Z_J(\R)\subset K_{J,\infty}$ acts trivially on $\mathfrak p_J^\pm$, and its actions on $\mathcal E$ and $\mathcal E^\vee$ cancel under contraction, the $K_{J,\infty}$-equivariance of Harris's lift (and hence of $m$) implies that $[f\wedge g]$ transforms under $Z_J(\R)$ by $\lambda(0,0,-m')$. The change from $m'$ to $-m'$ follows from Definition \ref{Definition:Harris:Cocycles}.
\end{rk}
\begin{theorem}\label{Harris:Main:Theorem}
For classes $[\alpha]\in H^{n-q}(X_J,[\mathcal E'_{m'}](-D_J))$ and $[\beta]\in H^q(X_J,[\mathcal E])$, choose rapidly decreasing and slowly increasing Dolbeault representatives $\alpha$ and $\beta$, respectively, and put $f:=\underline{\mathrm{Lift}}_J(\alpha)$ and $g:=\underline{\mathrm{Lift}}_J(\beta)$. With $[\alpha]$ first in the cup product, Serre duality is given, in the $e_{m'/2}$-coordinate, by
\[
\langle[\alpha],[\beta]\rangle_{X_J}
\overset{.}{=}
\int_{Z_J(\R)J(\Q)\setminus J(\A)}|\nu(h)|_{\A}^{m'/2}[f(h)\wedge g(h)](1_J)\,dh,
\]
where
\[
1_J=(-1)^{n(n-1)/2}(2/i)^n\xi_1^+\wedge\cdots\wedge\xi_n^+\wedge\xi_1^-\wedge\cdots\wedge\xi_n^-.
\]
\end{theorem}
\begin{proof}
The identity $[\mathcal E'_{m'}](-D_J)=[\mathcal E]^\vee\otimes K_{X_J}\otimes C_{m'}$ identifies the cup product with a class in $H^n(X_J,K_{X_J}\otimes C_{m'})$. Its Dolbeault representative is the contraction of $\alpha\wedge\beta$; rapid decrease of the first factor makes its product with the slowly increasing second factor integrable and makes integration vanish on coboundaries. Harris's realization of Serre trace \cite[(3.7.1)--(3.8.4)]{Harris:Partial} therefore applies to this product. In the adelic frame the coefficient has central character $\lambda(0,0;m')$, so the lifted function transforms by its inverse. Passing to the $e_{m'/2}$-coordinate of Remark~\ref{Remark:Trivialization} multiplies it by $|\nu(h)|_{\A}^{m'/2}$. This accounts for the positive exponent. Evaluation on $1_J$ uses the complex orientation; the fixed comparison between the resulting volume form and $dh$ is included in the proportionality constant.
\end{proof}
\begin{comment}
The normalization of $1_J$ follows from the coordinate bases fixed above. Writing $z_j=x_j+iy_j$, we have $\partial_{z_j}=\tfrac12(\partial_{x_j}-i\partial_{y_j})$ and $\partial_{\bar z_j}=\tfrac12(\partial_{x_j}+i\partial_{y_j})$, so
\[
\xi_1^+\wedge\cdots\wedge\xi_n^+\wedge\xi_1^-\wedge\cdots\wedge\xi_n^-
=(-1)^{n(n-1)/2}(i/2)^n
\partial_{x_1}\wedge\partial_{y_1}\wedge\cdots\wedge\partial_{x_n}\wedge\partial_{y_n}.
\]
Thus $1_J$ is the real volume multivector for the complex orientation. This fixes the sign and the coordinate factor in its definition; the comparison with the chosen Haar measure is included in the proportionality in the Theorem.\\
\end{comment}
Finally, let $\iota:Y_H\to Y_G$ be the morphism of \S\ref{subsubsec:shimura-varieties}.
For $\alpha\in\mathcal A^{0,q}(Y_G,[\mathcal E])$, we have 
\begin{equation*}
\underline{\mathrm{Lift}}_H(\iota^*\alpha)(h)(w_1,\ldots,w_q)
=
\underline{\mathrm{Lift}}_G(\alpha)(\iota(h))
\bigl(d\iota(w_1),\ldots,d\iota(w_q)\bigr),
\end{equation*}
for $h\in H(\A)$ and $w_i\in\mathfrak p_H^-$. Furthermore, $\underline{\mathrm{Lift}}$ is also compatible with
exterior products and projections, i.e. the following diagram is commutative:
\begin{equation}\label{Pullback:And:Lift}
\begin{tikzcd}[column sep=9em, row sep=4em]
[H]
    \arrow[r, "{\underline{\mathrm{Lift}}_G(\alpha)|_H}"]
    \arrow[dr, swap,
        "{\underline{\mathrm{Lift}}_H(p_\diamond^*(\alpha))}"]
&
V_\diamond \otimes
\bigwedge^{i_\diamond}(\mathfrak p_G^-)^\vee
    \arrow[d, "p_\diamond\otimes\bigwedge^{i_\diamond}\iota^*"]
\\
&
M_\diamond \otimes
\bigwedge^{i_\diamond}(\mathfrak p_H^-)^\vee
\end{tikzcd}
\end{equation}
\begin{lemma}\label{First:Formula:Pairing}
    Let $[\alpha]\in H^{i_\diamond}(X_G,E_\diamond(-D_G))$ and $[\beta]\in H^{2-i_\diamond}(X_H,H_\diamond)$, where $\alpha$ and $\beta$ are rapidly decreasing and slowly increasing Dolbeault representatives. We have that 
    \[P_{\diamond}([\alpha],[\beta]) \overset{.}{=} 
    \int_{Z_{H}(\R)H(\Q)\setminus H(\A)}|\nu(h)|^{m}_{\A}[\underline{\mathrm{Lift}}_{H}(p_{\diamond}^*(\alpha))(h)\wedge \underline{\mathrm{Lift}}_{H}(\beta)(h)](1_H)\,dh,\]
    where $\underline{\mathrm{Lift}}_{H}(p_{\diamond}^*(\alpha))(h):H(\Q)\setminus H(\A)/K_{H,f}\to M_{\diamond}\otimes \bigwedge^{i_{\diamond}}(\mathfrak{p}_H^-)^{\vee}$.
\end{lemma}
\begin{proof}
Apply Theorem~\ref{Harris:Main:Theorem} with $n=2$, $q=2-i_\diamond$, $\mathcal E=W_\diamond$ and $m'=2m$. The identity \eqref{Condition:Serre:Dual}, with the values in \eqref{Table:Images:p}, gives $[\mathcal E'_{2m}]=S_\diamond$, and \eqref{Pullback:Vector:Bundles} places $p_\diamond^*[\alpha]$ in $H^{i_\diamond}(X_H,S_\diamond(-D_H))$. Formula \eqref{Pullback:And:Lift} identifies its representative with the one in the integral.
\end{proof}
Let $f_i\in C_{H,i,\mathcal{E}}$ and $g_j\in C_{H,j,\mathcal{E}'}$,
with $i+j=2$. By the definition of the exterior product, we have that $[f_i\wedge g_j](h)
(\xi_1^-\wedge\xi_2^-\wedge\xi_1^+\wedge\xi_2^+)$ equals
\begin{align}\label{Wedge:Evaluated}
\begin{cases}
m\left(
f_2(h)(\xi_1^-,\xi_2^-)\otimes g_0(h)
\right)(\xi_1^+,\xi_2^+),
& (i,j)=(2,0),\\[6pt]
\sum_{\sigma\in S_2}\operatorname{sgn}(\sigma)\,
m\left(
f_1(h)(\xi_{\sigma(1)}^-)
\otimes
g_1(h)(\xi_{\sigma(2)}^-)
\right)(\xi_1^+,\xi_2^+)
& (i,j)=(1,1),\\[6pt]
m\left(
f_0(h)\otimes g_2(h)(\xi_1^-,\xi_2^-)
\right)(\xi_1^+,\xi_2^+),
& (i,j)=(0,2).
\end{cases}
\end{align}
where $S_2$ is the symmetric group on two elements.
\begin{lemma}\label{Explicit:Expresion:Wedge:Lift}
    Let $[\alpha]$ and $[\beta]$ be as in Lemma~\ref{First:Formula:Pairing}, and suppose that $[\beta]$ corresponds to a pure tensor $f\boxtimes g$ in $\mathcal M_{a_F}\boxtimes\mathcal M_{b_F}$ in region $(F)$, or in $\mathcal M_{a_D}\boxtimes S_{b_D}^{\vee}$ in region $(D)$. Then there exists an automorphic form $\mathcal{F}$ on $G$ such that, in the $e_m$-coordinate fixed in Remark~\ref{Remark:Trivialization},
    \begin{align*}
        [\underline{\mathrm{Lift}}_{H}(p(\alpha))(h)\wedge \underline{\mathrm{Lift}}_{H}(\beta)(h)](1_H) &= |\nu(h)|_{\A}^{i_{\diamond}}\mathcal{F}(\iota(h))f(h_1)g(h_2).
    \end{align*}
\end{lemma}
\begin{proof}
The coefficient contraction in \eqref{Contraction:After:Projection} is induced by $M_\diamond\otimes W_\diamond\simeq\lambda(2,2;2m)$. Indeed, the sums of the central weights are $m+u_F+a_F+b_F-4=2m$ and $m+u_D+a_D-b_D-2=2m$, respectively. Under the compact realization of Remark~\ref{Representation:Siegel:Compact}, this is $\bigwedge^2(\mathfrak p_H^+)^\vee\otimes\lambda(0,0;2m)$. Including the antiholomorphic factors, the contraction is
\[
\begin{aligned}
m:\bigl(M_\diamond\otimes\bigwedge^{i_\diamond}(\mathfrak p_H^-)^\vee\bigr)\otimes\bigl(W_\diamond\otimes\bigwedge^{2-i_\diamond}(\mathfrak p_H^-)^\vee\bigr)
\longrightarrow\bigwedge^4(\mathfrak p_H^+\oplus\mathfrak p_H^-)^\vee\otimes\lambda(0,0;2m).
\end{aligned}
\]
The two top exterior factors have compact weights $(-2,-2)$ and $(2,2)$, so evaluation on $1_H$ leaves precisely the central line $\lambda(0,0;2m)$. Its chosen fibre generator corresponds to $e_m$ in the target $E\{m\}$.

We first consider region $(F)$. By Harris's lift,
\[
\underline{\mathrm{Lift}}_H(\beta):
H(\Q)\backslash H(\A)/K_{f,H}
\longrightarrow
W_F.
\]
The compatibility of Harris's lift with tensor products and the identification $\underline{\mathrm{Lift}}_{\GL_2}(f)\overset{.}{=}f$ allow us to choose a basis $e_F$ of this one-dimensional target such that
\[
\underline{\mathrm{Lift}}_H(\beta)(h)
=|\nu(h)|^2f(h_1)g(h_2)e_F,
\]
where the similitude factor arises due to the isomorphism given in \eqref{Isom:H0:F}. On the other hand, by Lemma \ref{First:Formula:Pairing}, we have
\[
\underline{\mathrm{Lift}}_H(p_F^*(\alpha)):
H(\Q)\backslash H(\A)/K_{f,H}
\longrightarrow
M_F
\otimes\bigwedge^2(\mathfrak p_H^-)^\vee.
\]
This target is one-dimensional, and the contraction just described is nonzero on it. We choose
\[
\widetilde e_F\in
M_F
\otimes\bigwedge^2(\mathfrak p_H^-)^\vee,\;s.t.\;[m(\widetilde e_F\otimes e_F)](1_H)=e_m.
\]
With this normalization, we can write
\[
\underline{\mathrm{Lift}}_H(p_F^*(\alpha))(h)
=\mathcal{F}(\iota(h))\widetilde e_F,
\]
where, by \eqref{Pullback:And:Lift}, the function $\mathcal{F}$ is a scalar-valued automorphic form on $G$. Hence, by \eqref{Wedge:Contraction}
\[
\begin{aligned}
&[\underline{\mathrm{Lift}}_H(p_F^*(\alpha))(h)
\wedge\underline{\mathrm{Lift}}_H(\beta)(h)](1_H)=|\nu(h)|_{\A}^{2}\mathcal{F}(\iota(h))f(h_1)g(h_2)e_m.
\end{aligned}
\]
Taking the $e_m$-coordinate in $\lambda(0,0;2m)$ gives the asserted scalar formula.
We now consider region $(D)$. By Harris's lift, we have
\[
\underline{\mathrm{Lift}}_H(\beta):
H(\Q)\backslash H(\A)/K_{H,f}
\longrightarrow
W_D\otimes(\mathfrak p_H^-)^\vee.
\]
Fix a basis $e_D$ of $W_D$, and let $\theta_1,\theta_2$ be dual to $\xi_1^-,\xi_2^-$, in the order of the two modular curves.
We can write
\[
\underline{\mathrm{Lift}}_H(\beta)(h)
=e_D\otimes\bigl(B_1(h)\theta_1+B_2(h)\theta_2\bigr),
\]
where $B_i:H(\Q)\backslash H(\A)/K_{H,f}\to\C$. Recall that, by \eqref{Kunneth:H1}, we are considering a class
\[
\beta=f\boxtimes g\in
H^0(X_1,\omega_{X_1}^{a_D})
\otimes H^1(X_2,\omega_{X_2}^{2-b_D}).
\]
Choosing local coordinates $(z_1,z_2)$ on $X_1\times X_2$,
we may choose a Dolbeault representative of $\beta$ of the form
$f(z_1)g(z_2)\,d\overline z_2$.
Since this representative has no $d\overline z_1$-component, by Harris's lift construction we obtain $B_1(h)=0$. Thus,
\[
\underline{\mathrm{Lift}}_H(\beta)(h)
=|\nu|_{\A}f(h_1)g(h_2)e_D\otimes\theta_2.
\]
Choose $\widetilde e_D\in M_D\otimes\mathbb C\theta_1$ so that $[m(\widetilde e_D\otimes(e_D\otimes\theta_2))](1_H)=e_m$. Write the $\theta_1$-component of $\underline{\mathrm{Lift}}_H(p_D^*\alpha)(h)$ as $\mathcal F(\iota(h))\widetilde e_D$. Its $\theta_2$-component contributes zero, since $\theta_2\wedge\theta_2=0$. Thus \eqref{Wedge:Evaluated}, with the Siegel factor first, gives $|\nu(h)|_{\A}\mathcal F(\iota(h))f(h_1)g(h_2)e_m$.
\end{proof}
\begin{rk}\label{How:Is:F}
    Note that $\mathcal F$ is obtained from a $K_{G,\infty}$-equivariant map in
    \[
    \mathrm{Hom}_{K_{G,\infty}}
    \left(
    V_\diamond^\vee\otimes
    \bigwedge^{i_\diamond}(\mathfrak p_G^-),
    \mathcal C^\infty(G(\Q)\backslash G(\A)/K_{G,f})
    \right).
    \]
    This allows us to conclude that $\mathcal F$ is right $K_{G,f}$-invariant and, by $K_{G,\infty}$-equivariance, that the $K_{G,\infty}$-representation generated by $\mathcal F$ is a quotient of a subrepresentation of $V_\diamond^\vee\otimes
    \bigwedge^{i_\diamond}(\mathfrak p_G^-)$. In particular, $\mathcal F$ is
    $K_{G,\infty}$-finite.
\end{rk}
\subsection{Characterizing the automorphic form $\mathcal{F}$}
By combining Lemmas \ref{First:Formula:Pairing} and \ref{Explicit:Expresion:Wedge:Lift}, we can express these pairings as integrals of automorphic functions. However, this description is not yet sufficiently precise for our purposes. In particular, we need to identify the specific automorphic form $\mathcal{F}$ on $G$ appearing in Lemma \ref{Explicit:Expresion:Wedge:Lift} refining Remark \ref{How:Is:F}, i.e. we have to determine its minimal $K_{G,\infty}$-type and the corresponding vector within that $K_{G,\infty}$-type. This identification will play a crucial role in establishing a non-vanishing result later on. To this end, we will use \cite[Theorem 5.6(a) and Lemma 6.5]{LPSZ1} and \cite{Su24}, which provide a way to identify the relevant automorphic forms in terms of the data encoded in the cocycles $C_{J,q,\mathcal{E}}$ introduced by Harris.
\subsubsection{Factorization of Harris's lift}
First, as in Remark \ref{How:Is:F}, by composing with linear functionals on the target representation, we may identify elements of $C_{G,q,\mathcal E}$ with $K_{G,\infty}$-equivariant maps in
    \begin{equation}\label{Changing:The:Hom}\mathrm{Hom}_{K_{G,\infty}}(\mathcal{E}^{\vee}\otimes \bigwedge^q(\mathfrak{p}^{-}_G),\mathcal{C}^{\infty}(G(\Q)\setminus G(\A)/K_{G,f})).\end{equation}
Given a cohomology class $\eta\in H^{i_{\diamond}}(X_G,E_{\diamond}(-D_G))$, we will denote by $\widetilde{\mathcal{F}}_{\eta}$ its corresponding homomorphism after applying the map $\underline{\mathrm{Lift}}$ and the above identification.
\begin{defn}\label{Definition:Constant:Projection}
    Let $
p_\diamond: V_\diamond\longrightarrow
 M_\diamond$ be the projection specified in Table \eqref{Table:Images:p}.
Fix a nonzero vector $e_\diamond\in M_\diamond$.
Since $M_\diamond$ is one-dimensional, there is a
unique linear functional
$Q_\diamond\in V_\diamond^\vee$ such that
\[
p_\diamond(v)=Q_\diamond(v)e_\diamond,
\qquad v\in V_\diamond.
\]
Note that $Q_{\diamond}$ depends on the choice of $e_{\diamond}$, although we suppress this dependence from the notation.
\end{defn}
\begin{lemma}\label{Explicit:Expression:F}
    With $Q_\diamond$ as in Definition~\ref{Definition:Constant:Projection}, we have
    \begin{align*}
    \mathcal{F}(h) =\begin{cases} \widetilde{F}_{\alpha}(Q_{\diamond}\otimes d\iota(\xi_1^-\wedge\xi_2^-))(\iota(h))&\textit{region }(F),\\
     \widetilde{F}_{\alpha}(Q_{\diamond}\otimes d\iota(\xi_1^-))(\iota(h))&\textit{region }(D).\end{cases}
    \end{align*}
\end{lemma}
\begin{proof}
    Recall that, as in the proof of Lemma \ref{Explicit:Expresion:Wedge:Lift}, we have
    \[\mathcal{F}(h)e_{\diamond} = \underline{\mathrm{Lift}}_{H}(p^*_{\diamond}(\alpha))(h)(\xi_1^-\wedge\xi_2^-)\]
    By commutativity of the diagram \eqref{Pullback:And:Lift}, this equals
    \[p_{\diamond}(\underline{\mathrm{Lift}}_{G}(\alpha)(\iota(h))(d\iota(\xi_1^-\wedge\xi_2^-))).\]
    Hence, by Definition \ref{Definition:Constant:Projection}, one has that 
    \[\mathcal{F}(h)e_{\diamond} = Q_{\diamond}(\underline{\mathrm{Lift}}_{G}(\alpha)(\iota(h))(d\iota(\xi_1^-\wedge\xi_2^-)))e_{\diamond}.\]
    The argument for region $(D)$ is analogous. 
\end{proof}
In what follows, let
\[
\tau_{i_{\diamond}}^{\mathrm{arith}}
:= \tau_{i_{\diamond}}\boxtimes |\nu|^{-\frac{m+u_{\diamond}}{2}}
\]
denote the minimal $K_{\infty}$-type of the archimedean component of
$\Pi^{(\diamond),\mathrm{arith}}$, as defined in \eqref{Normalization:Pi}.
\begin{theorem}[\cite{Su24}]\label{thm:harris:su}
There are isomorphisms of $\mathcal H_{K_f}$-modules of the form
\begin{align*}
    H^{i_{\diamond}}(\Pi_f)_{\C}\to \mathrm{Hom}_{K_{G,\infty}}\left(\tau_{i_{\diamond}}^{\mathrm{arith}},\Pi^{(\diamond),\mathrm{arith},K_f}\right),\quad\eta\mapsto F_{\eta}.
\end{align*}
\end{theorem}
Assume that $\Pi$ is a unitary cuspidal automorphic representation
of general or Yoshida type, globally generic and non-CAP. We have two realizations of the classes
$\eta\in H^{i_\diamond}(\Pi_f)_{\C}$. On the one hand, combining Theorem \ref{thm:Yoshida:Cohomology}
with \cite[Proposition~3.6, Proposition~4.5 (b),
Theorem~4.6.2 \& Lemma~5.2.3]{Harris:Partial},
we may represent each $\eta$ by a unique harmonic cuspidal
differential form. Applying Harris's lift to this representative and using the $G(\A_f)$-equivariance given by \cite[Corollary 3.4]{Harris:Partial}, we have a $K_{G,\infty}$-equivariant map
\[
\widetilde F_\eta:
\mathcal V_{i_\diamond}
\longrightarrow
\Pi^{(\diamond),\mathrm{arith},K_f},
\]
where
\[
\mathcal V_2
:=V_F^\vee\otimes\bigwedge^2\mathfrak p_G^-,
\qquad
\mathcal V_1
:=V_D^\vee\otimes\mathfrak p_G^-.
\]
On the other hand, the preceding theorem associates to $\eta$ a
$K_{G,\infty}$-equivariant homomorphism
\[
F_\eta:
\tau_{i_\diamond}^{\mathrm{arith}}
\longrightarrow
\Pi^{(\diamond),\mathrm{arith},K_f}.
\]
To compare these realizations, we will construct a
$K_{G,\infty}$-equivariant surjection $C_\diamond$, unique up to
scalars, and a nonzero scalar $s_\diamond$, independent of $\eta$,
such that the following diagram commutes:
\[
\begin{tikzcd}[column sep=large, row sep=large]
\mathcal{V}_{i_{\diamond}} \arrow[rr, "\widetilde F_\eta"]
    \arrow[dr, dotted, "C_{\diamond}"']
&& \Pi^{(\diamond),\mathrm{arith},K_f} \\
& \tau_{i_{\diamond}}^{\mathrm{arith}} \arrow[ur, "F_\eta"'] &
\end{tikzcd}
\]
More precisely, an analysis of the irreducible constituents of
$\mathcal V_{i_\diamond}$ will show that
$\tau_{i_\diamond}^{\mathrm{arith}}$ occurs with multiplicity one
and that $\widetilde F_\eta$ vanishes on all other isotypical
components. This determines the $K_{G,\infty}$-equivariant
surjection $C_\diamond$, unique up to scalars, through which
$\widetilde F_\eta$ factors. Once $C_\diamond$ is fixed,
Hecke equivariance and Schur's Lemma, applied to the relevant
irreducible Hecke module, will give
\[
\widetilde F_\eta
=
s_\diamond F_\eta\circ C_\diamond,
\]
where $s_\diamond\in\C^\times$ is independent of $\eta$.\\

The centre of $K_{G,\infty}$ acts by the same character on $\mathcal V_{i_\diamond}$ and on the arithmetic realization defined in \eqref{Normalization:Pi}. It therefore suffices to decompose $\mathcal V_{i_\diamond}$ into irreducible $U(2)$-representations and compare these with the $U(2)$-types of $\Pi_\infty$. This identifies the components used to define the projections $C_\diamond$.

\begin{lemma}\label{lem:twisted-domains}
Let $U$ be the standard representation of $U(2)$ and $n:=m+2$, where $m = r_1+r_2$. We have that 
\begin{align*}
\mathcal{V}_2|_{U(2)}\simeq \bigoplus_{\ell=0}^{2}\Sym^{n+2-2\ell}(U)\otimes \mathrm{det}^{-r_1-3+\ell},\quad
\mathcal{V}_1|_{U(2)}\simeq\bigoplus_{\ell=0}^{2}\mathrm{Sym}^{n+2-2\ell}(U)\otimes \mathrm{det}^{-r_2-1+\ell}.
\end{align*}
\end{lemma}
\begin{proof}
By Remark
\ref{Representation:Siegel:Compact}, the representation
$W_G(A,B;C)$ of $M_{\mathrm{Si}}$ corresponds to a representation
of $K_{G,\infty}$ with highest weight $(-B,-A;C)$. Consequently, we obtain
\[
V_F^\vee|_{U(2)}
\simeq \Sym^n U\otimes\mathrm{det}^{-r_1},
\qquad
V_D^\vee|_{U(2)}
\simeq \Sym^n U\otimes\mathrm{det}^{1-r_2}.
\]
By the proof of \cite[Lemma 2.14]{Horawa:Prasanna}, one has that $\mathfrak{p}_G^-\simeq \mathrm{Sym}^2(U)\otimes \mathrm{det}^{-2}$ as a representation of $U(2)$, therefore one can deduce that 
$\bigwedge^2\mathfrak{p}^-_{G}\simeq \mathrm{Sym}^2(U)\otimes \mathrm{det}^{-3}$, and hence
\begin{equation}\label{First:Description:V:cal}\mathcal{V}_2|_{U(2)}\simeq \mathrm{Sym}^2(U)\otimes\mathrm{Sym}^n(U)\otimes \mathrm{det}^{-r_1-3},\quad \mathcal{V}_1|_{U(2)}\simeq \mathrm{Sym}^2(U)\otimes \mathrm{Sym}^n(U)\otimes \mathrm{det}^{-1-r_2}.\end{equation}
Using the Clebsch--Gordan formula, one gets the decomposition asserted in the Lemma
\end{proof}
\begin{rk}\label{Coincidence:Minimal:K:types}
   By the preceding Lemma, the minimal $U(2)$-types of Definition \ref{Minimal:U:Types},
\[
\tau_2=\Sym^{n+2}U\otimes\mathrm{det}^{-r_1-3},
\qquad
\tau_1=\Sym^{n+2}U\otimes\mathrm{det}^{-r_2-1},
\]
occur with multiplicity one in $\mathcal V_{i_\diamond}|_{U(2)}$.
Moreover, by the Clebsch--Gordan decomposition and Schur's Lemma,
the projection
\[
C_\diamond:\mathcal V_{i_\diamond}|_{U(2)}
\longrightarrow \tau_{i_\diamond}
\]
is given, up to a nonzero scalar, by the multiplication map
\[
\begin{aligned}
\Sym^2 U\otimes\Sym^n U &\longrightarrow \Sym^{n+2}U,\\
P\otimes Q &\longmapsto PQ,
\end{aligned}
\]
under the description in \eqref{First:Description:V:cal}.
Here multiplication is taken in the polynomial realization
of the symmetric powers of $U$ (see \eqref{Poly:Model}). We conclude that then $C_{\diamond}$ extends trivially to a map of the form
\[C_{\diamond}:\mathcal{V}_{i_{\diamond}}\to \tau_{i_{\diamond}}^{\mathrm{arith}}.\]
\end{rk}
We now turn to the study of the $K_{G,\infty}$-equivariant maps
$\widetilde{F}_{\alpha}:\mathcal{V}_{i_\diamond}\to
\Pi^{(\diamond),\mathrm{arith}}$.
Since $K_{G,\infty}< G(\R)^+$, we have to consider the decomposition
of $\Pi^{(\diamond),\mathrm{arith}}_\infty|_{G(\R)^+}$
derived from \S\ref{Archimedean:Component:Siegel}.
We denote their irreducible constituents by
\[
\Pi_{\infty,i_{\diamond}}^{(\diamond),\mathrm{arith}}
:=\Pi_{\infty,i_{\diamond}}
\otimes|\nu|_\infty^{-\frac{m+u_{\diamond}}{2}},
\]
where $\Pi_{\infty,i_{\diamond}}$ was defined in \S\ref{Archimedean:Component:Siegel}.
\begin{lemma}\label{lem:twisted-unique-types}
For every $0\leq q\leq3$, we have
\[
\dim\Hom_{K_{G,\infty}}
\left(
\bigwedge^q\mathfrak p_G^-\otimes V_\diamond^\vee,
\Pi_{\infty,i_{\diamond}}^{(\diamond),\mathrm{arith}}
\right)
=
\begin{cases}
1,&q=i_{\diamond},\\
0,&q\neq i_{\diamond},
\end{cases}
\]
For $q=i_\diamond$, the source identifies with $\mathcal V_{i_\diamond}$ by the canonical interchange of tensor factors, and the unique contributing $U(2)$-type is $\tau_{i_\diamond}$.
\end{lemma}

\begin{proof}
Recall that $V_\diamond=L_{i_\diamond}\otimes C_{u_\diamond}$, hence
\begin{equation}\label{Removing:Real:Twist}
\begin{aligned}
\Hom_{K_{G,\infty}}
\left(
\bigwedge^q\mathfrak p_G^-
\otimes V_{\diamond}^\vee,
\Pi_{\infty,i_{\diamond}}^{(\diamond),\mathrm{arith}}
\right)\simeq
\Hom_{K_{G,\infty}}
\left(
\bigwedge^q\mathfrak p_G^-
\otimes L_{i_{\diamond}}^\vee,
\Pi_{\infty,i_{\diamond}}\otimes|\nu|_\infty^{-\frac{m}{2}}
\right).
\end{aligned}
\end{equation}
and analogously for $(D)$. The central characters on the right-hand side agree, so it remains to compare their $U(2)$-types. We proceed as in the proof of \cite[Lemma 6.5]{LPSZ1}.
Recall that the highest weights of $L_2^\vee|_{U(2)}$ and $L_1^\vee|_{U(2)}$ are $(r_2+2,-r_1)$ and $(r_1+3,1-r_2)$ respectively. Hence, by the Clebsch--Gordan decomposition, if $R$ is an
irreducible $U(2)$-representation of highest weight
$(A,B)$, with $A-B\geq 2$, the highest weights of
$\bigwedge^q\mathfrak p_G^-\otimes R$ are
\[
\begin{array}{c|c}
q & \text{highest weights}\\ \hline
0 & (A,B)\\
1 & (A-\ell,B-2+\ell),\quad 0\leq\ell\leq2\\
2 & (A-1-\ell,B-3+\ell),\quad 0\leq\ell\leq2\\
3 & (A-3,B-3).
\end{array}
\]
Each constituent occurs with multiplicity one. For region $(F)$, \cite[Lemma 6.5]{LPSZ1} shows
that the $U(2)$-types of $\Pi_{\infty,2}$ lie in the cone
\[
(r_2+1,-r_1-3)
+M(1,0)+N(-1,-1),
\qquad M,N\geq0.
\]
Taking $(A,B)=(r_2+2,-r_1)$, the only common type
occurs for $q=2$ and $\ell=0$, and is precisely
$\tau_2$. Similarly, the $U(2)$-types of $\Pi_{\infty,1}$ lie in the cone
\[
(r_1+3,-r_2-1)
+M(0,-1)+N(1,1),
\qquad M,N\geq0.
\]
Taking $(A,B)=(r_1+3,1-r_2)$, the only common type
occurs for $q=1$ and $\ell=0$, giving $\tau_1$. Since these types occur with multiplicity one in both the source and the target, the result follows.
\end{proof}
\begin{prop}\label{prop:twisted-factorization}
One has that 
\[\widetilde F_\eta=s_{\diamond} F_\eta\circ C_{\diamond},\]
for some $s_F,s_D\in\C^\times$.
\end{prop}
\begin{proof}
By Remark \ref{Coincidence:Minimal:K:types} and
Lemma \ref{lem:twisted-unique-types}, $\widetilde F_\eta$
factors through $C_\diamond$. Since $C_\diamond$ is surjective,
there is a unique
\[
B_\eta\in\Hom_{K_{G,\infty}}
\left(
\tau_{i_\diamond}^{\mathrm{arith}},
\Pi^{(\diamond),\mathrm{arith},K_f}
\right)
\quad\text{such that}\quad
\widetilde F_\eta=B_\eta\circ C_\diamond.
\]
Let $\mathcal H_{K_f}$ denote the Hecke algebra of
$(G(\A_f),K_f)$. By Lemma \ref{lem:twisted-unique-types},
the above Hom-space is isomorphic, as an $\mathcal H_{K_f}$-module, to $\Pi_f^{K_f}\otimes|\nu|_f^{-(m+u_\diamond)/2}$. The map $F_\bullet:\eta\mapsto F_\eta$ is a Hecke-equivariant
isomorphism by Theorem \ref{thm:harris:su}.
Hecke equivariance of Harris's lift and uniqueness of the
factorization imply that $B_\bullet:\eta\mapsto B_\eta$
is also Hecke-equivariant. Hence
\[
B_\bullet\circ F_\bullet^{-1}
\in\mathrm{End}_{\mathcal H_{K_f}}(\Pi_f^{K_f}\otimes|\nu|_f^{-(m+u_\diamond)/2}).
\]
Since $\Pi_f$ is irreducible and admissible, $\Pi_f^{K_f}$ is
an irreducible finite-dimensional $\mathcal H_{K_f}$-module.
Schur's Lemma therefore gives
$B_\bullet=s_\diamond F_\bullet$.
\end{proof}
\subsubsection{The vectors selected by the factorization}
Lemma \ref{Explicit:Expression:F} and Proposition \ref{prop:twisted-factorization} allow us to express
\begin{align*}
    \mathcal{F}(h) =\begin{cases} s_F F_{\alpha}(C_{F}(Q_{F}\otimes d\iota(\xi_1^-\wedge\xi_2^-)))(\iota(h))&\textit{region }(F),\\
     s_D F_{\alpha}(C_{D}(Q_{D}\otimes d\iota(\xi_1^-)))(\iota(h))&\textit{region }(D).\end{cases}
    \end{align*}
To determine precisely which automorphic forms arise in this way, it remains to identify the vectors
\[C_{F}(Q_{F}\otimes d\iota(\xi_1^-\wedge\xi_2^-)),\qquad C_{D}(Q_{D}\otimes d\iota(\xi_1^-)),\]
in $\tau_{i_{\diamond}}^{\mathrm{arith}}$.\\ 

Before stating the next Lemma, choose a basis $\{x,y\}$ of $U$ of weight vectors for the diagonal torus. We identify
\begin{equation}\label{Poly:Model}
\mathrm{Sym}^n(U)\simeq\bigoplus_{j=0}^n\C x^{n-j}y^j
\end{equation}
with the space of homogeneous polynomials of degree $n$ in variables $x$ and $y$. Under this identification, the monomial $x^{n-j}y^j$ has weight $(n-j,j)$. Recall that
\[
\begin{aligned}
V_F|_{U(2)}\simeq\mathrm{Sym}^{m+2}(U)\otimes\mathrm{det}^{-r_2-2},\qquad
V_D|_{U(2)}\simeq\mathrm{Sym}^{m+2}(U)\otimes\mathrm{det}^{-r_1-3},
\end{aligned}
\]
where $\R^{\times}\subset K_{G,\infty}$ acts by the character $z\mapsto z^{m+u_{\diamond}}$.
Therefore the weights of the monomials $x^{m+2-j}y^j$ are 
\[(r_1-j,j-r_2-2),\quad (r_2-1-j,j-r_1-3)\]
respectively.
\begin{lemma}\label{Explicit:Description:Q}
Let $\{x^\vee,y^\vee\}$ be the basis of $U^\vee$ dual to $\{x,y\}$. Then there exists a choice of generators $e_\diamond$ of $M_\diamond$ such that
\[
Q_F=(x^\vee)^{r_2+a_F}(y^\vee)^{r_1+2-a_F},
\qquad
Q_D=(x^\vee)^{r_1+1+a_D}(y^\vee)^{r_2+1-a_D}.
\]
\end{lemma}

\begin{proof}
Let $T_U\subset U(2)$ denote the diagonal maximal torus, and write $\C_{(a,b)}$ for its one-dimensional representation of weight $(a,b)$. By Remark \ref{Representation:Siegel:Compact} and \eqref{RestrictionLi}, we have the following decompositions into $T_U$-weight spaces:
\[
V_F|_{T_U}\simeq M_F\oplus
\bigoplus_{\substack{0\leq j\leq m+2\\j\neq j_F}}
\C_{(j-A_F,-B_F-j)},\quad
V_D|_{T_U}\simeq M_D\oplus
\bigoplus_{\substack{0\leq j\leq m+2\\j\neq j_D}}
\C_{(j-A_D,-B_D-j)}.
\]
Here $M_F$ and $M_D$ are the one-dimensional $T_U$-weight spaces of weights $(a_F-2,b_F-2)$ and $(a_D-2,-b_D)$, respectively. By the polynomial realizations of $V_F$ and $V_D$, these weight spaces are generated by
\[
r_F=x^{r_2+a_F}y^{r_1+2-a_F},
\qquad
r_D=x^{r_1+1+a_D}y^{r_2+1-a_D}.
\]
Since $p_\diamond$ is the projection onto $M_\diamond$, the functional $Q_\diamond$ vanishes on every other weight space. Taking $e_\diamond=p_\diamond(r_\diamond)$, we have $Q_\diamond(r_\diamond)=1$, so $Q_\diamond$ is precisely the dual coordinate functional corresponding to $r_\diamond$. Under the natural identification $(\mathrm{Sym}^{m+2}U)^\vee\simeq\mathrm{Sym}^{m+2}(U^\vee)$, this coordinate functional is a nonzero multiple of the stated dual monomial. Rescaling $e_\diamond$ gives the asserted normalization.
\end{proof}

\begin{defn}
We denote by $v_{\diamond}$ the vectors of $\tau_{i_{\diamond}}$ of weights $(-a_F,-b_F)$ and $(-a_D,b_D)$, respectively.
\end{defn}

\begin{lemma}\label{C:Evaluated}
For nonzero constants $\lambda_F,\lambda_D$ determined by the chosen bases, we have
\[
C_F(Q_F\otimes d\iota(\xi_1^-\wedge\xi_2^-))
    \overset{.}{=} v_{F},\quad
C_D(Q_D\otimes d\iota(\xi_1^-))
    \overset{.}{=} v_{D}.\]
\end{lemma}

\begin{proof}
By Lemma \ref{lem:twisted-domains}, we have
\[
\mathfrak p_G^-\simeq\mathrm{Sym}^2(U^\vee),
\qquad
\bigwedge^2\mathfrak p_G^-
\simeq\mathrm{Sym}^2(U^\vee)\otimes\mathrm{det}^{-1},
\]
as representations of $U(2)$. Since
$d\iota(\xi_1^-)$ and $d\iota(\xi_2^-)$ have
weights $(-2,0)$ and $(0,-2)$, respectively,
they correspond, up to nonzero scalars, to
$(x^\vee)^2$ and $(y^\vee)^2$ in the polynomial
model. Consequently,
\[
d\iota(\xi_1^-\wedge\xi_2^-)
\doteq x^\vee y^\vee\otimes\mathrm{det}^{-1}.
\]
By Remark \ref{Coincidence:Minimal:K:types}
and Lemma \ref{Explicit:Description:Q}
we obtain
\begin{align*}
C_F(Q_F\otimes d\iota(\xi_1^-\wedge\xi_2^-))
&\doteq
(x^\vee)^{r_2+a_F+1}
(y^\vee)^{r_1+3-a_F},\\
C_D(Q_D\otimes d\iota(\xi_1^-))
&\doteq
(x^\vee)^{r_1+3+a_D}
(y^\vee)^{r_2+1-a_D}.
\end{align*}
Their respective weights are $(-a_F,-b_F)$ and $(-a_D,b_D)$.
Since the weight spaces of an irreducible
$U(2)$-representation are one-dimensional,
the result follows.
\end{proof}
\begin{theorem}\label{Integral:New:Pairing}
Let $\Pi$ be globally generic, cuspidal and non-CAP, of general or Yoshida type, and let $\eta\in H^{i_\diamond}(\Pi_f)_{\C}$. Given $f\boxtimes g\in H^{2-i_{\diamond}}(X_H,H_{\diamond})$, we have
\[
P_\diamond(\eta,f\boxtimes g)\overset{.}{=}\int_{Z_H(\R)H(\Q)\backslash H(\A)}|\nu|_{\A}^{m+i_{\diamond}}F_\eta(v_{\diamond})(\iota(h))f(h_1)g(h_2) dh,
\]
where $F_\eta(v_{\diamond})\in \Pi^{(\diamond),\mathrm{arith}}$. When $(a_\diamond,b_\diamond)=(0,c_0)$ and the unitary representation $\Sigma$ generated by $\varphi_g$ satisfies $\omega_\Pi\omega_\Sigma=1$, we have 
\[P_{\diamond}(\eta,1\boxtimes g) \overset{.}{=}\mathcal{I}(\Psi_{\eta},\varphi_g),\]
where $\Psi_{\eta} := F_{\eta}(v_{\diamond})\otimes |\nu|_{\A}^{\frac{m+u_{\diamond}}{2}}\in \Pi$ and $\varphi_g := g\otimes|\mathrm{det}|_{\A}^{\frac{m-u_{\diamond}}{2}+i_{\diamond}}$.
\end{theorem}
\begin{proof}
The first assertion follows from a combination of Lemmas \ref{First:Formula:Pairing}, \ref{Explicit:Expresion:Wedge:Lift} and \ref{Explicit:Expression:F}, Proposition \ref{prop:twisted-factorization} and Lemma \ref{C:Evaluated}. Since $F_{\eta}(v_{\diamond})\in \Pi^{(\diamond),\mathrm{arith}}$, we can write $F_{\eta}(v_{\diamond}) = \Psi_{\eta}\otimes |\nu|_{\A}^{-\frac{m+u_{\diamond}}{2}}$, with $\Psi_{\eta}(v_{\diamond})\in \Pi$. Since $\nu(\iota(h))=\mathrm{det} (h_2)$, we obtain
\begin{align*}
|\nu|_{\A}^{m+i_{\diamond}}F_\eta(v_{\diamond})(\iota(h))g(h_2) &= |\nu(\iota(h))|_{\A}^{\frac{m-u_{\diamond}}{2}+i_{\diamond}}\Psi_{\eta}(\iota(h))g(h_2)\\
&=
\Psi_\eta(\iota(h))\varphi_g(h_2).
\end{align*}
The condition $\omega_\Pi\omega_\Sigma=1$ makes this product
invariant under $Z_G(\A)$. Passing from the quotient by
$Z_H(\R)$ to the quotient by $Z_G(\A)$ contributes only the
fixed nonzero central volume factor, proving the last identity.
\end{proof}
\begin{rk}
    Since
\[
\frac{m-u_\diamond}{2}+i_\diamond
=
\begin{cases}
c_0/2, & \diamond=F,\\
-c_0/2, & \diamond=D,
\end{cases}
\]
and the central character of $g$ is a unitary character times
$|\mathrm{det}|_{\A}^{-c_0/2}$ in region $(F)$ and
$|\mathrm{det}|_{\A}^{c_0/2}$ in region $(D)$, the form $\varphi_g$
is the unitarization of $g$. We denote by $\Sigma$ the unitary
automorphic representation generated by $\varphi_g$, commonly referred to as the unitary automorphic representation associated to $g$.
\end{rk}
\begin{rk}
    Note that, with the notation of the preceding theorem,
$\Psi_\eta$ lies in the $U(2)$-type $\tau_{i_\diamond}$
and corresponds to the vector $v_\diamond$. In particular,
it has $U(2)$-weight $(0,-c_0)$ in region $(F)$ and
$(0,c_0)$ in region $(D)$.
\end{rk}

\begin{cor}\label{cor:coherent-nonvanishing}\label{lem:infinity}
Let $\Pi$ be globally generic, non-CAP and cuspidal, of general
or Yoshida type, with the archimedean component fixed above.
Let $g$ be a nonzero cuspidal modular input of weight
$c_0$, realized as in
Theorem~\ref{Integral:New:Pairing}, and let $\Sigma$ be the
irreducible unitary automorphic representation generated by
$\varphi_g$. Assume $\omega_\Pi\omega_\Sigma=1$ and
Assumption~\ref{eq:local-compatibility}.
For every nonzero $\eta\in H^{i_\diamond}(\Pi_f)_{\C}$,
\[
\operatorname{Res}_{s=1}\Lambda(s,\Pi\times\Sigma)\neq0
\Longleftrightarrow 
P_\diamond(t\eta,1\boxtimes g)\neq0
\quad\text{for some }t\in G(\A_f).
\]
\end{cor}

\begin{proof}
A nonzero pairing gives a nonzero distinguished period by Theorem~\ref{Integral:New:Pairing}, hence a pole by Theorem~\ref{thm:complex-period}. For the converse, suppose $\Lambda(\Pi\times \Sigma,s)$ has a pole. For suitable test data, Propositions \ref{Global:Zeta:Integral:Result} and \ref{Period:Equals:Residue} show that 
\[
P_\diamond(\eta,1\boxtimes g)\overset{.}{=}\operatorname{Res}_{s=1}\Lambda(s,\Pi\times\Sigma)\,\prod_{v\in S}Q_v(1),
\]
where $S$ contains the infinite place, the ramified finite places and every finite place with nonspherical data and $Q_v(s)=Z_v(s,W_{\Psi_{\eta,v}},W_{\varphi_{g,v}},\varphi_v)/L_v(s,\Pi_v\times\Sigma_v)$. Let us assume we are in region $(F)$. We now prove that $Q_\infty(1)\ne0$. Taking the Schwartz function with archimedean component $\varphi_{\infty} = e^{-\pi(x^2+y^2)}$, one has that 
\[Z_{\infty}(s,W_{\Psi_{\eta,\infty}},W_{\varphi_{g,\infty}},\varphi_{\infty}) = \int_{Z_G(\R)N_{H}(\R)\setminus H(\R)}W_{c_0}(\iota(h))W'_{c_0}(h_2)f_{\varphi_{\infty},1}(h_1,s)dh,\]
where $W_{c_0}$ is the vector in the Whittaker model of $\Pi_{\infty}$ of compact weight $(0,-c_0)$, $W'_{c_0}$ is the vector of compact weight $c_0$ in the Whittaker model of $\Sigma_{\infty}$ and $f_{\varphi_{\infty},1}(h_1,s)$ is the function defined in \eqref{Section:f:Archimedean}. By \cite[\S 6.4]{LR2}, we then have
\[Z_{\infty}(s,W_{\Psi_{\eta,\infty}},W_{\varphi_{g,\infty}},\varphi_{\infty})\overset{.}{=}L_\infty(\Pi\times\Sigma,s),\]
concluding that $Q_{\infty}(1)\neq 0$. 

For region $(D)$, let $w_0\in H(\A)$ be trivial at finite places and
have archimedean component
$(\operatorname{diag}(1,-1),\operatorname{diag}(1,-1))$.
Simultaneous right translation of $\Psi_\eta$ by $\iota(w_0)$
and $\varphi_g$ by $w_0$ preserves the integral,
sends $\tau_2$ to $\tau_1$, and changes the compact weights
$(0,-c_0)$ and $c_0$ to $(0,c_0)$ and $-c_0$, respectively.
Thus the translated vectors realize the region $(F)$ pairing,
which is nonzero by the previous discussion and Theorem \ref{Integral:New:Pairing}.

When $v$ is a finite prime, there is a choice of data so that $Q_v(1)\neq 0$ by the GCD property of the $L$-function. Irreducibility of $\Pi_f$ and $\Sigma_f$ then implies that
this period is nonzero on some pair of finite translates
$\Pi(t)\Psi_{\eta}$ and $\Sigma(a)\varphi_g$.
Put $h_a=(\operatorname{diag}(\det a,1),a)\in H(\A_f)$.
Invariance of the period gives
\[
0\neq
\mathcal I\bigl(\Pi(t)\Psi_{\eta},\Sigma(a)\varphi_g\bigr)
=
\mathcal I\bigl(
\Pi(\iota(h_a)^{-1}t)\Psi_{\eta},\varphi_g
\bigr),
\]
since the first modular factor is the constant $1$.
Theorem~\ref{Integral:New:Pairing} now gives the required
nonvanishing, with $g$ fixed; the normalization in
\eqref{Normalization:Pi} contributes only a nonzero scalar.
Passing to a common level preserves nonvanishing by
compatibility of pullback and trace.
\end{proof}

\section{\texorpdfstring{$p$-adic $L$-functions for $\GSp_4\times\GL_2$}{p-adic L-functions for GSp4 x GL2}}\label{sec:review}

To interpolate the period formula of Section~\ref{subsec:coherent-periods}, we use the ordinary families and coherent pairings of \cite{LZvista,LZ-BK,LR3}. We construct the family inputs in this section and compare them with their classical realizations in Section~\ref{subsec:classical-comparisons}. The comparison with the classical pairing at level $p^2$ is made in Section~\ref{subsec:level-comparison}.

\subsection{$p$-adic families for $\GSp(4)$ and cohomology}\label{Section:P-adic:Families:Cohom}

We use the locally analytic cohomology of \cite[\S\S9--10]{LZ-BK} and \cite[\S\S4--5]{LR3}, based on \cite{boxerpilloni20}. Throughout, $K$ is a finite extension of $\mathbb Q_p$ containing the coefficients of our families, and $\mathcal W$ is the rigid analytic space of continuous characters of $\mathbb Z_p^\times$. An integer $n$ denotes the character $z\mapsto z^n$. All tame levels and finite-order character components are fixed.

Let $U\subset\mathcal W^2$ be an affinoid neighbourhood, with universal characters $\hr_1,\hr_2$, and put $\nu_U=(\hr_1,\hr_2;\hr_1+\hr_2)$. As in \cite{LZ-BK}, we write weight characters additively, so their sum denotes their product. The classical Siegel weights are $k_i=r_i+3$. Higher Coleman theory constructs an eigenvariety $\mathscr E\to\mathcal W^2$ and coherent cohomology sheaves $H^i(\mathcal M^{\bullet,-,\mathrm{fs}}_{\mathrm{cusp},w_j})$ on it; the superscript $\bullet$ refers to a complex, whereas $H^i$ denotes its cohomology. These sheaves record Hecke eigensystems in locally analytic overconvergent cohomology. A family is specified by its map to this eigenvariety, together with the conditions below.

\begin{defn}\label{def:automorphic}
A family $\underline\Pi$ over $U$ consists of a smooth finite flat covering $\widetilde U\to U$ and a map $j_{\underline\Pi}:\widetilde U\to\mathscr E$ over weight space, with fixed tame level and finite-order central character, satisfying the conditions of \cite[Definition 10.4.1]{LZ-BK}. In particular, the pulled-back cohomology sheaves vanish unless $i+j=3$, and \[ S^i(\underline\Pi):=j_{\underline\Pi}^*H^i(\mathcal M^{\bullet,-,\mathrm{fs}}_{\mathrm{cusp},w_{3-i}}) \] are locally free of rank one for $i=1,2$. They have rank one also for $i=0,3$ in general type, and vanish for $i=0,3$ in Yoshida type. After shrinking, we trivialize the nonzero line bundles. The central character is normalized as in loc.\,cit., namely $|\cdot|^{-(r_1+r_2)}\widehat\chi_{\Pi}$ at classical weights. We assume throughout that the chosen refinement is Borel ordinary.
\end{defn}

A classical point $P\in\widetilde U$ is called good in degree $i$ if its eigensystem comes from a classical, globally generic cuspidal representation with the chosen refinement and the specialization map identifies $S^i(\underline\Pi)\otimes_{\mathcal O(\widetilde U)}K(P)$ with its classical coherent-cohomology eigenline at the specified level. We use this condition in both degrees $1$ and $2$ (it is the condition of \cite[Definition 10.4.2]{LZ-BK} in degree $2$, and its degree-$1$ counterpart in \cite[\S5.2]{LR3}).

We will use ordinary control to specialize the selected Siegel classes and their tame translates.

\begin{lemma}\label{lem:ordinary-control}
Let $P$ be an unramified Borel-ordinary, globally generic point of a family as in Definition~\ref{def:automorphic}, at quasiparamodular new level away from $p$. After shrinking the given cover $\widetilde U$, the localized complex in each of the two Weyl realizations is a free rank-one module concentrated in its coherent degree, and its formation commutes with coefficient specialization. The same holds for the images of these lines under the fixed tame translations and pullbacks to a common larger level.
\end{lemma}
\begin{proof}
The ordinary control result follows from \cite[Proposition 1.4.7 and Theorems 1.4.9--1.4.10]{Boxer:Pilloni:Siegel}. The small-slope bounds for the normalized operators are those of \cite[Proposition 3.3.2]{LZ-BK} and \cite[Proposition 6.1]{LR3}:
\[
\begin{aligned}
v_p(\mathcal U'_{\Si})<r_2+1,\quad v_p(\mathcal U'_{\Kl})<r_1-r_2+1&\qquad\text{in degree }2,\\
v_p(\mathcal U'_{\Si})<r_1+2,\quad v_p(\mathcal U'_{\Kl})<r_1-r_2+1&\qquad\text{in degree }1.
\end{aligned}
\]
Both ordinary slopes are zero, so these bounds apply for $r_1\geq r_2\geq0$. At new tame level, the classical localization is concentrated in degree $i$ and is one-dimensional at the $p$-regular refinement \cite[Proposition 2.7.2 and \S2.8]{LZ-BK}; its four root valuations are $0,r_2+1,r_1+2,r_1+r_2+3$. In a finite free representative of the localized perfect complex, cancel the differential entries invertible at $P$. Its remaining fibre has a single rank-one term in degree $i$, and hence, after shrinking, so does the complex. This free line commutes with coefficient specialization. Tame translation is an isomorphism at the conjugate level, and pullback to a common larger level is injective because trace composed with pullback is multiplication by the degree. Their images therefore have the same specialization property.
\end{proof}

\begin{lemma}\label{lem:Yoshida-family}
Let $P$ be a Yoshida point at quasiparamodular new tame level. Suppose the localized finite-slope complexes in the degree-zero and degree-three Weyl realizations satisfy derived classical control at $P$. Then, after shrinking the family neighbourhood, both complexes vanish. In particular, every nearby non-CAP classical point satisfying control in these realizations is of Yoshida type.
\end{lemma}
\begin{proof}
By \cite[Theorem 5.2]{LPSZ1}, the Yoshida eigensystem has zero localization in degrees zero and three at this level. Derived control identifies these zero spaces with the corresponding derived fibres of the localized perfect complexes. A perfect complex with zero derived fibre is acyclic in a neighbourhood: over the local ring one successively cancels each differential entry that is a unit, and Nakayama's lemma leaves no terms. Thus both complexes vanish after shrinking. At a nearby controlled point of general type, the same theorem gives a nonzero classical summand in each of these degrees, contradicting this vanishing.
\end{proof}

The pole criterion uses only the degree-one and degree-two lines of Lemma~\ref{lem:ordinary-control}; the additional control hypotheses in Lemma~\ref{lem:Yoshida-family} concern variation within the Yoshida locus.

Similarly, let $\widetilde U'\to U'\subset\mathcal W$ carry an ordinary $\GL_2$ family $\underline\Sigma$ of modular forms $\mathcal G$. We write $c$ for the modular weight; thus the cohomological weight variable used in \cite[\S10.5]{LZ-BK} is $c-2$. Write $G(\chi_\Sigma^{-1})$ for the Gauss sum of the fixed nebentype. We use the line $S^0(\underline\Sigma)$ generated by $G(\chi_\Sigma^{-1})\mathcal G$, and a locally free rank-one module $S^1(\underline\Sigma)$ as in \cite[\S5.2]{LR3}. At the specializations used below we require the relevant classicality and base-change isomorphisms. These hypotheses also specify the rational classical eigenlines used to define periods.

The ordinary weight algebras of \cite{TU99} give the same coordinates: on a fixed finite-character component their parameters specialize to $(1+p)^{r_i}-1$ on the Siegel side and $(1+p)^{c-2}-1$ on the $\GL_2$ side; the finite covers also record the Hecke eigenvalues.

\subsection{The triple product and its Eisenstein degeneration}\label{subsec:triple-product}

The pairings in regions (F) and (D) arise from the same construction for $\GSp_4\times\GL_2\times\GL_2$. Let $\Sigma_1,\Sigma_2$ be cuspidal representations of modular weights $c_1,c_2$, with $\chi_\Pi\chi_{\Sigma_1}\chi_{\Sigma_2}=1$. For forms $\varphi\in\Pi$ and $f_i\in\Sigma_i$, the Gross--Prasad period is
\[
\mathcal P(\varphi,f_1,f_2)=\int_{Z_G(\A)H(\Q)\backslash H(\A)}\varphi(\iota(h))f_1(h_1)f_2(h_2)\,dh.
\]
For tempered representations satisfying the local distinction hypotheses, the global Gross--Prasad conjecture relates the nonvanishing of this period to the central value of the degree-sixteen $L$-function $\Lambda(\Pi\times\Sigma_1\times\Sigma_2,\tfrac12)$; its refined form expresses the squared absolute value of the period in terms of this value, adjoint $L$-values and local factors; see \cite[\S2.2]{LZvista}. Up to interchanging the two modular factors, the six regions (a)--(f) of \cite[\S2.3]{LZvista} describe the possible interlacings of the archimedean weights. For fixed $\Pi$, the expected $p$-adic interpolation of these periods in the regions with the appropriate local signs, and the corresponding Euler factors, are described in \cite[Conjecture 4.1.2 and \S4.3]{LZvista}.

The two coherent realizations used here lie on edges of regions (f) and (d). When $c_1+c_2=k_1-k_2+2$, the restriction of a degree-two Siegel class pairs with two degree-zero modular classes, giving a product of degrees $(2,0,0)$ as in \cite{LPSZ1}. When $c_2-c_1=k_1-k_2+2$, the degrees are instead $(1,0,1)$: the second modular input is a Serre-dual class, and the Siegel class has degree one \cite[\S3.4]{LR2}. These are the two branching relations of Section~\ref{Section:Identification:Modular:Forms}; the second requires the degree-one elliptic realization used in region (D).

Replacing $f_1$ by an Eisenstein series of weight $c_1=t$, and writing $c_2=c$, turns the Gross--Prasad integral into a Novodvorsky integral for $\Pi\times\Sigma_2$; see \cite[\S1]{LR3}. This is the origin of the trilinear local functionals in Section~\ref{subsec:local-improvement}, where the first modular representation is a principal series, possibly reducible. On the Eisenstein component used below, the ordinary arithmetic root is $\mathfrak a_1=1$. We compare its ordinary vector with its $p$-depletion; their local integrals differ by the Euler factor removed by improvement. The local parameters will be renormalized in Section~\ref{subsec:local-improvement}, so the root $1$ becomes $p^{1-t/2}$ there. The two weight relations give $t=k_1-k_2+2-c$ and $t=c-(k_1-k_2+2)$, respectively. Both therefore meet the weight-zero Eisenstein specialization at $c=k_1-k_2+2$, where the residue replaces the first modular form by the constant function.

Write $(\alpha,\beta,\gamma,\delta)$ for the motivic Hecke parameters of $\Pi_p$ and $(\mathfrak a_i,\mathfrak b_i)$ for those of $\Sigma_{i,p}$, ordered by valuation. The central motivic argument of the triple tensor product is $s'=(k_1+k_2+c_1+c_2-4)/2$. In the unramified ordinary case, the sixteen products have the following valuations:
\begin{center}
\small
\begin{tabular}{c|cccc}
&$\mathfrak a_1\mathfrak a_2$&$\mathfrak a_1\mathfrak b_2$&$\mathfrak b_1\mathfrak a_2$&$\mathfrak b_1\mathfrak b_2$\\\hline
$\alpha$&$0$&$c_2-1$&$c_1-1$&$c_1+c_2-2$\\
$\beta$&$k_2-2$&$k_2+c_2-3$&$k_2+c_1-3$&$k_2+c_1+c_2-4$\\
$\gamma$&$k_1-1$&$k_1+c_2-2$&$k_1+c_1-2$&$k_1+c_1+c_2-3$\\
$\delta$&$k_1+k_2-3$&$k_1+k_2+c_2-4$&$k_1+k_2+c_1-4$&$k_1+k_2+c_1+c_2-5$
\end{tabular}
\end{center}
An interpolation factor $1-p^{s'-1}/(xyz)$ can vanish only if the corresponding entry equals $s'-1$; the unit parts must agree as well. For example, $x=\beta$, $y=\mathfrak a_1$, $z=\mathfrak b_2$ requires $k_1+c_1=k_2+c_2$. Whether this factor occurs in a particular interpolation formula depends on its critical region. We refer to \cite[\S4.3]{LZvista} for that choice.

\subsection{The degree-eight case and its critical regions}\label{subsec:degree-eight}

In the degree-eight problem, the variables are $(k_1,k_2,c,s)$: the two Siegel weights, the remaining elliptic weight and the cyclotomic variable. The complex function $L(\Pi\times\Sigma,s)$ is attached to the tensor product of the spin representation of the $L$-group of $\GSp_4$ and the standard representation for $\GL_2$. For our cohomological weights $(k_1,k_2)$ and modular weight $c$ of $\Sigma$, the motivic weights are $k_1+k_2-3$ and $c-1$, respectively. In motivic notation the normalization is
\[
L(\Pi\times\Sigma,s)=L\left(M(\Pi)\otimes M(\Sigma),s+\frac w2\right),\qquad
w=k_1+k_2+c-4.
\]
For $s$ to be critical for $L(\Pi \times \Sigma,s)$, we must have $s = -\tfrac{w}{2} \bmod \mathbb Z$, and $s$ must satisfy various inequalities depending on how the weights of $\Pi$ and $\Sigma$ interlace. For fixed $(k_1, k_2)$, the allowed pairs $(s,c)$ form three disjoint polygonal regions in the plane (all symmetric about the line $s = \tfrac{1}{2}$).

We denote these regions by (A), (D) and (F), following \cite{LR3}. Their inequalities are as follows; the diagram in the introduction shows their relative positions. We work with the ordinary families and good classical points fixed above.
\begin{itemize}
    \item {\bf Region (A).} Described by the set of points of weights $(k_1,k_2,c,s)$ with \[ \begin{cases} k_1+k_2-1 \leq c, \\ |2s-1| \leq c-(k_1+k_2-1). \end{cases} \]
    \item {\bf Region (D).} Described by the set of points of weights $(k_1,k_2,c,s)$ with \[ \begin{cases} k_1-k_2+3 \leq c \leq k_1+k_2-3, \\ |2s-1| \leq \min(k_1+k_2-3-c, c-(k_1-k_2+3)). \end{cases} \]
    \item {\bf Region (F).} Described by the set of points of weights $(k_1,k_2,c,s)$ with \[ \begin{cases} 1 \leq c \leq k_1-k_2+1, \\ |2s-1|  \leq k_1-k_2+1-c. \end{cases} \]
\end{itemize}

We use the higher-cohomology construction in region (F) of \cite{LPSZ1,LZ-BK}, and the degree-$1$ construction for region (D) in \cite{LR3}. The extensions in \cite{GrahamRockwood} and the automorphic constructions of \cite{liu1,liu2} give further approaches to these functions. Here we work with the three-variable cohomological pairings.

The change between the unitary automorphic variable $s$ and the motivic variable $s'$ is
\begin{equation}\label{eq:normalizations}
s'=s+\frac{k_1+k_2+c-4}{2}.
\end{equation}
The two domains studied below therefore have the following descriptions:
\[
\begin{array}{c|c|c|c}
 &\text{Eisenstein weight}&s'&s\\ \hline
(F)&t_F=k_1-k_2-c+2&k_1-1&t_F/2\\
(D)&t_D=-k_1+k_2+c-2&c+k_2-3&t_D/2
\end{array}
\]
At $c_0=k_1-k_2+2$, both domains have $t=0$, $s'=k_1-1$ and $s=0$. At this boundary, the residue of the Eisenstein family is a constant modular form. The complex pole at $s=1$, corresponding to $s'=k_1$, is connected to this boundary through the automorphic period: Section~\ref{sec:improved} compares the residue of the $p$-adic pairing with its complex realization. The Rankin functional equation will enter separately in the factorization discussion of Section~\ref{sec:factorization}.

\subsection{Duality in families}\label{subsec:family-duality}
To form the pairings in families, let $\widetilde U$ and $\widetilde U'$ be the parameter spaces of the ordinary Siegel and modular families of Section~\ref{Section:P-adic:Families:Cohom}. After shrinking, put $A=\mathcal O(\widetilde U)\widehat\otimes_K\mathcal O(\widetilde U')$. We retain the chosen character components, refinements and tame translates, and set
\[
(i_F,w_F,w_{G,F})=(2,\mathrm{id},w_1),\qquad
(i_D,w_D,w_{G,D})=(1,w_{01},w_2).
\]
The locally analytic cuspidal complexes for $G$ are those of \cite[\S9.5]{LZ-BK} and \cite[\S4.8]{LR3}. For sufficiently large auxiliary level $p^T$, translated pullback and evaluation against the branching kernel map their degree-$i_\diamond$ classes to the negative realization on $H$. Denote this composite by $r_{\diamond,A}$.

Write $\nu_A=(\nu_{1,A},\nu_{2,A};m_A)$ for the pullback of the Siegel weight to $A$, with $m_A=\nu_{1,A}+\nu_{2,A}$. Let $b_A$ be the weight of $\underline\Sigma$, and let $a_A$ be the Eisenstein weight, determined by $a_A+b_A=\nu_{1,A}-\nu_{2,A}+2$ in $(F)$ and $b_A-a_A=\nu_{1,A}-\nu_{2,A}+2$ in $(D)$. The character of the negative coefficient on $H$ is
\[
\tau_A=\begin{cases}(a_A-2,b_A-2;m_A),&\diamond=F,\\(a_A-2,-b_A;m_A),&\diamond=D.\end{cases}
\]
The calculation in Lemma~\ref{lem:D-coefficient-character} verifies the second formula in the branching convention of \cite{LR3}.

Here $S_{H,\mathrm{Iw}(p^T)}$ denotes the adic realization of the toroidal compactification of $X_H$ at the indicated $p$-level, with the fixed tame level understood. The symbols $+$ and $-$ distinguish the analytic induction and distribution realizations; $\dagger$ denotes the limit over the neighbourhoods and supports used in the cited constructions.

\begin{defn}\label{def:analytic-pairing}
We use the $A$-bilinear pairing
\begin{equation}\label{p-adic:Serre}
\begin{aligned}
H^{i_\diamond}_{w_\diamond,\mathrm{an}}(S_{H,\mathrm{Iw}(p^T)},\tau_A,\mathrm{cusp})^{-,\dagger}
\times H^{2-i_\diamond}_{w_\diamond,\mathrm{an}}(S_{H,\mathrm{Iw}(p^T)},\tau_A)^{+,\dagger}
&\longrightarrow A,\\
(x,h)&\longmapsto\langle x,h\rangle_{\mathrm{an}},
\end{aligned}
\end{equation}
of \cite[Theorem 9.6.2]{LZ-BK} in region $(F)$ and \cite[Theorem 4.12]{LR3} in region $(D)$. The negative factor is placed first in the cup product. For a coefficient specialization $Q:A\to K$, write $B_{\diamond,Q}$ for its scalar extension: $B_{\diamond,Q}(x\otimes u,h\otimes v)=uv\,Q(\langle x,h\rangle_{\mathrm{an}})$.
\end{defn}

The modular inputs have degrees $(0,0)$ in $(F)$ and $(0,1)$ in $(D)$. We now construct the degree-one input and its coefficient transport.

\subsubsection{The degree-one modular input}\label{subsubsec:GL2-input}
Put $B=\mathcal O(\widetilde U')$, and fix a generator $g^\vee$ of the holomorphic ordinary line for the contragredient modular family. Write $\mathfrak a$ for the ordinary $U_p$-eigenvalue of $\underline\Sigma$. The modified Serre duality of \cite[\S5.7.3]{BP:Modular} realizes the dual line on the \'etale ordinary neighbourhood, in the positive long-Weyl realization. We take this as our model for $S^1(\underline\Sigma)$ and normalize its generator $\underline\upsilon$ by $P(g^\vee,\underline\upsilon)=1$. Here $P(f,n)=P_0(f,w^*n)$, where $P_0$ is ordinary Serre duality and $w^*$ is the normalized Atkin--Lehner map. Thus $\underline\upsilon$ is the image under $(w^*)^{-1}$ of the $P_0$-dual generator.

\begin{lemma}\label{lem:modular-adjunction}
For the normalized modular operators, one has
\begin{equation}\label{eq:elliptic-adjoint}
P(U_pf,n)=P(f,\langle p\rangle^{-1}U_pn),\qquad
\langle p\rangle^{-1}U_p\underline\upsilon=\mathfrak a^\vee\underline\upsilon,
\quad \mathfrak a^\vee=\chi_\Sigma(p)^{-1}\mathfrak a.
\end{equation}
\end{lemma}
\begin{proof}
The first equality is \cite[Lemma 5.22]{BP:Modular}, with the normalized isogeny and trace maps used there. The holomorphic vector $g^\vee$ has ordinary eigenvalue $\mathfrak a^\vee$. Applying adjunction to $g^\vee$ and using the perfect pairing on the selected lines proves the second equality. Equivalently, $U_p$ on the degree-one $\Sigma$-line has eigenvalue $\mathfrak a$, and its diamond operator acts by $\chi_\Sigma(p)$.
\end{proof}

\subsubsection{Coefficients and supports}\label{subsubsec:family-support}

In region $(D)$ write $R_1=r_1+u$ and $R_2=r_2+v$ for the Siegel weights and $b$ for the weight of the cuspidal modular family. The Eisenstein weight is then $a=b-R_1+R_2-2$. On the chosen smooth modular branch, the family space $X=\operatorname{Sp}(A)=\widetilde U\times\widetilde U'$ is reduced and finite flat over the weight coordinates $(R_1,R_2,b)$. The change of variables from $b$ to $a$ gives the coordinates used in Section~\ref{subsec:level-comparison}.

Write $\omega_i$ for the Hodge bundle on the $i$-th modular curve. The degree-one class of Section~\ref{subsubsec:GL2-input} has Hodge weight $2-b$. To pair its product with the Eisenstein input against the Siegel class, we must also specify the determinant character. The two factors of $H$ have the same determinant $\nu$, and we place this character on the second curve. Indeed, in our inverse-frame coordinates the character of $\omega_H(a,2-b;-R_1-R_2)$ is
\[
t_1^{-a}t_2^{b-2}\nu^{(a+2-b-R_1-R_2)/2}
=t_1^{-a}t_2^{b-2}\nu^{-R_1},
\]
where the equality follows from $b-a=R_1-R_2+2$. Thus, writing $\mathrm{det}_2$ for the determinant coefficient line on the second curve, the product coefficient is
\[
\omega_1^a\boxtimes\bigl(\omega_2^{2-b}\otimes(\mathrm{det}_2)^{-R_1}\bigr).
\]

The power $(\mathrm{det}_2)^{-R_1}$ is defined in families by its character on changes of frame. If $t\in\mathbb Z_p^\times$ is the determinant of such a change on the second curve, the character on the weight component through $r_1$ is
\begin{equation}\label{eq:determinant-family}
t\longmapsto t^{-r_1}\delta_u(t)^{-1},
\qquad
\delta_u(t)=\langle t\rangle^u,
\end{equation}
where $\langle t\rangle$ denotes the principal-unit part of $t$. We use its analytic extension to the reduced torsor as in \cite[Definition 4.5]{LR3}. At an integral weight $R_1=r_1+n$ on this component, the characters $t^{r_1+n}$ and $t^{r_1}$ agree on the torsion subgroup of $\mathbb Z_p^\times$. Hence $\langle t\rangle^n=t^n$, and \eqref{eq:determinant-family} specializes to $t^{-R_1}$.

The translated Siegel branching quotient has character $t_1^{a-2}t_2^{-b}\nu^{R_1+2}$. Its contraction with the modular coefficient has character
\begin{equation}\label{eq:canonical-family}
\bigl(t_1^{a-2}t_2^{-b}\nu^{R_1+2}\bigr)
\bigl(t_1^{-a}t_2^{b-2}\nu^{-R_1}\bigr)
=t_1^{-2}t_2^{-2}\nu^2.
\end{equation}
By Kodaira--Spencer this is the character of $\omega_H(2,2;0)=K_{X_H}(D_H)$. The boundary map \eqref{Pullback:Vector:Bundles} supplies the factor $\mathcal O(-D_H)$ from the cuspidal Siegel coefficient, so the contracted product takes values in $K_{X_H}$.

We also need this determinant transport to agree with the normalized isogeny maps. For the elliptic $U_p$-correspondence, work at full $p$-power level and let $(e_1,e_2)$ be a basis of $T_pE$. If $\pi:E\to E'$ is the quotient by the subgroup generated by $e_2+je_1$ modulo $p$, with $0\leq j<p$, the target Tate basis is
\[
e'_1=\pi(e_1),\qquad e'_2=p^{-1}\bigl(\pi(e_2)+j\pi(e_1)\bigr)
\]
by the proof of \cite[Lemma 3.2.4]{RC:EichlerShimura}. Thus $e'_1\wedge e'_2=p^{-1}(\bigwedge^2\pi)(e_1\wedge e_2)$. The degree-$p$ isogeny multiplies the Weil pairing by $p$, so $\langle e'_1,e'_2\rangle_{E'}=\langle e_1,e_2\rangle_E$. Thus each branch preserves the normalized determinant frame. The action on degree-one cohomology uses the normalized operator of Lemma~\ref{lem:modular-adjunction}; Proposition~\ref{prop:classical-specialization} compares its classical normalization with that of Boxer--Pilloni.

The pairing also uses the modular degeneracy $(h_1,h_2)\mapsto(h_1d_T,h_2)$, where $d_T=\operatorname{diag}(p^T,1)$; see \cite[\S\S3.5 and 4.1]{LR3}. The induced change of structure group is conjugation by $(d_T,1)$, as described after Proposition~4.4 of that reference. It acts on the first curve, so its coefficient map is the normalized isogeny map on $\omega_1^a$ tensored with the identity on $\omega_2^{2-b}\otimes(\mathrm{det}_2)^{-R_1}$. In particular, the determinant coefficient on the second curve remains unchanged as $T$ varies.

Fix a common analytic radius $n$ and an integer $M>3n$. In region $(D)$ take $T\geq3M+1$ as in \cite[\S4.8]{LR3}. By Proposition~3.8 of that reference, the inverse image of the common support under the translated map is closed in the ambient $H$-surface and lies in the neighbourhood on which the cup product is defined. Excision therefore allows us to take its trace on $X_H$. In region $(F)$ we use the ordinary neighbourhoods and supports of \cite[\S\S6 and 9.5--9.6]{LZ-BK}. In both regions the Siegel pullback is the negative input of \eqref{p-adic:Serre} and the modular product is the positive input.

\begin{prop}\label{prop:supported-pairing}
At such an auxiliary level, translated restriction, the branching quotient, cup product and trace define an $A$-multilinear scalar pairing of the selected Siegel class and the two modular inputs. Its formation commutes with coefficient specialization. This applies to the depleted input and to the regularized ordinary input of Lemma~\ref{lem:regularized-eisenstein}.
\end{prop}
\begin{proof}
Let $C_G,C_1,C_2$ denote the supported coefficient complexes for the Siegel class and the two modular inputs. Let $V$ be their common neighbourhood on the $H$-surface, and $W\subset V$ the closed support supplied by the translated pullback. The pairing is induced by
\[
\begin{aligned}
C_G\widehat\otimes_A^{\mathbf L}C_1\widehat\otimes_A^{\mathbf L}C_2
&\longrightarrow R\Gamma_W(V,K_{X_H}\widehat\otimes A)\\
&\xrightarrow{\mathrm{excision}}R\Gamma_W(X_H,K_{X_H}\widehat\otimes A)
\longrightarrow R\Gamma(X_H,K_{X_H}\widehat\otimes A)
\xrightarrow{\Tr}A[-2].
\end{aligned}
\]
The first arrow is translated pullback, evaluation on the branching kernel and cup product, with the Siegel factor first. The coefficient calculation above, including the boundary map \eqref{Pullback:Vector:Bundles}, puts its image in $K_{X_H}$. The remaining arrows are defined because $W$ is closed in the surface and contained in $V$. The base-change compatibility of \cite[Theorem 9.6.2]{LZ-BK} and \cite[Theorem 4.12]{LR3} applies to this composite. Lemma~\ref{lem:regularized-eisenstein} places the regularized input in the same coefficient module, while depletion is a finite combination of isogeny maps commuting with specialization. Proposition~\ref{prop:classical-specialization} will identify the resulting specialized classes with classical classes.
\end{proof}

\subsection{Exceptional zeros and analytic Euler factors}\label{subsec:analytic-euler-factors}

In the degree-eight setting we write $\mathcal W_F^{\imp}$ and $\mathcal W_D^{\imp}$ for the three-dimensional subspaces described by the table in Section~\ref{subsec:degree-eight}. At algebraic points these have $s'=k_1-1$ and $s'=c+k_2-3$, respectively. The character equations, including the finite-order parts, define these subspaces. On either domain, $s'$ is determined by the three weights.

On the ordinary Siegel branch the normalized unit parameters $\beta^\circ$ and $\gamma^\circ$ have classical values $\beta/p^{k_2-2}$ and $\gamma/p^{k_1-1}$. They are defined by normalized Hecke eigenvalues: division of the ordinary Klingen eigenvalue $\alpha\beta/p^{k_2-2}$ by the ordinary Siegel eigenvalue $\alpha$ gives $\beta^\circ$, and
\[
\gamma^\circ=\chi_\Pi(p)/\beta^\circ.
\]

To determine the analytic factors removed by improvement, suppose that $\chi_\Pi\chi_1\chi_2=1$ and $\mathfrak a_i\mathfrak b_i=\chi_i(p)p^{c_i-1}$. For (F), take $(c_1,c_2)=(t,c)$; for the following (D) calculation take $(c_1,c_2)=(c,t)$, so the Eisenstein factor is the second one. Interchanging the two $\GL_2$ factors changes these subscripts accordingly.

\begin{prop}\label{prop:analytic-euler-factors}
The factors
\[
1-\frac{p^{s'-1}}{\beta\mathfrak b_1\mathfrak b_2},
\qquad
1-\frac{p^{s'-1}}{\gamma\mathfrak a_1\mathfrak b_2}
\]
extend analytically along the domains for (F) and (D), respectively. Their restrictions are
\[
1-\gamma^\circ\mathfrak a_1\mathfrak a_2,
\qquad
1-\frac{\chi_1(p)\beta^\circ\mathfrak a_2}{\mathfrak a_1}.
\]
\end{prop}

\begin{proof}
Before imposing the weight relation, the two factors are respectively
\[
1-\frac{\chi_\Pi(p)\mathfrak a_1\mathfrak a_2}{\beta^\circ}
 p^{(k_1-k_2-c-t+2)/2},
\qquad
1-\frac{\chi_\Pi(p)\chi_1(p)\mathfrak a_2}{\gamma^\circ\mathfrak a_1}
 p^{(-k_1+k_2+c-t-2)/2}.
\]
The displayed powers of $p$ have exponent zero on their respective domains. The remaining expressions involve only fixed characters and analytic units.
\end{proof}

For the Eisenstein refinement with ordinary root $1$, denote these analytic factors by
\begin{equation}\label{eq:removed-factors}
e_F=1-\gamma^\circ\mathfrak a,\qquad e_D=1-\frac{\beta^\circ}{\mathfrak a^\vee},\qquad \mathfrak a^\vee=\chi_\Sigma(p)^{-1}\mathfrak a.
\end{equation}
Here $\mathfrak a$ is the cuspidal ordinary root and $\mathfrak a^\vee$ is the dual root of Lemma~\ref{lem:modular-adjunction}. Proposition~\ref{teo:interpolation} will show that $e_F$ and $e_D$ are precisely the factors removed by improvement.

The smaller Yoshida constituent has modular weight $k_1-k_2+2$, so the common possible Eisenstein pole lies at its weight. The other constituent has weight $k_1+k_2-2$ and is not detected by these two weight-zero Eisenstein specializations. At this boundary, a zero of the Euler factor may cancel a pole of the improved function. After fixing the Siegel point and the other family data, any identity $L_p=e_\diamond L_p^{(\diamond),\imp}$ gives the one-variable order formula $\operatorname{ord}(L_p)=\operatorname{ord}(e_\diamond)+\operatorname{ord}(L_p^{(\diamond),\imp})$ at points where these functions are not identically zero. We will use this order formula on the fixed Siegel fibre when studying the poles.

\section{\texorpdfstring{Improved $p$-adic $L$-functions and Yoshida lifts}{Improved p-adic L-functions and Yoshida lifts}}~\label{sec:improved}

The ordinary Eisenstein family has a simple pole at weight zero, with constant residue. Inserting it into the pairings of Section~\ref{sec:review} gives the improved functions. We compare their specializations with the distinguished period and their depleted counterparts with the original $p$-adic $L$-functions.

\subsection{Construction of improved $p$-adic $L$-functions}\label{subsec:improved-construction}

To replace the depleted Eisenstein input of \cite{LPSZ1,LR3} by the ordinary family of \cite{LZ25-poles}, we first place its regularization in the analytic coefficient module. We keep the generators, Gauss sums and tame translations of Definition~\ref{def:analytic-pairing}. As in \cite[Definition 10.6.2]{LZ-BK}, the global functions include the tame volume multiplier $\nu_K=\operatorname{vol}(\overline K_H^p)$, which we keep separate from the geometric Serre trace $\Tr$.

Let $\Lambda = \mathbb Z_p[\![\mathbb Z_p^{\times}]\!]$ be the Iwasawa algebra of $\mathbb Z_p^{\times}$, and $\Lambda' = \mathbb Z_p[\![\mathbb Z_p^{\times} \times \mathbb Z_p^{\times}]\!]$. We also write $\kappa$ for the canonical character $\mathbb Z_p^{\times} \rightarrow \Lambda^{\times}$, and $\kappa_1, \kappa_2$ for the two canonical characters into $(\Lambda')^{\times}$. Fix a positive integer $N \geq 1$ coprime to $p$. We begin by defining the usual two-parameter Eisenstein series. Write $\mathcal M_{\kappa}(N,\Lambda')$ for the space of $p$-adic modular forms of weight $\kappa$, tame level $N$ and with coefficients in $\Lambda'$. We use additive notation for characters: $\kappa-1$ means $z\mapsto\kappa(z)z^{-1}$, and $\kappa=0$ means the trivial character. Coefficients are extended to contain the required roots of unity, with the embeddings fixed as in Section~\ref{sec:GL2xGL2}.

\begin{defn}\label{def:depleted-eisenstein}
For $N \geq 1$ coprime to $p$, let \[ \mathcal E_{1/N}(\kappa_1,\kappa_2) \in \mathcal M_{\kappa_1+\kappa_2+1}(N,\Lambda') \] stand for Katz's two-parameter family of Eisenstein series of level $N$, with $q$-expansion given by \[ \sum_{\substack{u,v \in \mathbb Z \\ uv>0, \, p \nmid uv}} u^{\kappa_1}v^{\kappa_2}\sign(u)\exp(2\pi iv/N)q^{uv}. \] 
\end{defn}

Relaxing the condition $p\nmid v$ gives the ordinary Eisenstein family that we use for improvement. Write $\mathcal Q(\Lambda)$ for the ring of fractions of $\Lambda$.

\begin{defn}\label{def:improved-eisenstein}
The improved Eisenstein family $\mathcal E^\flat_{1/N}(\kappa)\in\mathcal M_\kappa(N,\Lambda)\otimes_\Lambda\mathcal Q(\Lambda)$ has $q$-expansion \[ \mathcal E^\flat_{1/N}(\kappa)=\zeta_p(1-\kappa)+ \sum_{\substack{u,v\in\mathbb Z\\uv>0,\ p\nmid u}} u^{\kappa-1}\sign(u)\exp(2\pi iv/N)q^{uv}. \]
Its $p$-depleted form is $(\mathcal E^\flat_{1/N}(\kappa))^{[p]}=\mathcal E_{1/N}(\kappa-1,0)$.
\end{defn}

Write $\ell(\kappa)=\kappa'(1)$ for the logarithmic weight. We use the Kubota--Leopoldt normalization $\operatorname{Res}_{s=1}\zeta_p(s)=1-p^{-1}$; see \cite{Fox03}. Since its argument here is $1-\kappa$, the residue in the weight coordinate has the opposite sign:
\begin{equation}\label{eq:eisenstein-residue}
\left.\ell(\kappa)\mathcal E^\flat_{1/N}(\kappa)\right|_{\kappa=0}
=-(1-p^{-1})\mathbf1.
\end{equation}
For $N=1$ and even positive integral weight, pairing $(u,v)$ with $(-u,-v)$ gives the nonconstant coefficients $2\sum_{d\mid n,\ p\nmid d}d^{\kappa-1}$, consistent with this constant term. We use the trivial tame Eisenstein-character component, whose idempotent fixes $\mathbf1$. To insert the residue into our pairing, we need \eqref{eq:eisenstein-residue} to hold in the analytic coefficient module.

\begin{lemma}\label{lem:regularized-eisenstein}
On a sufficiently small affinoid $U$ around the trivial character, with ring $A_U=\mathcal O(U)$, the family
\[
R(\kappa)=\ell(\kappa)\mathcal E^\flat_{1/N}(\kappa)
\]
is a regular $A_U$-valued overconvergent modular form with canonical extension. Its specialization at zero is $-(1-p^{-1})\mathbf1$. The same holds after external product with either $\mathrm{GL}_2$ cuspidal input used below.
\end{lemma}
\begin{proof}
The measure-valued construction of \cite[Proposition 3.3]{LZ25-poles} gives the regularized family in the module of ordinary Katz forms. Its nonconstant coefficients vanish at zero after multiplication by $\ell(\kappa)$, and its constant coefficient has the value in \eqref{eq:eisenstein-residue}.

To place it in the analytic pairing module, choose $f>0$ with $p^f\equiv1\pmod N$. Applying $U_p$ replaces the tame root of unity $\zeta_N$ in the coefficient formula by $\zeta_N^p$: the condition $p\nmid u$ means that the power of $p$ removed from $uv$ comes entirely from $v$. The constant coefficient is unchanged. The $q$-expansion principle consequently gives $U_p^fR=R$. Thus $U_p$ is invertible on the finite span of these tame translates, and the ordinary projector fixes $R$.

By \cite[Proposition 6.2]{Pil13}, an ordinary Katz family is overconvergent. Its proof realizes the ordinary part in a finite projective module over a sufficiently small weight affinoid, with a map to Katz forms that remains injective on fibres. After shrinking $U$, choose finitely many generators of this module. Each is represented on a strict neighbourhood of the ordinary locus, so all are defined on a common neighbourhood. Restricting further gives a neighbourhood and analytic radius on which the coefficient maps for our pairings are defined.

The specialization of $R$ maps to $-(1-p^{-1})\mathbf1$ in the Katz fibre. Since the constant section is overconvergent, fibrewise injectivity gives the same equality in the analytic fibre. External product with either regular modular input commutes with evaluation $A_U\to A_U/(\ell)$ after a common restriction of neighbourhoods. The Eisenstein section has canonical extension; the Siegel cusp ideal provides the subcanonical factor required for the trace.
\end{proof}

We now fix the generators and translated embeddings used in the improved functions, together with the Fourier convention for their depleted counterparts.

\begin{defn}\label{def:family-generators}
Choose bases $\underline\eta$ of $S^2(\underline\Pi)$ and $\underline\xi$ of $S^1(\underline\Pi)$, and retain the normalized long-Weyl generator $\underline\upsilon$ of Section~\ref{subsubsec:GL2-input}. Their good classical fibres give the refined eigenline bases used in the period comparison.
\end{defn}

For region (F), the translated embedding is defined using the matrices
\[
\tau = \begin{smatrix} 1 \\ 1&1 \\ &&\hfill 1 \\ &&\hfill -1&1 \end{smatrix},\qquad \hat \tau = \begin{smatrix} 1\\ 1&&1 \\ &\hfill -1\\ &\hfill 1&&1 \end{smatrix}.
\]
The matrix $\tau$ is denoted $\gamma$ in \cite{LPSZ1}.
\begin{itemize}
    \item We write $\hat \iota$ for the embedding of flag varieties sending $B_H h$ to $P_{\Si} h \hat \tau$, with the notations introduced in \S\ref{subsec:parabolic} (in particular, it is a closed immersion of $\mathbb Z_p$-schemes).
    \item We write $\gamma_{0,S}$ for the tame test data, as defined in \cite[\S10.2]{LZ-BK}.

\end{itemize}

We choose each tame translation $\gamma_{0,S}$ once and apply it to the whole family, pulling the resulting classes to a common tame level $\hat K^p$. This preserves the ordinary refinements and the prescribed archimedean vectors. The choices may differ between regions (F) and (D).

The construction of \cite{LPSZ1} imposes the character relation which, at classical weights, reads
\[
k_1-k_2+2=c_1+c_2.
\]
On this subspace, the family class $\underline\eta$ pairs with a product of two degree-zero elliptic classes through the Serre-duality pairing of \cite[Theorem 9.6.2]{LZ-BK}. This is the three-variable construction we use. The differential operators of \cite{GrahamRockwood} allow an additional variable, but we do not consider this approach here.

For region (D), the pairing instead uses a degree-one class on each side. With the conventions about flag varieties introduced in \cite[\S3]{LR3}, define the element $\tau^{\sharp} = \iota(w_{01})^{-1} \tau w_2 \in G(\mathbb Z_p)$. We will consider the translated embedding $\iota^\sharp: H \to G$ given by $h \mapsto \iota(h) \tau^\sharp$. The map $\FL_H \to \FL_G$ induced by $\iota^{\sharp}$ by construction sends $[w_{01}]$ to $[w_2]$. We also have projection maps $\pi_i: \FL_H \to \FL_{\GL_2} \cong \mathbb{P}^1$.

The constructions of \cite{LPSZ1,LR3} use Katz's ordering $\mathcal E(0,\mathbf t-1)$, whereas the $p$-depletion of our family is $\mathcal E(\mathbf t-1,0)$. To identify them, we Fourier-transform the Schwartz data. We use the symplectic Fourier transform $T\Phi(x,y)=\int\Phi(u,v)\psi(xv-yu)\,du\,dv$, with $\psi_\infty(x)=\exp(-2\pi i x)$ and self-dual measures. For $\Psi^{(p)}=T\Phi^{(p)}$, the identity $(T\Phi)'(u,v)=\Phi'(v,u)$ gives
\begin{equation}\label{eq:katz-transport}
\mathcal E_{\Phi^{(p)}}(\mathbf t-1,0)
=\mathcal E_{\Psi^{(p)}}(0,\mathbf t-1).
\end{equation}
By \cite[(8.1) and Remark 8.2]{LPSZ1}, this transport sends $s=1-t/2$ to $s=t/2$, fixes the archimedean weight-$t$ Gaussian, and gives the ordinary and depleted local Whittaker functions $p^{-mt/2}\mathbf1_{m\geq0}$ and $\mathbf1_{m=0}$ at $\operatorname{diag}(p^m,1)$. Thus they are the vectors $W_{\mathfrak a_1}$, with $\mathfrak a_1=p^{1-t/2}$, and $W_1^{\mathrm{dep}}$ used below, both normalized to have value $1$ at the identity.

In region (F) we use \cite[Definition 10.6.2]{LZ-BK}. In region (D) we use the pairing of Proposition~\ref{prop:supported-pairing}, normalized as in Proposition~\ref{prop:running-level}, so that its classical specializations are the pairings $B_2$ at level $p^2$ of \eqref{eq:classical-B2}; this normalization is the same for both Eisenstein inputs. The Eisenstein weight characters $\htt_F,\htt_D$ have classical values
\[
t_F=k_1-k_2-c+2,\qquad t_D=-k_1+k_2+c-2.
\]
\begin{defn}\label{def:joint-functions}
With the bases, Gauss sums and tame data fixed above, put
\begin{align*}
L_p^{(F)}
&=\big\langle\hat\iota^*(\gamma_{0,S}\underline\eta),
 \mathcal E(\htt_F-1,0)\boxtimes G(\chi_\Sigma^{-1})\mathcal G^{[p]}\big\rangle_{\mathrm{an}},\\
L_p^{(F),\imp}
&=\big\langle\hat\iota^*(\gamma_{0,S}\underline\eta),
 \mathcal E^\flat(\htt_F)\boxtimes G(\chi_\Sigma^{-1})\mathcal G\big\rangle_{\mathrm{an}},\\
L_p^{(D)}
&=\big\langle r_{D,A}(\gamma_{0,S}\underline\xi),
 \mathcal E(\htt_D-1,0)\boxtimes G(\chi_\Sigma^{-1})\underline\upsilon\big\rangle_{\mathrm{an}},\\
L_p^{(D),\imp}
&=\big\langle r_{D,A}(\gamma_{0,S}\underline\xi),
 \mathcal E^\flat(\htt_D)\boxtimes G(\chi_\Sigma^{-1})\underline\upsilon\big\rangle_{\mathrm{an}}.
\end{align*}
Here $\imp$ denotes improvement by the ordinary Eisenstein input; the superscript ``imp'' in \cite{LR3} instead denotes imprimitive tame data. In region (F), the bracket is $\nu_K$ times \eqref{p-adic:Serre}, with its branching map and translated pullback. In region (D), it is the normalized family pairing, including the same tame multiplier. The map $r_{D,A}$ is the translated pullback and branching map of Section~\ref{subsec:family-duality}, and the modular input includes the determinant transport \eqref{eq:determinant-family}. Section~\ref{subsubsec:modular-adjoint} compares this construction with the classical modular adjoint. We suppress the common tame Eisenstein subscript and work on the character components on which these maps are defined.
\end{defn}

These functions are defined on $\mathcal W_F^{\imp}$ and $\mathcal W_D^{\imp}$. By Lemma~\ref{lem:regularized-eisenstein}, the improved functions are meromorphic, with at most a simple pole along the divisor where the Eisenstein character is trivial.

\subsection{Classical specialization and comparison of pairings}\label{subsec:classical-comparisons}

Taking the Eisenstein residue leaves a constant in the first modular factor. We compare the resulting analytic pairing with coherent Serre duality, and then with the pairing $P_\diamond$ of Section~\ref{Section:Identification:Modular:Forms}. The comparison is made at a fixed auxiliary level; the passage to level $p^2$ is made in Section~\ref{subsec:level-comparison}.

\subsubsection{Classical inputs and coefficient maps}
\label{subsubsec:classical-specialization}
Let $Q_0:A\to K$ be a good classical specialization of the fixed families, with Siegel weight $(r_1,r_2)$ and modular weights $(0,c_0)$, where $m=r_1+r_2$ and $c_0=r_1-r_2+2$. Enlarge $K$ to contain its residue field, and put $\nu_0=Q_0\circ\nu_A=(r_1,r_2;m)$ and $\tau_0=Q_0\circ\tau_A$.

Coefficient specialization induces maps
\[
\begin{aligned}
\operatorname{sp}_{G,Q_0}:\;&
H^{i_\diamond}_{w_{G,\diamond},\mathrm{an}}
(\nu_A,\mathrm{cusp})^{-,\mathrm{fs}}
\otimes_{A,Q_0}K
\longrightarrow
H^{i_\diamond}_{w_{G,\diamond},\mathrm{an}}
(\nu_0,\mathrm{cusp})^{-,\mathrm{fs}},\\
\operatorname{sp}_{H,Q_0}:\;&
H^{2-i_\diamond}_{w_\diamond,\mathrm{an}}
(S_{H,\mathrm{Iw}(p^T)},\tau_A)^{+,\dagger}
\otimes_{A,Q_0}K
\longrightarrow
H^{2-i_\diamond}_{w_\diamond,\mathrm{an}}
(S_{H,\mathrm{Iw}(p^T)},\tau_0)^{+,\dagger}.
\end{aligned}
\]
We specialize the Siegel class on the line of Lemma~\ref{lem:ordinary-control}. For the modular input, we first specialize the two elliptic factors and then take their external product. Let $\eta_A$ denote the Siegel family class, namely $\eta$ in $(F)$ and $\xi$ in $(D)$, and let $h_A$ denote the modular input. Set
\[
\eta_{Q_0}:=\operatorname{sp}_{G,Q_0}(\eta_A\otimes1),
\qquad
h_{Q_0}:=\operatorname{sp}_{H,Q_0}(h_A\otimes1).
\]

To identify the modular coefficient in $(D)$, write $\rho_H$ for the half-sum of positive roots of $H$, $w_{0,M}$ for the longest Weyl element of the Siegel Levi, and $K^\lambda(z)=(1+z)^\lambda$ for the branching kernel.

\begin{lemma}\label{lem:D-coefficient-character}
Evaluation against the kernel $K^{\nu_{2,A}+1-a_A}$ gives the branching character
\[
\tau_A=(a_A-2,-b_A;m_A).
\]
At $Q_0$, the positive coefficient is $\omega_H(0,2-c_0;-m)$.
\end{lemma}

\begin{proof}
The defining Weyl formula of \cite[\S4.4]{LR3} gives the polynomial-model weight
\[
\kappa_{\mathrm{pol},A}
=(\nu_{1,A}+3,1-\nu_{2,A};-m_A).
\]
For $\lambda_A=\nu_{2,A}+1-a_A$, we have
\[
w_{0,M}\kappa_{\mathrm{pol},A}
+(\lambda_A,-\lambda_A;0)
=(2-a_A,b_A;-m_A).
\]
By \cite[Lemma 4.7]{LR3}, this is the character of the kernel; evaluation on distributions therefore has the opposite character $\tau_A$. The positive coefficient is
\[
V^H_{-\tau_A-2\rho_H}
=V^H_{(-a_A,b_A-2;-m_A)}.
\]
At $Q_0$, this becomes $V^H_{(0,c_0-2;-m)}=\omega_H(0,2-c_0;-m)$. Moreover, the specialized kernel has degree $r_2+1\le m+2$, so it belongs to the algebraic polynomial subrepresentation.
\end{proof}

In the Weyl-dependent convention of \cite[\S6.3]{boxerpilloni20}, the analytic weight is
\[
\sigma_A=(a_A-2,b_A-2;m_A),\qquad
\tau_A=w_{01}(\sigma_A+\rho_H)-\rho_H.
\]
This gives the same positive coefficient $V^H_{-\tau_A-2\rho_H}$, with the branching convention of \cite[Definition 4.10]{LR3}.

In the two regions we consequently have
\[
\tau_0=
\begin{cases}
(-2,c_0-2;m),&\diamond=F,\\
(-2,-c_0;m),&\diamond=D,
\end{cases}
\qquad
H_\diamond^{\mathrm{sc}}=
\begin{cases}
\omega_H(0,c_0;-m),&\diamond=F,\\
\omega_H(0,2-c_0;-m),&\diamond=D.
\end{cases}
\]
The coefficient identification in region $(F)$ is \cite[\S9.6]{LZ-BK}.

Choose modular levels at which the two elliptic inputs are defined, and enlarge $K$ so that the determinant components in \eqref{Union:Modular:Curves} are defined over $K$. Let $Y$ be the disjoint union of the product components used by the modular map, and write $X_1\times X_2$ for each component. The following formulas are understood componentwise. The second-factor coefficient $\mathcal F_{2,\diamond}$ is the automorphic line bundle associated with
\[
\chi_{2,\diamond}(\operatorname{diag}(x,y))=
\begin{cases}
x^{1-r_2}y^{-r_1-1},&\diamond=F,\\
x^{-r_1}y^{-r_2},&\diamond=D.
\end{cases}
\]
Thus the modular product coefficient on $Y$ is $\mathcal O_{X_1}\boxtimes\mathcal F_{2,\diamond}$. Its pullback, with the coefficient transport of Section~\ref{subsubsec:family-support}, is $H_\diamond^{\mathrm{sc}}$ on the common-level $H$-surface.

For sufficiently small rational $\epsilon>0$, put
\[
U_i=X_i[1-\epsilon,1],\qquad Z_2=X_2[0,\epsilon[,
\]
where these loci are pulled back from Iwahori level. The comparison on the second curve takes place between
\[
\mathcal N_\diamond=
\begin{cases}
H^0(U_2,\mathcal F_{2,F}^{\mathrm{an}})^{\mathrm{ord}},&\diamond=F,\\
H^1_{Z_2}(X_2^{\mathrm{an}},\mathcal F_{2,D}^{\mathrm{an}})^{\mathrm{ord}},&\diamond=D,
\end{cases}
\qquad
\mathcal C_\diamond=H^{2-i_\diamond}(X_{2,K},\mathcal F_{2,\diamond})^{\mathrm{ord}}.
\]
In these definitions we take the subspaces selected by the fixed cuspidal eigensystem and the prescribed old-vector line; $\mathrm{ord}$ denotes slope zero for the normalized elliptic operator with the indicated coefficient character. The degree-one input is constructed with coefficient $\mathcal F_{2,D}(-D_2)$, where $D_2$ is the cusp divisor; in $\mathcal N_D$ we use its image under $\mathcal F_{2,D}(-D_2)\to\mathcal F_{2,D}$. We refer to this map as forgetting the boundary twist.

The first factor needs no control theorem: the section $1\in H^0(X_1,\mathcal O_{X_1})$ restricts to $1$ on $U_1$. Its external product with $u\in\mathcal N_\diamond$ lies in
\[
\begin{cases}
H^0(U_1\times U_2,\mathcal O_{U_1}\boxtimes\mathcal F_{2,F}^{\mathrm{an}}),&\diamond=F,\\
H^1_{U_1\times Z_2}(U_1\times X_2^{\mathrm{an}},\mathcal O_{U_1}\boxtimes\mathcal F_{2,D}^{\mathrm{an}}),&\diamond=D.
\end{cases}
\]
Passing to the cofinal neighbourhoods and supports of Section~\ref{subsubsec:family-support} gives the positive realization used in \eqref{p-adic:Serre}; see \cite[Proposition 9.6.1]{LZ-BK} and \cite[\S5.4.1 and Theorems 5.4.12 and 6.4.8]{boxerpilloni20}. We form these products on the curves before applying the modular level maps.

Write $e_{Q_0}\boxtimes u_{Q_0}$ for the product representative of $h_{Q_0}$, with $u_{Q_0}\in\mathcal N_\diamond$. For the constant input $e_{Q_0}=1$; for the regularized Eisenstein input, Lemma~\ref{lem:regularized-eisenstein} gives $e_{Q_0}=-(1-p^{-1})1$. We retain this scalar when forming the classical input.

Set
\[
\kappa_\diamond=
\begin{cases}
(r_1,-r_2-2;m),&\diamond=F,\\
(r_2-1,-r_1-3;m),&\diamond=D.
\end{cases}
\]

On the Siegel side, restriction of distributions to algebraic polynomials defines the moment map $\operatorname{mom}_{G,Q_0}$, from weight $\nu_0$ to weight $\kappa_\diamond$, on the prescribed Weyl supports. Write $\operatorname{glob}_{G,Q_0}$ for the subsequent identification with global coherent cohomology on the selected ordinary eigenspace. This is the edge map of \cite[Theorems 5.5.1 and 5.12.3]{boxerpilloni20}: the other terms vanish by \cite[Theorem 1.4.10]{Boxer:Pilloni:Siegel} under the hypotheses of Lemma~\ref{lem:ordinary-control}. Both maps include the fixed level maps.

On the second modular curve, let $\operatorname{res}_2$ be restriction of global sections to $U_2$, and let $\operatorname{forg}_2$ be the map forgetting supports in $Z_2$, followed by proper comparison with coherent cohomology on $X_{2,K}$.

\begin{prop}
\label{prop:classical-specialization}
On the selected ordinary Siegel eigenspaces, the moment and global comparison maps induce isomorphisms
\[
\begin{aligned}
H^{i_\diamond}_{w_{G,\diamond},\mathrm{an}}
(\nu_0,\mathrm{cusp})^{-,\mathrm{ord}}
\xrightarrow{\operatorname{mom}_{G,Q_0}}
H^{i_\diamond}_{w_{G,\diamond}}
(\kappa_\diamond,\mathrm{cusp})^{-,\mathrm{ord}}
\xrightarrow{\operatorname{glob}_{G,Q_0}}
H^{i_\diamond}
(X_{G,K},[L_{i_\diamond}](-D_G))^{\mathrm{ord}}.
\end{aligned}
\]
Moreover, restriction and forgetting supports give isomorphisms
\[
\begin{cases}
\operatorname{res}_2:\mathcal C_F
\xrightarrow{\;\sim\;}\mathcal N_F,
&\diamond=F,\\
\operatorname{forg}_2:\mathcal N_D
\xrightarrow{\;\sim\;}\mathcal C_D,
&\diamond=D.
\end{cases}
\]
\end{prop}

\begin{proof}
The coefficient identifications are
\[
V^G_{(r_1,-r_2-2;m)}=[L_2],
\qquad
V^G_{(r_2-1,-r_1-3;m)}=[L_1];
\]
see \cite[\S2.5.3 and Lemma 9.2.2]{LZ-BK} and \cite[\S4.5]{LR3}. The Siegel isomorphisms follow from \cite[Theorem 5.12.3 and Corollary 6.8.4]{boxerpilloni20}, strengthened by \cite[Theorem 1.4.10]{Boxer:Pilloni:Siegel}. Here $\mathrm{ord}$ denotes simultaneous slope zero for the normalized Siegel and Klingen operators.

For an algebraic elliptic character $\kappa=(k'_1,k'_2)$, put $k=k'_1-k'_2$. By \cite[Remark 3.2.14 and Definition 3.2.16]{RC:EichlerShimura}, the normalized classical operator is
\[
U_p^{\mathrm{good}}
=p^{-\min\{1-k'_1,-k'_2\}}p^{-k'_1}U_{p,k}^{\mathrm{naive}}
=p^{-\min\{1,k\}}U_{p,k}^{\mathrm{naive}}.
\]
This is the normalization of \cite[Theorem 5.13]{BP:Modular}. A determinant twist shifts both entries of $\kappa$ by the same integer, so the normalization depends only on the Hodge weight $k$. For $\mathcal F_{2,F}$ and $\mathcal F_{2,D}$ these weights are $c_0$ and $2-c_0$, respectively. Parts (2) and (3) of that theorem therefore apply to restriction in degree zero and forgetting supports in degree one, with respective bounds $<c_0-1$ and $<1-(2-c_0)=c_0-1$. Since $c_0\geq2$, both bounds contain the selected slope-zero lines. At deeper levels we use the positive finite-slope pullbacks of \cite[Theorems 5.4.14 and 6.4.10]{boxerpilloni20}; the prescribed old vectors and their classical realizations are transported by the same level maps.
\end{proof}

Put
\[
c_{2,\diamond}:=
\begin{cases}
\operatorname{res}_2^{-1},&\diamond=F,\\
\operatorname{forg}_2,&\diamond=D.
\end{cases}
\]
The comparison of modular inputs is obtained by applying $c_{2,\diamond}$ to the second factor and leaving the constant factor unchanged. On the product of curves this gives
\begin{equation}\label{eq:modular-product-comparison}
\begin{aligned}
\left.(1\boxtimes c_{2,F}(u))\right|_{U_1\times U_2}
&=1|_{U_1}\boxtimes u,\\
\operatorname{forg}_{U_1\times Z_2}(1|_{U_1}\boxtimes u)
&=\left.(1\boxtimes c_{2,D}(u))\right|_{U_1\times X_2^{\mathrm{an}}}.
\end{aligned}
\end{equation}
These equalities follow by applying restriction or forgetting supports to the second factor of the external product. They remain valid after forgetting the boundary twist, since that operation is induced by a morphism of coefficient sheaves.

To place these classes on the surface where the pairing is evaluated, let $\pi_{\diamond,T}:X_{H,K}\to Y$ be the modular map in the fixed-level construction. In $(F)$ it is induced by the two projections of $H$; in $(D)$ its first factor includes the translation $d_T=\operatorname{diag}(p^T,1)$ of Section~\ref{subsubsec:family-support}. The source is the common level at which this map and the translated Siegel embedding are both defined. We keep all determinant components, and write
\[
\theta_{\diamond,T}:\pi_{\diamond,T}^*(\mathcal O_{X_1}\boxtimes\mathcal F_{2,\diamond})\longrightarrow H_\diamond^{\mathrm{sc}}
\]
for the coefficient map of this modular pullback. In $(D)$ it is the normalized map on the first coefficient tensored with the identity on $\mathcal F_{2,D}$, as in Section~\ref{subsubsec:family-support}. At weight zero it fixes the constant section on the first curve; at other algebraic weights we use the specialization of the same family map. Any prescribed finite-character projector is applied on both sides of \eqref{eq:modular-product-comparison}. Its deck transformations commute with restriction, forgetting supports and the coefficient maps, so the comparison is preserved.

\begin{defn}\label{def:classical-inputs}
For the chosen family inputs, their \emph{classical realizations} are
\[
\begin{aligned}
\eta_{Q_0}^{\mathrm{cl}}&:=\operatorname{glob}_{G,Q_0}\bigl(\operatorname{mom}_{G,Q_0}(\eta_{Q_0})\bigr),\\
h_{Q_0}^{\mathrm{cl}}&:=\theta_{\diamond,T}\pi_{\diamond,T}^*\bigl(e_{Q_0}\boxtimes c_{2,\diamond}(u_{Q_0})\bigr).
\end{aligned}
\]
The second formula uses the specified elliptic input $u_{Q_0}$, with its fixed translates and, when present, the same finite-character projector as on the analytic side.
\end{defn}

Thus
\[
\eta_{Q_0}^{\mathrm{cl}}\in
\begin{cases}
H^2(X_{G,K},[L_2](-D_G)),&\diamond=F,\\
H^1(X_{G,K},[L_1](-D_G)),&\diamond=D,
\end{cases}
\qquad
h_{Q_0}^{\mathrm{cl}}\in
\begin{cases}
H^0(X_{H,K},\omega_H(0,c_0;-m)),&\diamond=F,\\
H^1(X_{H,K},\omega_H(0,2-c_0;-m)),&\diamond=D.
\end{cases}
\]

\begin{rk}\label{rk:selected-modular-control}
Since the elliptic classes are compared on the modular curves before they are inserted into the cohomology of $H$, we only need control for the selected elliptic lines, and not for the analytic cohomology of $X_H$. The same construction applies to a depleted input $g^{[p]}=(1-V_pU_p)g$: we apply it to both terms at a common level, using that restriction and coefficient specialization commute with the isogeny maps.
\end{rk}
\begin{rk}
After restricting the Siegel class and contracting its coefficients, the two inputs have coefficients $(H_\diamond^{\mathrm{sc}})^\vee\otimes K_{X_H}$ and $H_\diamond^{\mathrm{sc}}$, and degrees $i_\diamond$ and $2-i_\diamond$. Their cup product therefore admits Serre trace. The coefficient contraction is described below.
\end{rk}
\begin{lemma}\label{lem:elliptic-supported-control}
At a good crystalline specialization of weights $(r_1,r_2;a,b)$, with $b\geq2$ and $b-a=c_0$, the long-Weyl family class $\underline\upsilon$ of Section~\ref{subsubsec:GL2-input} specializes to the classical degree-one class dual to $g^\vee$ under the modified Serre pairing, with value $P(g^\vee,\underline\upsilon)=1$. After the determinant transport of \eqref{eq:determinant-family} and forgetting the cuspidal boundary twist, its external product with the first $\GL_2$ input has coefficient $\omega_H(a,2-b;-m)$.
\end{lemma}
\begin{proof}
Apply \cite[Theorem 5.13(3)]{BP:Modular} with $k=2-b$ to the supported cuspidal modular cohomology. The ordinary slope is less than $1-k=b-1$, so forgetting supports identifies the selected line with its classical realization. The modified pairing of \cite[\S5.7.3]{BP:Modular} is the normalized Atkin--Lehner transport followed by cup product and trace. Its specialization is the same pairing on the classical classes, so the chosen generators still pair to $1$. Since the family class already lies in the positive long-Weyl model, the Atkin--Lehner transport is not applied again. The determinant transport is \eqref{eq:determinant-family}, and the map forgetting the boundary twist commutes with forgetting supports. At $(a,b)=(0,c_0)$ this gives the degree-one input used in Definition~\ref{def:classical-inputs}.
\end{proof}

\subsubsection{The algebraic scalar pairing}\label{subsubsec:scalar-pairing}
We describe the coefficient contraction before comparing the traces. For integral weights in the ranges of Definition~\ref{Definition:Vector:Bundles}, the scalar coefficients are
\begin{equation}\label{eq:scalar-coefficients}
\begin{array}{c|c|c|c}
\diamond&H_\diamond^{\mathrm{sc}}&\mathcal Q_\diamond&\text{weight relation}\\\hline
F&\omega_H(a,b;-m)&\omega_H(2-a,2-b;m)&a+b=c_0\\
D&\omega_H(a,2-b;-m)&\omega_H(2-a,b;m)&b-a=c_0.
\end{array}
\end{equation}
In particular, at $(a,b)=(0,c_0)$ the quotients are $\mathcal Q_F=\omega_H(2,2-c_0;m)$ and $\mathcal Q_D=\omega_H(2,c_0;m)$. The branching law and the boundary map \eqref{Pullback:Vector:Bundles} give
\begin{equation}\label{eq:scalar-quotient}
q_\diamond^{\mathrm{sc}}:\iota^*([L_{i_{\diamond}}](-D_G))
\longrightarrow\mathcal Q_\diamond(-D_H)
=(H_\diamond^{\mathrm{sc}})^\vee\otimes K_{X_H}.
\end{equation}
The projection selects the indices $j_F=r_2+a$ and $j_D=r_1+1+a$ of \eqref{Table:Images:p}, with the character twist removed from the Siegel coefficient. For translated embeddings, the projection and its frame are transported by the same torsor map.

\begin{defn}\label{def:scalar-algebraic-pairing}
At the chosen common level, set
\begin{equation}\label{eq:scalar-algebraic-pairing}
\begin{aligned}
B_\diamond^{\mathrm{sc}}:
H^{i_\diamond}(X_G,[L_{i_{\diamond}}](-D_G))
\times H^{2-i_\diamond}(X_H,H_\diamond^{\mathrm{sc}})&\longrightarrow K,\\
(\eta,h)&\longmapsto\Tr_{X_H}(q_\diamond^{\mathrm{sc}}(\iota^*\eta)\cup h).
\end{aligned}
\end{equation}
Here $\Tr_{X_H}$ is the geometric Serre trace, as in Definition~\ref{Pairing:With:Projection}. We use the same translated maps and component sums on the analytic and algebraic sides, with the Siegel factor first. The tame volume multiplier in the functions of Definition~\ref{def:joint-functions} is kept separate.
\end{defn}

\begin{lemma}\label{lem:boundary-branching}
For $r_1\geq r_2\geq0$, the analytic branching quotient specializes to $q_\diamond^{\mathrm{sc}}$ in the ranges
\[
a+b=c_0,\quad-r_2\leq a\leq r_1+2\quad(F),\qquad
b-a=c_0,\quad-r_1-1\leq a\leq r_2+1\quad(D).
\]
Here the quotient generators are normalized as in the proof. In particular, the statement applies at $a=0$.
\end{lemma}
\begin{proof}
Put $N=m+2$. The dual polynomial models have weights $(r_2+2,-r_1;-m)$ and $(r_1+3,1-r_2;-m)$. The translated torus acts on $K_\lambda(z)=(1+z)^\lambda$ with weight $(B+\lambda,A-\lambda;C)$ for a polynomial model of highest weight $(A,B;C)$; see \cite[Lemma 8.3.2]{LZ-BK} and \cite[Lemma 4.7]{LR3}. Evaluation on distributions has the opposite character. Matching it with \eqref{eq:scalar-coefficients} gives
\[
\lambda_F=r_1+2-a,\qquad\lambda_D=r_2+1-a.
\]
The stated ranges are exactly $0\leq\lambda_\diamond\leq N$, so this evaluation factors through the algebraic moment quotient.

To fix its scale, use the monomial basis $1,z,\ldots,z^N$ in the polynomial model and its coordinate-dual basis in the distribution quotient. The Levi part of the translated torsor map is $\gamma=\left(\begin{smallmatrix}1&0\\1&1\end{smallmatrix}\right)$, and its inverse sends $z^\lambda$ to $(1+z)^\lambda$ in the induced-model convention of the cited lemmas. Consequently both the transported algebraic projection and the analytic evaluation send a moment functional $\mu$ to
\[
\mu((1+z)^\lambda)=\sum_{j=0}^{\lambda}\binom{\lambda}{j}\mu(z^j).
\]
We identify the selected quotient generator with $1$ through this coordinate projection and use its tensor transport for the projection defining $P_\diamond$. This fixes the generator allowed in Lemma~\ref{Explicit:Description:Q}; since we use coordinate-dual monomials, no further binomial factor appears. The identification descends from torsor charts to the associated sheaves, and both constructions use the same boundary morphism.
\end{proof}

\subsubsection{Comparison of the pairings}\label{subsubsec:pairing-comparisons}
We can now compare the analytic pairing with the scalar algebraic pairing at the same level.

\begin{theorem}\label{thm:analytic-algebraic}
Let $Q$ be a good algebraic specialization in the range of Lemma~\ref{lem:boundary-branching}. Suppose the chosen modular inputs have classical realizations under the comparison maps above; at $(a,b)=(0,c_0)$ these are the inputs of Proposition~\ref{prop:classical-specialization}. Let $\eta_A,h_A$ be the selected family classes, and denote their classical realizations by $\eta_Q^{\mathrm{cl}},h_Q^{\mathrm{cl}}$, as in Definition~\ref{def:classical-inputs} at $Q_0$. Then
\begin{equation}\label{eq:analytic-algebraic}
B_{\diamond,Q}(r_{\diamond,A}\eta_A\otimes1,h_A\otimes1)
=B_\diamond^{\mathrm{sc}}(\eta_Q^{\mathrm{cl}},h_Q^{\mathrm{cl}}).
\end{equation}
The level, finite translates and coefficient generators on the two sides are the same, and both pairings use geometric Serre trace. At an Eisenstein pole the assertion applies to the regularized input: at weight zero its first factor is $-(1-p^{-1})\mathbf1$. It also applies to prescribed depleted classical inputs, using the finite isogeny maps of Remark~\ref{rk:selected-modular-control}.
\end{theorem}
\begin{proof}
By Proposition~\ref{prop:supported-pairing}, specialization commutes with the pullbacks, coefficient maps and cup product defining the analytic pairing. The Siegel comparison of Proposition~\ref{prop:classical-specialization} applies at the Siegel weight of $Q$, under the ordinary-control hypotheses of Lemma~\ref{lem:ordinary-control}. Lemma~\ref{lem:boundary-branching} identifies the specialized coefficient quotient with $q_\diamond^{\mathrm{sc}}$.

For the modular input at $Q_0$, use its product representative $e_{Q_0}\boxtimes u_{Q_0}$. In $(F)$, \eqref{eq:modular-product-comparison} identifies it with the restriction of $e_{Q_0}\boxtimes c_{2,F}(u_{Q_0})$. In $(D)$ the same equation identifies its image after forgetting support in the second factor with the restriction of $e_{Q_0}\boxtimes c_{2,D}(u_{Q_0})$. This is the positive supported realization of the degree-one class; its normalization is fixed by Lemma~\ref{lem:elliptic-supported-control}. At the other algebraic points of the statement, use the assumed classical realizations of the two modular inputs.

The cup-product pairings of \cite[Theorem 9.6.2]{LZ-BK} and \cite[Theorem 4.12]{LR3} are compatible with classical Serre duality under restriction and forgetting supports. Applied to the preceding representatives, their coefficient evaluation is
\[
\bigl((H_\diamond^{\mathrm{sc}})^\vee\otimes K_{X_H}\bigr)\otimes H_\diamond^{\mathrm{sc}}\longrightarrow K_{X_H}.
\]
On the common-level surface, pull back the modular representative by $\pi_{\diamond,T}$ and apply $\theta_{\diamond,T}$, as in Definition~\ref{def:classical-inputs}. Pullback commutes with cup product and with the maps forgetting supports, and the coefficient map is the specialization of the one used analytically. Thus the supported product, after excision and forgetting supports, represents
\[
q_\diamond^{\mathrm{sc}}(\iota^*\eta_Q^{\mathrm{cl}})\cup h_Q^{\mathrm{cl}}\in H^2(X_H,K_{X_H}).
\]
Equivalently, when the analytic pairing is written on the product $Y$, let $(\pi_{\diamond,T})_!$ be the adjoint on the negative factor: apply the dual of $\theta_{\diamond,T}$ and then the relative Serre trace. For a negative class $\alpha$ on the common-level surface and a modular class $v$ on $Y$, the projection formula reads
\[
\Tr_Y\bigl((\pi_{\diamond,T})_!\alpha\cup v\bigr)
=\Tr_{X_H}\bigl(\alpha\cup\theta_{\diamond,T}\pi_{\diamond,T}^*v\bigr).
\]
Both traces are geometric, that is, not divided by the degree, and finite-character projectors are treated term by term. The Siegel factor remains first throughout, so no interchange of cohomology classes introduces a sign. This proves \eqref{eq:analytic-algebraic}.

For the regularized input, Lemma~\ref{lem:regularized-eisenstein} supplies the scalar $-(1-p^{-1})$ in the first factor on both sides. For depleted inputs, apply the comparison to the finite isogeny maps of Remark~\ref{rk:selected-modular-control} and use linearity.
\end{proof}

\begin{defn}\label{def:tensor-transports}
Let $J\in\{G,H\}$, and let $\mathcal V$ be an automorphic vector bundle with its canonical or subcanonical extension on $X_J$. Using the trivializations fixed in Remark~\ref{Remark:Trivialization}, define, for $s\in\mathbb Z$,
\[
T_s:H^\bullet(X_J,\mathcal V)\longrightarrow H^\bullet(X_J,\mathcal V\otimes C_{2s}),
\qquad v\longmapsto v\otimes e_s.
\]
\end{defn}

Put
\[
(d_F,h_F):=(r_2+1,r_1-1),
\qquad
(d_D,h_D):=(r_1+2,r_2-2).
\]
In either case, we have $d_\diamond+h_\diamond=m$.

\begin{theorem}\label{thm:scalar-period-comparison}\label{lem:scalar-pairing-comparison}
Fix integral weights in the range of Lemma~\ref{lem:boundary-branching}, and use its quotient generators. For
\[
\eta\in H^{i_\diamond}(X_G,[L_{i_{\diamond}}](-D_G)),\qquad
h\in H^{2-i_\diamond}(X_H,H_\diamond^{\mathrm{sc}}),
\]
at the same level and with the same finite translates, one has
\begin{equation}\label{eq:boundary-branching}
P_\diamond(T_{d_\diamond}\eta,T_{h_\diamond}h)
=B_\diamond^{\mathrm{sc}}(\eta,h)e_m.
\end{equation}
Both sides use geometric Serre trace; over $K$, the target is the scalar extension of $E\{m\}$.
\end{theorem}
\begin{proof}
Definition~\ref{Definition:Vector:Bundles} and \eqref{eq:scalar-coefficients} give
\begin{equation}\label{eq:coefficient-conventions}
E_\diamond=[L_{i_{\diamond}}]\otimes C_{2d_\diamond},\qquad
H_\diamond=H_\diamond^{\mathrm{sc}}\otimes C_{2h_\diamond}.
\end{equation}
Indeed, the central coordinates of $H_F,H_D$ are respectively $c_0-4=-m+2(r_1-1)$ and $-c_0-2=-m+2(r_2-2)$. Their noncentral coordinates agree with \eqref{eq:scalar-coefficients}. The projection defining $P_\diamond$ is $q_\diamond^{\mathrm{sc}}$ tensored with $C_{2d_\diamond}$; the chosen quotient generators are related by the same transport. Restriction and cup product therefore give the scalar contraction tensored with $e_{d_\diamond}\otimes e_{h_\diamond}=e_m$. Taking the same trace proves \eqref{eq:boundary-branching}.
\end{proof}

\begin{rk}\label{rk:classical-complex-period}
Suppose the classes are defined over a number field $E$ and $(a,b)=(0,c_0)$. Put $\eta^\sharp=T_{d_\diamond}\eta$ and $h^\sharp=T_{h_\diamond}h$. If these satisfy the hypotheses of Theorem~\ref{Integral:New:Pairing}, with $h^\sharp=1\boxtimes g$, then, in the $e_m$-coordinate,
\[
B_\diamond^{\mathrm{sc}}(\eta,h)\overset{.}{=}\mathcal I(\Psi_{\eta^\sharp},\varphi_g).
\]
Here the left-hand side is evaluated under the fixed embedding $E\hookrightarrow\mathbf C$, and the automorphic vectors and proportionality convention are those of that theorem. Section~\ref{subsec:poles-comparison} applies this formula to the residues of the improved functions.
\end{rk}

\begin{rk}\label{rk:tame-trace-normalization}
The functions of Definition~\ref{def:joint-functions} use $\nu_K$ times geometric trace. Their algebraic specialization is therefore $\nu_K B_\diamond^{\mathrm{sc}}$, after the level normalization of Section~\ref{subsec:level-comparison}; tensor transport gives $\nu_K P_\diamond$. This multiplier is the same for the ordinary and depleted inputs.
\end{rk}

\begin{rk}\label{rk:twisted-Hecke}
Tensor transport by $T_{d_\diamond}$ is not equivariant for the raw Hecke correspondences. If $t$ has similitude $p^j$, Remark~\ref{Remark:Trivialization} gives
\begin{equation}\label{eq:twisted-U}
[t]_{E_\diamond}T_{d_\diamond}
=p^{jd_\diamond}T_{d_\diamond}[t]_{[L_{i_{\diamond}}]},\qquad
\mathcal U'_{t,E_\diamond}
=T_{d_\diamond}\mathcal U'_{t,\mathrm{sc}}T_{d_\diamond}^{-1}.
\end{equation}
Equivalently, the scalar normalization of the raw operator on $E_\diamond$ is multiplied by $p^{-jd_\diamond}$; the transported normalized operators have the same eigenvalues $\alpha$ and $\alpha\beta^\circ$. These are formulas at algebraic weights; the analytic family operators are defined by the normalized torsor and trace maps.
\end{rk}

\subsection{Local integrals and test vectors}\label{subsec:local-improvement}

The local integrals at $p$ have two uses: they determine the Euler factor removed by improvement and show that the prescribed ordinary vectors detect the distinguished period at weight zero. We follow \cite{LZ-zeta} in $(F)$ and \cite[\S7]{LR2} in $(D)$.

We begin with the local conventions of these references. Let $\pi$ be an irreducible generic unramified representation of $G(\mathbb Q_p)$, and let $\sigma_1,\sigma_2$ be unramified principal series of $\GL_2(\mathbb Q_p)$, with
\[
\omega_\pi\omega_{\sigma_1}\omega_{\sigma_2}=1.
\]
The first $\GL_2$ factor will be the Eisenstein input. We first make the calculations for generic inducing characters and then specialize the resulting formulas wherever the normalized vectors are defined. Write $\alpha,\beta,\gamma,\delta$ for the Hecke parameters of $\pi$, and $\mathfrak a_i,\mathfrak b_i$ for those of $\sigma_i$. In the local normalization, these are the Satake parameters of $\pi\otimes|\cdot|^{-3/2}$ and $\sigma_i\otimes|\cdot|^{-1/2}$, respectively. Thus
\[
\alpha\delta=\beta\gamma=p^3\omega_\pi(p),
\qquad
\mathfrak a_i\mathfrak b_i=p\omega_{\sigma_i}(p),
\]
and the relation between central characters gives
\begin{equation}\label{eq:local-central-characters}
\beta\gamma\mathfrak a_1\mathfrak b_1\mathfrak a_2\mathfrak b_2=p^5.
\end{equation}
We return to the motivic normalization of Section~\ref{sec:review} at the end of the subsection, where we give the required substitutions.

Let $\mathfrak z$ be the $H(\mathbb Q_p)$-invariant trilinear functional normalized to take the value $1$ on the normalized spherical vectors. Its existence and uniqueness follow from \cite[Theorem 2.3]{LZ-zeta}. We use an unramified Whittaker character $\psi$ for $\pi$ and $\psi^{-1}$ for the two $\GL_2$ factors, and normalize $d^\times x$ so that $\mathbb Z_p^\times$ has volume $1$.

For each $i$, let $W_{\mathfrak a_i}$ be the Iwahori eigenvector with eigenvalue $\mathfrak a_i$, normalized by $W_{\mathfrak a_i}(1)=1$, and let $W_i^{\dep}$ be its $p$-depletion. The torus restrictions of these vectors are
\[
W_{\mathfrak a_i}\!\begin{pmatrix}p^m&0\\0&1\end{pmatrix}
=\begin{cases}(\mathfrak a_i/p)^m&m\geq0,\\0&m<0,\end{cases}
\qquad
W_i^{\dep}\!\begin{pmatrix}x&0\\0&1\end{pmatrix}
=\mathbf1_{\mathbb Z_p^\times}(x).
\]
These are the vectors of \cite[\S3.3 and Definition 5.14]{LZ-zeta}. We also put $W'_n=W_{\alpha,\beta}^{\prime,\Iw}[n]$, using the trace-compatible Iwahori vectors of loc.\,cit.

The degree-eight Euler polynomials for the depleted data are
\begin{align}
\mathcal B_F
&=\prod_{x\in\{\alpha,\beta\}}
  \prod_{y\in\{\mathfrak a_1,\mathfrak b_1\}}
  \prod_{v\in\{\mathfrak a_2,\mathfrak b_2\}}
  \left(1-\frac{p^2}{xyv}\right),\label{eq:local-euler-F}\\
\mathcal B_D
&=\prod_{y\in\{\mathfrak a_1,\mathfrak b_1\}}
  \prod_{v\in\{\mathfrak a_2,\mathfrak b_2\}}
  \left(1-\frac{p^2}{\alpha yv}\right)
  \prod_{x\in\{\beta,\gamma\}}
  \prod_{y\in\{\mathfrak a_1,\mathfrak b_1\}}
  \left(1-\frac{p^2}{xy\mathfrak a_2}\right).\label{eq:local-euler-D}
\end{align}
We will show that the factors removed by the improved constructions are
\begin{equation}\label{eq:local-removed-factors}
q_F=1-\frac{p^2}{\beta\mathfrak b_1\mathfrak b_2},
\qquad
q_D=1-\frac{p^2}{\gamma\mathfrak b_1\mathfrak a_2}.
\end{equation}
Put
\[
\mathcal B_F^{\imp}=\frac{\mathcal B_F}{q_F},\qquad
\mathcal B_D^{\imp}=\frac{\mathcal B_D}{q_D},\qquad
C_p=\frac{p^3}{(p+1)^2(p-1)}.
\]
Here division means deleting the indicated factor from the product. Thus $\mathcal B_F^{\imp}$ and $\mathcal B_D^{\imp}$ are degree-seven Euler polynomials, defined also where $q_F$ or $q_D$ vanishes.

\begin{prop}\label{prop:local-F}
Let $u_{\Iw}$ be the element of \cite[Definition 5.12]{LZ-zeta}, denoted $\hat\tau$ in our construction. For $n\geq2$,
\begin{align*}
\mathfrak z(u_{\Iw}W'_n,W_1^{\dep},W_2^{\dep})
&=\mathfrak z(u_{\Iw}W'_n,W_1^{\dep},W_{\mathfrak a_2})
=C_p\mathcal B_F,\\
\mathfrak z(u_{\Iw}W'_n,W_{\mathfrak a_1},W_{\mathfrak a_2})
&=C_p\mathcal B_F^{\imp}.
\end{align*}
\end{prop}
\begin{proof}
The first line is \cite[Theorem 5.15]{LZ-zeta}, and the second is Proposition 5.13 of the same reference. The polynomial $\mathcal B_{\Kl}$ appearing there is \eqref{eq:local-euler-F}, and both formulas use the same normalization of $\mathfrak z$ and $W'_n$.
\end{proof}

In (F), replacing the depleted cuspidal input by its ordinary stabilization leaves the value unchanged, while replacing the Eisenstein input removes $q_F$. These are the two changes made in Definition~\ref{def:joint-functions}.

In $(D)$, \cite[\S7]{LR2} reduces the calculation first to the Siegel Jacquet module and then to a $\GL_2$ torus integral. Write
\[
\pi=\Ind_{P_{\Si}}^G(\rho\times\lambda),\qquad
\sigma_0=\rho\otimes\lambda,
\]
with normalized induction and $\lambda(p)=p^{-3/2}\alpha$. The Hecke parameters of $\sigma_0$ are then $\beta/p$ and $\gamma/p$. Let $\pi[\lambda]$ be the preimage of the canonical subrepresentation $(\rho\times\lambda)\delta_{P_{\Si}}^{-1/2}$ in the unnormalized Siegel Jacquet module. By \cite[Proposition 7.7]{LR2}, for $x\in\pi[\lambda]$ the sequence
\[
\left(\frac{\lambda(p)}{p^{3/2}}\right)^{-n}
\mathfrak z\bigl(\hat\tau\diag(1,1,p^n,p^n)x,y_1,\diag(-1,1)y_2\bigr)
\]
is eventually constant. Its eventual value factors through the Jacquet module and defines
\begin{equation}\label{eq:homsetmap}
\partial_{\Si}:\Hom_H(\pi\otimes\sigma_1\otimes\sigma_2,\mathbb C)
\longrightarrow\Hom_{\GL_2}(\sigma_0\otimes\sigma_1\otimes\sigma_2,\mathbb C).
\end{equation}
Here $\GL_2$ is embedded in the Siegel Levi by $A\mapsto(A,\mathrm{det} A)$. If $\mathfrak y$ denotes the spherical trilinear functional on $\sigma_0\otimes\sigma_1\otimes\sigma_2$, Proposition 7.10 of loc.\,cit. gives
\begin{align*}
\partial_{\Si}(\mathfrak z)&=\Delta\mathfrak y,\;\;
\Delta=\frac{p^2}{p^2-1}
L(\lambda^\vee\times\sigma_1^\vee\times\sigma_2^\vee,\tfrac12)^{-1}
=\frac{p^2}{p^2-1}
\prod_{y\in\{\mathfrak a_1,\mathfrak b_1\}}
\prod_{v\in\{\mathfrak a_2,\mathfrak b_2\}}
\left(1-\frac{p^2}{\alpha yv}\right).
\end{align*}
This supplies the first four factors of $\mathcal B_D$.

For the second reduction, take $W_0\in\mathcal W(\sigma_0)$ and $W_1\in\mathcal W(\sigma_1)$ fixed by the lower unipotent subgroup with lower-left entry in $p^r\mathbb Z_p$, for some $r\geq1$. By \cite[Proposition 7.11]{LR2}, the value
\[
\mathfrak y(W_0,W_1,W'_{\mathfrak a_2}[n])
\]
is independent of $n\geq r$. It is obtained by evaluating at $u=1/2$ the normalized integral
\begin{equation}\label{eq:second-local-reduction}
\frac{p}{p+1}L(\sigma_0\times\sigma_1\times\tau_2,u)^{-1}
\int_{\mathbb Q_p^\times}
W_0\!\begin{pmatrix}x&0\\0&1\end{pmatrix}
W_1\!\begin{pmatrix}x&0\\0&1\end{pmatrix}
\tau_2(x)|x|^{u-1}\,d^\times x,
\end{equation}
continued from its domain of convergence. Here $u$ is an auxiliary complex variable and $\tau_2$ is the unramified character with $\tau_2(p)=p^{-1/2}\mathfrak b_2$. At $u=1/2$, the reciprocal $L$-factor is
\[
\prod_{x\in\{\beta,\gamma\}}
\prod_{y\in\{\mathfrak a_1,\mathfrak b_1\}}
\left(1-\frac{p^2}{xy\mathfrak a_2}\right),
\]
which supplies the remaining four factors of $\mathcal B_D$. Thus the two reductions give $C_p\mathcal B_D$ multiplied by the torus integral.

\begin{prop}\label{prop:local-D}
Let $\tau^\sharp=w_{01}^{-1}\tau w_2$ be the translate of \cite[Proposition 7.13]{LR2}, and put
\[
d_n=\diag(p^n,1),\qquad
A_n=\left(\frac{p^2}{\beta\mathfrak b_2}\right)^n.
\]
For $n\geq2$,
\begin{align*}
\mathfrak z(\tau^\sharp W'_n,d_nW_1^{\dep},W_{\mathfrak a_2})
=A_nC_p\mathcal B_D,\quad
\mathfrak z(\tau^\sharp W'_n,d_nW_{\mathfrak a_1},W_{\mathfrak a_2})
=A_nC_p\mathcal B_D^{\imp}.
\end{align*}
\end{prop}
\begin{proof}
The first formula is \cite[Proposition 7.13]{LR2}, with the level exponent written as $n$, as in its proof. The second is the modification described in Remark 7.14 there, which we make explicit.

The first Jacquet reduction and the translations in that proof give the common factor $A_n\Delta$. After cancelling the simultaneous translations by $d_n$ using $\GL_2$-invariance, it remains to compare \eqref{eq:second-local-reduction} for
\[
W_0=W_{\beta/p},\qquad
W_1=W_1^{\dep}\quad\text{or}\quad W_{\mathfrak a_1}.
\]
The first vector satisfies $W_0(\diag(p^m,1))=(\beta/p^2)^m$ for $m\geq0$, and vanishes for $m<0$.

If $W_1=W_1^{\dep}$, only $x\in\mathbb Z_p^\times$ contributes, so the torus integral is $1$. If $W_1=W_{\mathfrak a_1}$, the integral becomes the geometric series
\[
\sum_{m\geq0}\left(\beta\mathfrak a_1\mathfrak b_2p^{-5/2-u}\right)^m
=\left(1-\beta\mathfrak a_1\mathfrak b_2p^{-5/2-u}\right)^{-1},
\]
initially for $\Re(u)$ sufficiently large. At $u=1/2$, relation~\eqref{eq:local-central-characters} gives
\[
\frac{\beta\mathfrak a_1\mathfrak b_2}{p^3}
=\frac{p^2}{\gamma\mathfrak b_1\mathfrak a_2}.
\]
The ordinary vector therefore contributes the denominator $q_D$, while the other factors are unchanged. Cancelling it against $\mathcal B_D$ gives the second formula, including at specializations where $q_D=0$.
\end{proof}

Put $w_G=k_1+k_2-3$. In the motivic normalization used elsewhere in the paper,
\[
\beta\gamma=\chi_\Pi(p)p^{w_G},\qquad
\mathfrak a\mathfrak b=\chi_\Sigma(p)p^{c-1},
\]
where $\mathfrak a,\mathfrak b$ are the parameters of the cuspidal $\GL_2$ form. We use the ordinary Eisenstein refinement with arithmetic root $1$, as in Section~\ref{sec:review}; this is the refinement on the untwisted Eisenstein component relevant to the pole. Its other arithmetic root at weight $t$ is $(\chi_\Pi\chi_\Sigma)(p)^{-1}p^{t-1}$.

The following table gives the substitutions at crystalline points. The left column refers to the local parameters used above; the right column is written in terms of the motivic parameters.
\[
\begin{array}{c|c}
\text{Local parameters}&\text{Substitution in motivic parameters}\\ \hline
(\alpha,\beta,\gamma,\delta)
&p^{(3-w_G)/2}(\alpha,\beta,\gamma,\delta)\\[2pt]
(\mathfrak a_1,\mathfrak b_1)
&\bigl(p^{1-t/2},(\chi_\Pi\chi_\Sigma)(p)^{-1}p^{t/2}\bigr)\\[2pt]
(\mathfrak a_2,\mathfrak b_2)
&p^{(2-c)/2}(\mathfrak a,\mathfrak b).
\end{array}
\]
From now on, $\alpha,\beta,\gamma,\delta$ again denote motivic parameters. Substituting $t=k_1-k_2-c+2$ in (F) and $t=-k_1+k_2+c-2$ in (D), the removed factors become
\begin{align}
q_F&=1-\frac{\chi_\Pi(p)\mathfrak a p^{k_2-2}}{\beta}
=1-\frac{\gamma\mathfrak a}{p^{k_1-1}},\label{eq:local-to-global-F}\\
q_D&=1-\frac{\chi_\Pi(p)\chi_\Sigma(p)p^{k_1-1}}{\gamma\mathfrak a}
=1-\frac{\chi_\Sigma(p)\beta}{\mathfrak a p^{k_2-2}}
=1-\frac{\beta}{\mathfrak a^\vee p^{k_2-2}},\label{eq:local-to-global-D}
\end{align}
where $\mathfrak a^\vee=\chi_\Sigma(p)^{-1}\mathfrak a$ is the ordinary arithmetic root of the dual form. This explains the character factor in (D) when the formula is written using the root $\mathfrak a$ of the original form.

\begin{prop}\label{prop:boundary-euler}
For a globally generic, non-CAP, unramified Borel-ordinary Siegel specialization of weight $(k_1,k_2)$, $k_1\geq k_2\geq3$, and an unramified ordinary $\GL_2$ specialization of weight $c=k_1-k_2+2$, the two polynomials $\mathcal B_F^{\imp}$ and $\mathcal B_D^{\imp}$ are nonzero at $t=0$, $c=k_1-k_2+2$.
\end{prop}
\begin{proof}
Use the local normalization of this subsection. Both products contain the three factors
\[
T=\left(1-\frac{p^2}{\alpha\mathfrak b_1\mathfrak a_2}\right)
  \left(1-\frac{p^2}{\alpha\mathfrak b_1\mathfrak b_2}\right)
  \left(1-\frac{p^2}{\beta\mathfrak b_1\mathfrak a_2}\right).
\]
The valuations of the quantities subtracted from $1$ are, respectively,
$k_1-1$, $k_2-2$ and $k_1-k_2+1$. They are positive, so $T\ne0$.
The four remaining factors of $\mathcal B_F^{\imp}$ are indexed by
$x\in\{\alpha,\beta\}$, $v\in\{\mathfrak a_2,\mathfrak b_2\}$ and have the form
$1-p^2/(x\mathfrak a_1v)$. For $\mathcal B_D^{\imp}$ the remaining factors are
\[
\prod_{v\in\{\mathfrak a_2,\mathfrak b_2\}}
\left(1-\frac{p^2}{\alpha\mathfrak a_1v}\right)
\prod_{x\in\{\beta,\gamma\}}
\left(1-\frac{p^2}{x\mathfrak a_1\mathfrak a_2}\right).
\]
At $t=0$, $|\mathfrak a_1|_\infty=p$. By \cite{Weissauer:Endoscopy} in general type and \cite{Del74} for the Yoshida constituents and the elliptic form, the complex absolute values in our local normalization are
$|x|_\infty=p^{3/2}$ and $|v|_\infty=p^{1/2}$. Each quantity subtracted from $1$ in these four factors has complex absolute value $p^{-1}$, and hence cannot be $1$.
\end{proof}

At weight zero, the ordinary Eisenstein vector maps to the trivial quotient of its principal series. We compute its image and use it to pass from the trilinear functional to the distinguished period.

\begin{lemma}\label{lem:quotient}
Let $\pi$ and $\sigma$ be unramified generic representations of $\GSp_4(\Q_p)$ and $\GL_2(\Q_p)$. Assume $\omega_\pi\omega_\sigma=1$ and
\[
 \operatorname{Hom}_{H(\Q_p)}(\pi\otimes\mathbf1\otimes\sigma,\C)\ne0.
\]
Consider the normalized induction
\[
 0\longrightarrow\mathrm{St}\longrightarrow
 I^+=\operatorname{Ind}_{B}^{\GL_2}(|\cdot|^{1/2}\boxtimes|\cdot|^{-1/2})
 \xrightarrow{q}\mathbf1\longrightarrow0.
\]
Normalize $q(v^{\mathrm{sph}})=1$. In the local normalization fixed above, its roots are $\{p,1\}$. The root-$p$ ordinary vector satisfies
\begin{equation}\label{eq:quotient-value}
 v_p=(1-p^{-1}s_1)v^{\mathrm{sph}},\quad s_1=\operatorname{diag}(1,p),
 \qquad q(v_p)=1-p^{-1}.
\end{equation}
The spherical-normalized trilinear functional is
\[
 \mathfrak z(v,w,u)=\lambda(v,u)q(w),
\]
where $\lambda$ is the unique local period with value $1$ on spherical vectors.
\end{lemma}
\begin{proof}
The spherical vector generates $I^+$, so its image under $q$ is nonzero. The ordinary-vector formula is \cite[\S3.3]{LZ-zeta}; its image under the trivial quotient is \eqref{eq:quotient-value}. The roots there are those after twisting by $|\cdot|^{-1/2}$, explaining the values $p,1$.

A nonzero local period $\lambda_0$ gives a nonzero trilinear form $\lambda_0\otimes q$. Here $\pi$ and $\sigma$ are irreducible generic and unramified, $I^+$ is an unramified reducible principal series with one-dimensional quotient, and the product of the three central characters is trivial. These are among the hypotheses covered by \cite[Theorem 2.3]{LZ-zeta}, which asserts that the trilinear space is one-dimensional and every nonzero form is nonzero on spherical vectors. Normalizing $\lambda_0$ on the spherical vectors gives $\lambda$, and uniqueness identifies $\lambda\otimes q$ with $\mathfrak z$.
\end{proof}

To apply Propositions~\ref{prop:local-F} and~\ref{prop:local-D} at weight zero, we must also check that their normalized functional and vectors specialize regularly at the reducible principal series.

\begin{lemma}\label{lem:regular-local-functional}
Keep $\pi,\sigma$ as in Lemma~\ref{lem:quotient}, with unitary central characters, and vary
\[
I(u)=\Ind_B^{\GL_2}(|\cdot|^{u-1/2}\boxtimes|\cdot|^{1/2-u}).
\]
In its fixed compact model, the spherical-normalized trilinear functional is regular at $u=1$. Its specialization is the functional on $\pi\otimes I^+\otimes\sigma$ of Lemma~\ref{lem:quotient} whenever the local period exists. The ordinary vectors and their fixed translates specialize regularly as well.
\end{lemma}
\begin{proof}
The compact model is $C^\infty((B\cap K)\backslash K)$, where $K=\GL_2(\mathbb Z_p)$. For a vector $f$ in this model, define a Schwartz function $\Phi_f$ on primitive rows by
\[
\Phi_f((0,1)k)=f(k),\qquad k\in K,
\]
and put it equal to zero on nonprimitive rows. Left $B\cap K$-invariance makes this well defined. Only unit values of $t$ contribute to its section, giving $f_{\Phi_f}(k,u)=f(k)$ for every $u$.

In applying the local gcd theorem, the representations $\pi$ and $\sigma$ remain fixed and unramified; only the section parameter varies. Thus the unramified case of \cite[Theorem 8.9(i)]{LPSZ1} applies to every fixed Schwartz function $\Phi_f$, including the primitive-row functions above, and shows that
\[
\mathfrak z_u(W,f,U)=\frac{1}{1-p^{-2u}}\,
\frac{Z_p(u,W,U,\Phi_f)}{L(\pi\times\sigma,u)}
\]
is regular at $u=1$: the second quotient is a Laurent polynomial in $p^{-u}$. It depends only on the section, and its generic $H$-invariance continues to $u=1$. For the spherical vector, $\Phi_f=1_{\mathbb Z_p^2}-1_{p\mathbb Z_p^2}$. The section attached to $1_{\mathbb Z_p^2}$ has value $(1-p^{-2u})^{-1}$ on the maximal compact subgroup, whereas the primitive-row section has value $1$. Thus the unramified zeta integral for $\Phi_f$ is $(1-p^{-2u})L(\pi\times\sigma,u)$, and the displayed normalization gives $\mathfrak z_u=1$ on spherical data. Uniqueness in \cite[Theorem 2.3]{LZ-zeta}, which includes $I^+$, identifies the specialization with Lemma~\ref{lem:quotient}.

In this model the ordinary vector is $(1-p^{-u}s_1)v^{\mathrm{sph}}$, with roots $p^u,p^{1-u}$. It and every fixed translate depend regularly on $u$. Hence the identities of Propositions~\ref{prop:local-F} and \ref{prop:local-D} specialize at $u=1$ after the indicated Euler factor has been deleted from the polynomial.
\end{proof}

Dividing the specialized local formulas by the quotient value $1-p^{-1}$ gives the following period values.

\begin{prop}\label{prop:ptest}
At a specialization as in Proposition~\ref{prop:boundary-euler}, suppose that the local period is nonzero, and normalize $\lambda$ as in Lemma~\ref{lem:quotient}. Its values on the ordinary vectors used in the two improved constructions are, for $n\geq2$,
\begin{align}
 \lambda(u_{\mathrm{Iw}}W'_n,W_{\mathfrak a_2})
 &=\frac{C_p}{1-p^{-1}}\mathcal B_F^{\imp},\label{eq:localF}\\
 \lambda(\tau^\sharp W'_n,W_{\mathfrak a_2})
 &=\frac{A_nC_p}{1-p^{-1}}\mathcal B_D^{\imp},\label{eq:localD}
\end{align}
at $\mathfrak a_1=p$, $\mathfrak b_1=1$. Both values are nonzero.
\end{prop}
\begin{proof}
By Lemma~\ref{lem:regular-local-functional}, the functional and vectors in Propositions~\ref{prop:local-F} and \ref{prop:local-D} specialize at $\mathfrak a_1=p$, $\mathfrak b_1=1$. Lemma~\ref{lem:quotient} then gives \eqref{eq:localF}--\eqref{eq:localD}; in region (D), $d_n=\operatorname{diag}(p^n,1)$ acts trivially on $\mathbf1$. The Euler products are nonzero by Proposition~\ref{prop:boundary-euler}, and so are $A_n$, $C_p$ and $1-p^{-1}$.
\end{proof}

Corollary~\ref{cor:coherent-nonvanishing} supplies the archimedean nonvanishing for the compact types selected by the coherent pairings. Together with Proposition~\ref{prop:ptest}, this gives the prescribed local vectors needed for the pole criterion.

\subsection{Level normalization and comparison with Loeffler--Rivero}\label{subsec:level-comparison}

\subsubsection{The modular degeneracy and its adjoint}\label{subsubsec:modular-adjoint}
The pairing in the interpolation formula of \cite[\S5.3]{LR3} is defined on a product of modular curves at level $p^2$. To compare it with ours, we use the modular degeneracy of \cite[\S4.8]{LR3} from the common-level $H$-surface to this product:
\[
d_2=\operatorname{diag}(p^2,1),\qquad j_2=\iota_p:\mathscr Z_2\longrightarrow \mathscr Y_2,\qquad (h_1,h_2)\longmapsto(h_1d_2,h_2).
\]
Here $\mathscr Z_2$ is the common-level $H$-surface and $\mathscr Y_2$ the union of product components at the prescribed modular levels: $(\delta,\delta)$ maps to $(p^2\delta,\delta)$. At classical weights $(r_1,r_2;a,b)$, let $\mathcal V=\omega_1^a\boxtimes\bigl(\omega_2^{2-b}\otimes(\mathrm{det}_2)^{-r_1}\bigr)$ on $\mathscr Y_2$, and let $\mathcal B_2=H_D^{\mathrm{sc}}$ be its transported coefficient on $\mathscr Z_2$. The coefficient map $\theta:j_2^*\mathcal V\to\mathcal B_2$ is the specialization of the map in Section~\ref{subsubsec:family-support}: the normalized isogeny map on the first coefficient tensored with the identity on the second. Its dual, followed by derived pushforward and Grothendieck--Serre trace, defines $(j_2)_!$ on the Siegel quotient $\mathcal B_2^\vee\otimes K_{\mathscr Z_2}$. The projection formula is
\begin{equation}\label{eq:modular-projection}
\Tr_{\mathscr Y_2}((j_2)_!\alpha\cup v)
=\Tr_{\mathscr Z_2}(\alpha\cup\theta(j_2^*v)),
\end{equation}
with the Siegel factor first. Write $\iota_2^\sharp$ for the translated embedding at this common level. Then $\mathscr R_D=(j_2)_!q_D^{\mathrm{sc}}(\iota_2^\sharp)^*$ gives the class on the product of curves used in \cite[\S5.3]{LR3}. We denote the resulting classical scalar pairing by
\begin{equation}\label{eq:classical-B2}
B_2(x;E,\upsilon)=\nu_K\Tr_{\mathscr Y_2}\bigl(\mathscr R_Dx\cup(E\boxtimes\upsilon)\bigr).
\end{equation}
Here $E$, $\upsilon$ and $x$ are the Eisenstein, degree-one modular and Siegel inputs, including the fixed tame translations and Gauss sums. Both traces in \eqref{eq:modular-projection} are geometric.

\subsubsection{Variation of the auxiliary level}\label{subsubsec:level-transition}

Theorem~\ref{thm:analytic-algebraic} compares pairings at the same auxiliary level. As in \cite{LR3}, the analytic construction requires a level deep enough for its supports, while $B_2$ fixes the classical normalization. We now compute how the scalar changes with the auxiliary level. This requires following the inverse Siegel trace, the trace on $H$, the modular translation $d_T=\operatorname{diag}(p^T,1)$ and its coefficient map. The same calculation will apply to both Eisenstein inputs.

Put $K_T=\operatorname{Iw}_G(p^T)$, let $\mathcal M_T\subset H(\mathbb Z_p)$ be the inverse image modulo $p^T$ of
\[
S_0=\{(\operatorname{diag}(c,d),\operatorname{diag}(d,c)):c,d\in\mathbb Z_p^\times\},
\]
and set $\mathcal L_2=\mathcal M_2$ and $\mathcal L_T=\mathcal L_2\cap\tau^\sharp K_T(\tau^\sharp)^{-1}$. The intermediate level $\mathcal L_T$ allows us to compare the Siegel traces with $d_2$ fixed; the analytic pairing uses $(\mathcal M_T,d_T)$. On the selected classical Siegel line, write $x_T=\Tr_G^{-1}x_2$, using the unnormalized geometric trace.

For $*\in\{\mathrm{dep},\mathrm{ord}\}$, let $B_T^*$ be the scalar pairing at $(\mathcal M_T,d_T)$ with class $x_T$, the fixed elliptic old vector and the indicated Eisenstein input. We use the coefficients $[L_1],H_D^{\mathrm{sc}}$ and the tame multiplier $\nu_K$ of $B_2$ throughout. Algebraically these pairings are defined for every $T\geq2$; the analytic construction requires sufficiently large $T$. At a pole of the ordinary input, specialization means specialization after regularization as in Lemma~\ref{lem:regularized-eisenstein}.

\begin{prop}\label{prop:running-level}
Put $\mathfrak r=p\mathfrak a^\vee/\beta^\circ$.
\begin{enumerate}[label=\textup{(\roman*)}]
\item At every good classical point and for $T\geq2$, both Eisenstein inputs satisfy
\begin{equation}\label{eq:running-level}
B_{T+1}^*=\mathfrak r B_T^*,\qquad *\in\{\mathrm{dep},\mathrm{ord}\}.
\end{equation}
\item The operator
\[
\mathscr V=p^{-1}\mathcal U'_{\Kl}(\mathcal U'_{\Si})^{-1}(\langle p\rangle^{-1}U_p)^{-1}
\]
acts by $\mathfrak r^{-1}$ on the tensor product of the refined cuspidal lines.
\item The normalized pairings $\widehat B_T^*=\mathfrak r^{2-T}B_T^*$ are independent of $T$. For sufficiently large $T$, the functions $\widehat B_T^{\mathrm{dep}}$ and $\ell(\mathbf t_D)\widehat B_T^{\mathrm{ord}}$ are analytic on $X$, and their good classical specializations are $B_2$ with the corresponding depleted or regularized ordinary input.
\end{enumerate}
\end{prop}
\begin{proof}
Fix a good classical point and put $a=t_D$. Only the data at $p$ change with $T$; their four contributions are as follows.

\smallskip\noindent\textit{(a) The level on $H$.} The level indices are
\[
[K_2:K_T]=[\mathcal L_2:\mathcal L_T]=p^{4(T-2)},\qquad
[\mathcal L_2:\mathcal M_T]=p^{5(T-2)}.
\]
The first quotient has the four opposite-root directions. In $\mathcal M_T/\mathcal M_{T+1}$ there are four off-diagonal directions and one diagonal direction transverse to $S_0$. In addition, representatives for $\mathcal L_T/\mathcal M_T$ are
\[
\left(\begin{pmatrix}1&0\\c&1\end{pmatrix},\begin{pmatrix}1&c\\0&1\end{pmatrix}\right),\qquad c\in p^2\mathbb Z_p/p^T\mathbb Z_p.
\]
The modular adjoint of \eqref{eq:modular-projection} uses the geometric trace, which is not divided by the degree. With the fixed Haar measure, the scalar pairing at level $\mathcal M_T$ is therefore expressed by $\operatorname{vol}(\mathcal M_T)^{-1}$ times the corresponding adelic integral, up to factors independent of $T$. Since $[\mathcal M_T:\mathcal M_{T+1}]=p^5$, raising the level contributes $p^5$.

\smallskip\noindent\textit{(b) The Siegel class.} The inverse trace is taken on the minus finite-slope cohomology. For levels greater than the analyticity radius, it exists and is Hecke-compatible by \cite[Theorem 6.4.10]{boxerpilloni20}; at classical weights we use its restriction to the refined line. The vectors $W'_T$ are the dual \emph{Iwahori} vectors of \cite[\S3.6]{LZ-zeta}. Their normalization uses averaged traces, as in \S3.4 of that reference:
\[
\frac{1}{[K_T:K_{T+1}]}\sum_{g\in K_T/K_{T+1}}gW'_{T+1}=W'_T.
\]
Consequently, if $x_2\leftrightarrow c_yW'_2$, the geometric trace condition $\Tr_G(x_T)=x_2$ gives
\[
x_T\longleftrightarrow c_y p^{-4(T-2)}W'_T.
\]
Thus the inverse geometric trace contributes $p^{-4}$ at each step.

\smallskip\noindent\textit{(c) The first modular input.} Let $\mathcal N_{d_T}$ be the normalized coefficient map for the first modular isogeny, and let $R_{d_T}^{\mathrm{un}}$ denote right translation by $d_T=\operatorname{diag}(p^T,1)$ in the unitary automorphic realization. On Tate curves the former preserves $dz/z$ and sends $f(q)$ to $f(q^{1/p^T})$. In weight $a$, the latter sends $f(\tau)$ to $p^{-Ta/2}f(\tau/p^T)$. Thus
\begin{equation}\label{eq:translation-models}
\mathcal N_{d_T}=p^{Ta/2}R_{d_T}^{\mathrm{un}}.
\end{equation}
The second modular coefficient contributes the same factor at levels $T$ and $T+1$, since $d_T$ acts on the first curve. By \eqref{eq:translation-models}, increasing $T$ by one therefore contributes $p^{a/2}$ relative to unitary translation.

\smallskip\noindent\textit{(d) The local trilinear form.} Proposition~\ref{prop:local-D} gives the same level factor
\[
A_T=\left(\frac{p^2}{\beta_{\mathrm{loc}}\mathfrak b_{2,\mathrm{loc}}}\right)^T
\]
for both Eisenstein inputs. The two Jacquet factors are $(p^2/\beta_{\mathrm{loc}})^T$ and $\mathfrak b_{2,\mathrm{loc}}^{-T}$; after the simultaneous translations are cancelled, the remaining torus integral is independent of $T$. The substitutions at the end of Section~\ref{subsec:local-improvement} give
\[
\beta_{\mathrm{loc}}\mathfrak b_{2,\mathrm{loc}}=p^{2+a/2}\frac{\beta^\circ}{\mathfrak a^\vee},\qquad
\frac{A_{T+1}}{A_T}=p^{-a/2}\frac{\mathfrak a^\vee}{\beta^\circ}.
\]
Combining (a)--(d), we obtain
\[
\begin{array}{c|l|c}
&\text{change from }T\text{ to }T+1&\text{factor}\\\hline
\text{(a)}&\text{modular geometric trace}&p^5\\
\text{(b)}&\text{inverse Siegel geometric trace}&p^{-4}\\
\text{(c)}&\text{coefficient map relative to unitary translation}&p^{a/2}\\
\text{(d)}&\text{local trilinear evaluation}&p^{-a/2}\mathfrak a^\vee/\beta^\circ
\end{array}
\]
The powers $p^{\pm a/2}$ cancel, so the product is $\mathfrak r$, which is $p$ times an analytic unit, although (c) and (d) are not separately analytic. The other local functionals and the coefficient quotient do not depend on $T$, which proves (i); the argument does not require the global period to be nonzero.

On the scalar Siegel coefficient $[L_1]$, the normalized operators $\mathcal U'_{\Kl}$ and $\mathcal U'_{\Si}$ have eigenvalues $\alpha\beta^\circ$ and $\alpha$, respectively. Equation~\eqref{eq:twisted-U} gives the same eigenvalues after transport to $E_D$. On the elliptic line, $\langle p\rangle^{-1}U_p$ acts by $\mathfrak a^\vee$ by \eqref{eq:elliptic-adjoint}. Hence $\mathscr V$ acts by $p^{-1}\beta^\circ/\mathfrak a^\vee=\mathfrak r^{-1}$, proving (ii). In particular, the normalization only involves operators which are invertible on the cuspidal lines.

For sufficiently large $T$, the pairings $B_T^{\mathrm{dep}}$ and $\ell(\mathbf t_D)B_T^{\mathrm{ord}}$ are analytic on a common affinoid by Proposition~\ref{prop:supported-pairing} and Lemma~\ref{lem:regularized-eisenstein}. Part (i) gives the transition identity at good classical points. Applying Lemma~\ref{lem:three-variable-uniqueness} below to the difference, after regularization in the ordinary case, gives the same identity on $X$. For the ordinary input we may now cancel $\ell(\mathbf t_D)$ in the ring of meromorphic functions. Thus both normalized pairings are independent of $T$, as a meromorphic identity in the ordinary case. At classical points, iterating (i) down to $T=2$ and using \eqref{eq:modular-projection} identifies its value with $B_2$. This proves (iii).
\end{proof}

At classical points, the calculation gives $A_T=p^{-Tt_D/2}(\mathfrak a^\vee/\beta^\circ)^T$. In particular, the combined local and coefficient factor at level $p^2$ is
\begin{equation}\label{eq:absolute-level-factor}
p^{t_D}A_2=\left(\frac{\mathfrak a^\vee}{\beta^\circ}\right)^2.
\end{equation}
The right-hand side is an analytic unit. At a matching ordinary Yoshida point, $t_D=0$ and $\mathfrak a^\vee=\beta^\circ=\alpha_g$, so $A_T=1$ and the local boundary value of Proposition~\ref{prop:ptest} is unchanged.

To obtain the level identity in families, we use uniqueness on the three-dimensional space $X$. We vary all three weights, since fixing the Siegel weights leaves only finitely many critical Eisenstein weights.

\begin{lemma}\label{lem:three-variable-uniqueness}
An analytic function on $X$ is determined by its values at the good crystalline points whose weights satisfy
\[
R_1-R_2\geq3,\qquad R_2\geq3,\qquad 4\leq a\leq R_2+1.
\]
The same assertion holds for meromorphic functions whose product with $\ell(\mathbf t_D)$ is analytic.
\end{lemma}
\begin{proof}
Ordinary control makes every geometric point above these weights good. First work on the weight polydisc and let a function vanish at the indicated weights. Fix integers $a\geq4$ and $R_2\geq\max(3,a-1)$ in the prescribed congruence classes. The allowed values of $R_1$ contain an integral tail with an accumulation point in its weight disc, so one-variable uniqueness gives vanishing for every $R_1$. With $a$ still fixed, apply the same argument to the allowed tail of $R_2$, and then to the allowed values of $a$. The function is therefore zero on the polydisc.

Now let a function on $X$ vanish at every geometric point above these weights. Since $X$ is finite flat over the polydisc, multiplication by this function has a characteristic polynomial with analytic coefficients on the base. On each of the indicated fibres all its eigenvalues are zero, so every nonleading coefficient vanishes there. The preceding argument makes these coefficients identically zero. Cayley--Hamilton then shows that the function is nilpotent, and reducedness of $X$ makes it zero. For the meromorphic assertion, apply this to the product with $\ell(\mathbf t_D)$.
\end{proof}

\begin{rk}\label{prop:LR-comparison}
With the same generators, Gauss sums, Fourier-related tame data and tame trace, Proposition~\ref{prop:running-level} identifies the depleted classical evaluations with those of \cite[Definition 5.5]{LR3}. Lemma~\ref{lem:three-variable-uniqueness} then identifies the three-variable functions. The same normalized pairing defines the ordinary extension. Changing the generators or tame trace multiplies both inputs by
\begin{equation}\label{eq:LR-bases}
c_{\mathrm{LR}}=c_Ku_\Pi\,\pi_\Sigma^*u_\Sigma,\qquad
\xi_{\mathrm{LR}}=u_\Pi\underline\xi,\qquad
u_\Sigma=P(g^\vee,\eta_{\mathrm{LR}}).
\end{equation}
Here $\pi_\Sigma$ is the map to the $\GL_2$ parameter space. By Lemma~\ref{lem:elliptic-supported-control}, the class $\eta_{\mathrm{LR}}$ in the positive long-Weyl realization is $u_\Sigma\underline\upsilon$; its Atkin--Lehner transport is already included. Lemma~\ref{lem:boundary-branching} and \eqref{eq:modular-projection} identify the coefficient quotient and modular adjoint, while \eqref{eq:katz-transport} identifies the Eisenstein data. Bilinearity gives \eqref{eq:LR-bases} at good classical points for either input. Applying Lemma~\ref{lem:three-variable-uniqueness}, after regularization for the ordinary input, gives the same comparison on $X$.

The tame trace ratio $c_K$ is $1$ when both conventions include $\nu_K$, and $\nu_K^{-1}$ when the LR pairing omits it; the volume convention is that of Section~\ref{subsec:improved-construction}, with $H(\mathbb Z_p)$ of volume one. Our common generators and trace give $c_{\mathrm{LR}}=1$.
\end{rk}

\subsection{Poles and the comparison of periods}\label{subsec:poles-comparison}

We now express the Eisenstein residue through the distinguished period. Fix a good Siegel point $P$ of weight $(k_1,k_2)$ and a crystalline point $Q_0$ of the $\GL_2$ family of weight $c_0=k_1-k_2+2$. Retain the refinements of Section~\ref{sec:review} and assume both representations are unramified at $p$. The elliptic weight map is \'etale at $Q_0$, so we may work on a small disc where the generators and specialization maps are regular. Choose the character component on which both Eisenstein characters are trivial at $Q_0$.

The parameter $z$ measures the change in elliptic weight from $Q_0$: at a classical point of weight $c$, it is simply $c-c_0$. In terms of weight characters this is
\begin{equation}\label{eq:coordinates}
 z(Q)=\ell(\mathbf c_Q\mathbf c_{Q_0}^{-1}),\qquad
 \ell(\mathbf t_F)=-z,\qquad \ell(\mathbf t_D)=z.
\end{equation}
Thus the two Eisenstein weights move in opposite directions as the elliptic weight varies, and the \'etale condition makes $z$ a local parameter on the family. Since $\ell(\mathbf t_\diamond)$ is not identically zero on the fibre over $P$, the meromorphic functions of Definition~\ref{def:joint-functions} restrict to this fibre. We denote them by $L_p^{(\diamond),\imp}(\Pi\times\underline\Sigma)$, keeping the specialized Siegel generator in their normalization.

We shrink the disc so that $Q_0$ is its only point where the Eisenstein characters, including their finite-order parts, are trivial. Residues are taken in $z$: for a function with at most a simple pole, $\operatorname{Res}_{z=0}L=(zL)(Q_0)$.

Let $\mathscr B_{F,Q}$ and $\mathscr B_{D,Q}$ denote the specialized pairings of Definition~\ref{def:joint-functions}, with the prescribed pullbacks at $p$ and Gauss sums. Write $\gamma_F,\gamma_D$ for the tame translations in the two regions, fixed on the disc as in Section~\ref{subsec:improved-construction}, and $g_Q,\upsilon_Q$ for the normalized elliptic inputs in degrees zero and one.

\begin{prop}\label{prop:residue}
With $\operatorname{Res}_{z=0}L=(zL)(Q_0)$, one has
\begin{align}
 \operatorname{Res}_{z=0}L_p^{(F),\imp}
 & =(1-p^{-1})\mathscr B_{F,Q_0}(\gamma_F\eta_P,\mathbf1\boxtimes g_{Q_0}),\label{eq:resF}\\
 \operatorname{Res}_{z=0}L_p^{(D),\imp}
 & =-(1-p^{-1})\mathscr B_{D,Q_0}(\gamma_D\xi_P,\mathbf1\boxtimes\upsilon_{Q_0}).\label{eq:resD}
\end{align}
\end{prop}
\begin{proof}
By multilinearity, multiplication by $z$ regularizes the Eisenstein input. Specialize at $Q_0$ using Lemma~\ref{lem:boundary-branching} and Theorem~\ref{thm:analytic-algebraic}; in region (D), Proposition~\ref{prop:running-level} identifies the normalized scalar with \eqref{eq:classical-B2}. The Eisenstein residue \eqref{eq:eisenstein-residue} and the two signs in \eqref{eq:coordinates} give the stated coefficients.
\end{proof}

For later use, abbreviate the two constant-input pairings as
\begin{equation}\label{eq:constant-pairing}
\begin{aligned}
B_F(\mathbf1)&:=\mathscr B_{F,Q_0}(\gamma_F\eta_P,\mathbf1\boxtimes g_{Q_0}),\\
B_D(\mathbf1)&:=\mathscr B_{D,Q_0}(\gamma_D\xi_P,\mathbf1\boxtimes\upsilon_{Q_0}).
\end{aligned}
\end{equation}

To compare the residue with a complex period, we evaluate a rational functional on two realizations of the same tensor product of cuspidal eigenlines, following \cite[Theorem 3.10]{LZ25-poles}. Write
\[
 \nu_{F,p}=\eta_P\otimes g_{Q_0},\qquad
 \nu_{D,p}=\xi_P\otimes\upsilon_{Q_0}.
\]
Let $E\subset\mathbb C$ contain the eigenvalues, refinement parameters and coefficients of the finite-level maps, and fix an embedding $\iota_p:E\hookrightarrow\overline{\mathbb Q}_p$. For each $\diamond=F,D$, use the scalar coefficients $[L_{i_{\diamond}}],H_\diamond^{\mathrm{sc}}$ of Lemma~\ref{lem:boundary-branching}. The tensor product of their two refined classical eigenlines is a one-dimensional $E$-space, denoted $\mathscr P_\diamond$. Its $p$-adic realization is generated by $\nu_{\diamond,p}$, while its complex realization is generated by the tensor $\nu_{\diamond,\infty}$ of the coherent classes of the corresponding automorphic forms.

Choose an $E$-basis $\nu_{\diamond,E}$ of this line. There are nonzero periods $\Omega_{\diamond,p}$ and $\Omega_{\diamond,\infty}$ such that
\begin{equation}\label{eq:algebraic-period}
 \frac{\nu_{\diamond,p}}{\Omega_{\diamond,p}}
 =\nu_{\diamond,E}\otimes_{E,\iota_p}1,
 \qquad
 \frac{\nu_{\diamond,\infty}}{\Omega_{\diamond,\infty}}
 =\nu_{\diamond,E}\otimes_E1.
\end{equation}
These periods are products of the Siegel and elliptic comparison periods. In degree one the elliptic line is Serre-dual to the holomorphic line for $\Sigma^\vee$, with its adjoint Hecke action.

By \eqref{eq:coefficient-conventions}, tensor transport replaces $\mathscr P_\diamond$ by $\mathscr P_\diamond\otimes E\{m\}$, with basis $\nu_{\diamond,E}\otimes e_m$, without changing these periods.

The rational functional we compare is obtained by inserting the constant form into $\nu_K B_\diamond^{\mathrm{sc}}$, with $\nu_K$ as in Remark~\ref{rk:tame-trace-normalization}. By \eqref{eq:boundary-branching}, it can equally be written in the character-valued realization as
\[
\begin{aligned}
x\otimes y\longmapsto\nu_K B_\diamond^{\mathrm{sc}}(\gamma_\diamond x,\mathbf1\boxtimes y)=\nu_K e_m^\vee P_\diamond\bigl(T_{d_\diamond}(\gamma_\diamond x),T_{h_\diamond}(\mathbf1\boxtimes y)\bigr).
\end{aligned}
\]
Here $e_m^\vee$ extracts the coefficient of $e_m$. The pullbacks use the finite-level conventions of Sections~\ref{subsec:classical-comparisons} and~\ref{subsec:level-comparison}.

The tame translate is applied before tensor transport. By Theorems~\ref{thm:analytic-algebraic} and~\ref{thm:scalar-period-comparison}, and Proposition~\ref{prop:running-level} in $(D)$, the $p$-adic value of this functional is $B_\diamond(\mathbf1)$ of \eqref{eq:constant-pairing}.

Over $\C$, Theorem~\ref{Integral:New:Pairing} identifies the same functional with
\begin{equation}\label{eq:analytic-period}
\nu_K A_\diamond\mathcal I_\diamond,\qquad
\mathcal I_\diamond=\int_{Z_G(\mathbb A)H(\mathbb Q)\backslash H(\mathbb A)}
\Psi_\diamond(\iota(h))f_\diamond(h_2)\,dh.
\end{equation}
Here $\Psi_\diamond,f_\diamond$ are the unitary representatives of Theorem~\ref{Integral:New:Pairing}, with the prescribed finite translates. Its proportionality constant $A_\diamond$, expressed in the $e_m$-coordinate, includes the cocycle evaluation of Lemma~\ref{C:Evaluated} and the central volume factor. The separate multiplier $\nu_K$ comes from Definition~\ref{def:joint-functions}.

\begin{prop}\label{lem:boundary-comparison}\label{cor:periodcomparison}
For $\diamond=F,D$, the normalized complex period satisfies
\begin{equation}\label{eq:comparisonc}
 \frac{\nu_K A_\diamond\mathcal I_\diamond}{\Omega_{\diamond,\infty}}\in E,
\end{equation}
and the residue is given by
\begin{equation}\label{eq:comparisonp}
 \frac{\operatorname{Res}_{z=0}L_p^{(\diamond),\imp}}{\Omega_{\diamond,p}}
 =\varepsilon_\diamond(1-p^{-1})\,
 \iota_p\!\left(
 \frac{\nu_K A_\diamond\mathcal I_\diamond}{\Omega_{\diamond,\infty}}
 \right),
 \qquad \varepsilon_F=1,\quad\varepsilon_D=-1.
\end{equation}
In particular, $L_p^{(\diamond),\imp}$ has a pole at $Q_0$ if and only if $\mathcal I_\diamond\ne0$.
\end{prop}
\begin{proof}
Let $b_\diamond\in E$ be the value of this functional on $\nu_{\diamond,E}$. The functional is defined over $E$, since this holds for the constant section, tame correspondence, coefficient contraction and trace. Evaluating in the two realizations and using \eqref{eq:algebraic-period} gives
\[
\nu_K A_\diamond\mathcal I_\diamond=\Omega_{\diamond,\infty}b_\diamond,
\qquad
B_\diamond(\mathbf1)
=\Omega_{\diamond,p}\,\iota_p(b_\diamond),
\]
The first equality proves \eqref{eq:comparisonc}, and the second, together with Proposition~\ref{prop:residue}, proves \eqref{eq:comparisonp}. The comparison periods are nonzero and $\iota_p$ is injective, so the residue vanishes exactly when the complex period does. Since the possible pole is simple, this proves the final assertion.
\end{proof}

\begin{rk}\label{rk:factor-cancellation}
The factor $1-p^{-1}$ in \eqref{eq:comparisonp} is the global counterpart of the image $q(v_p)=1-p^{-1}$ of the ordinary Eisenstein vector in the trivial quotient (Lemma~\ref{lem:quotient}). The local periods of Proposition~\ref{prop:ptest} contain the inverse factor, so the two cancel when $\mathcal I_\diamond$ is written in terms of them: at a matching Yoshida point, where $A_2=1$, the residue is governed by the improved local factor $C_p\mathcal B_\diamond^{\imp}$, as the improved interpolation formula suggests.
\end{rk}

For factorizable vectors, we can express the comparison directly in terms of the complex residue. Choose a Schwartz function $\varphi$ with $\widehat\varphi(0)\ne0$, spherical outside $S$. With the measures of Proposition~\ref{Period:Equals:Residue}, that proposition and Corollary~\ref{Global:Zeta:Integral:Result} give
\[
\frac{\operatorname{Res}_{z=0}L_p^{(\diamond),\imp}}{\Omega_{\diamond,p}}
=\varepsilon_\diamond(1-p^{-1})\,
\iota_p\!\left(
\frac{\nu_K A_\diamond}{\Omega_{\diamond,\infty}c_E\widehat\varphi(0)}
\operatorname{Res}_{s=1}\Lambda(s,\Pi\times\Sigma)\,
Q_\infty(1)\prod_{v\in S,\ v<\infty}Q_v(1)
\right).
\]
Equation~\eqref{eq:comparisonc} gives the algebraicity of the whole expression inside $\iota_p$. Here the $Q_v$ are the normalized local zeta integrals of Corollary~\ref{Global:Zeta:Integral:Result}, computed with the same automorphic and Schwartz data as the global integral. The formula extends by linearity to finite sums of factorizable vectors.

Proposition~\ref{cor:periodcomparison} detects a pole through the period on the chosen vectors. For a distinguished pair, the test-vector results at $p$ and infinity allow us to choose tame translations $\gamma_F,\gamma_D$ making both periods nonzero. We fix them at $(P,Q_0)$ and apply them to the family classes throughout the neighbourhood, using the same translation for the original and improved functions in each region.

Write $\Sigma=\underline\Sigma(Q_0)$ and retain the compatibility between integral and Langlands factors at ramified places in Theorem~\ref{intro:complex}. We obtain the following $p$-adic counterpart of the complex criterion.

\begin{theorem}\label{thm:criterion}
There are tame test data, fixed on a neighbourhood of $Q_0$, for which the following are equivalent:
\begin{enumerate}
\item $L_p^{(F),\imp}$ has a pole at $Q_0$;
\item $L_p^{(D),\imp}$ has a pole at $Q_0$;
\item $\Pi=\mathcal Y(\Sigma^\vee\boxtimes\Sigma')$, where $\Sigma^\vee$ is the constituent of weight $c_0$ and $\omega_{\Sigma'}=\omega_\Sigma^{-1}$;
\item the global $H$-period on $\Pi\otimes\Sigma$ is nonzero.
\end{enumerate}
When these conditions hold, both poles are simple. For every choice of tame data, a pole of either function implies the Yoshida condition in (3).
\end{theorem}
\begin{proof}
For any fixed tame data, Proposition~\ref{cor:periodcomparison} identifies a nonzero residue with a nonzero value of the global period. The complex pole criterion then gives $\Pi=\mathcal Y(\Sigma^\vee\boxtimes\Sigma')$, with $\Sigma^\vee$ the smaller constituent because its weight is $c_0$.

Conversely, suppose this Yoshida condition holds. Corollary~\ref{lem:infinity} gives a nonzero period with the archimedean representatives of each of our two coherent pairings. Fix the vectors away from $p$ in one such nonzero global period. This gives a nonzero $H(\Q_p)$-invariant functional on $\Pi_p\otimes\Sigma_p$. By Lemma~\ref{lem:quotient} it is a nonzero multiple of the local period $\lambda$. Proposition~\ref{prop:ptest} shows that $\lambda$ is nonzero on our prescribed ordinary vectors. Replacing the vectors at $p$ by these therefore preserves nonvanishing. Applying this argument separately in the two degrees gives $\mathcal I_F\ne0$ and $\mathcal I_D\ne0$.

It remains to realize these nonzero periods using tame translates of our fixed generators. At the finite places away from $p$, the translates of any nonzero vector span the irreducible representation. A nonzero period must therefore be nonzero on some pair of translates of the chosen Siegel and elliptic vectors. The translation identity in the proof of Corollary~\ref{cor:coherent-nonvanishing}, applied to $a\in\GL_2(\mathbb A_f^p)$, transfers the elliptic translate to the Siegel factor while leaving the vectors at $p$ and infinity fixed. We pull the translated classes to a common tame level, as in Lemma~\ref{lem:ordinary-control}. This gives the required $\gamma_F$ and $\gamma_D$. Applying each fixed translation to the corresponding family generator gives the functions on a neighbourhood of $Q_0$.

For these choices, \eqref{eq:comparisonp} gives a nonzero residue in both regions. The poles are simple by Proposition~\ref{prop:residue}, completing the equivalence.
\end{proof}

We now allow the Siegel point to vary near the matching Yoshida point $x_0=(P,Q_0)$, using the family classes of Lemma~\ref{lem:ordinary-control} and keeping the generators and tame translations fixed. On a small product neighbourhood $X$ of $x_0$, let $D$ be the divisor on which the Eisenstein characters $\htt_F$ and $\htt_D$ are trivial (these characters are inverse to each other). Regularizing the analytic pairing gives
\[
\mathcal R_\diamond=\ell(\mathbf t_\diamond)L_p^{(\diamond),\imp}\in\mathcal O(X),
\qquad
\varrho_\diamond=\mathcal R_\diamond|_D\in\mathcal O(D).
\]
The coefficient of the possible pole is therefore analytic along $D$. On the fibre over $P$, the coordinates \eqref{eq:coordinates} give
\[
\varrho_F(x_0)=-\operatorname{Res}_{z=0}L_p^{(F),\imp},
\qquad
\varrho_D(x_0)=\operatorname{Res}_{z=0}L_p^{(D),\imp}.
\]
Both values are nonzero, so after shrinking, both $\varrho_\diamond$ are units on $D$. At nearby good classical points satisfying the hypotheses of the complex criterion, Proposition~\ref{cor:periodcomparison} gives a nonzero distinguished period. The complex criterion therefore identifies these specializations as Yoshida lifts, with the dual of the corresponding elliptic representation as their smaller constituent.

For fixed $P$, the pole is isolated in the elliptic variable. As $P$ varies, these poles form $D$, with nonzero residue. Since the elliptic weight map is \'etale at $Q_0$, this divisor is locally the graph of the $\GL_2$ point with weight character $\mathbf k_1\mathbf k_2^{-1}x^2$, where $x$ is the identity character, over the Siegel neighbourhood.

Since $L_p^{(\diamond)}$ is analytic, Proposition~\ref{teo:interpolation} gives
\[
\ell(\mathbf t_\diamond)L_p^{(\diamond)}=e_\diamond\mathcal R_\diamond
\quad\Longrightarrow\quad
(e_\diamond|_D)\varrho_\diamond=0.
\]
As $\varrho_\diamond$ is a unit, $e_\diamond$ vanishes along $D$, cancelling the pole throughout this neighbourhood.

\subsection{Euler factors and interpolation}\label{subsec:improved-interpolation}

We now apply the local calculations to the original and improved functions, using their common level normalization. In the Fourier convention of \eqref{eq:katz-transport}, their Eisenstein inputs are
\[
(\mathcal E^\flat(\htt_\diamond))^{[p]}=\mathcal E(\htt_\diamond-1,0),\qquad
\mathcal E^\flat(\htt_\diamond),\qquad\diamond\in\{F,D\}.
\]
In (F), the change from the depleted cuspidal input to its ordinary stabilization is accounted for by Proposition~\ref{prop:local-F}.

The analytic factors $e_F,e_D$ of \eqref{eq:removed-factors} specialize to $q_F,q_D$ by \eqref{eq:local-to-global-F}--\eqref{eq:local-to-global-D}.

On the component $\chi_\Pi\chi_\Sigma=1$, the two factors, viewed on the same three weight variables, satisfy
\[
e_D=-\frac{e_F}{\gamma^\circ\mathfrak a}.
\]
Since $\gamma^\circ\mathfrak a$ is an analytic unit, the two removed Euler factors have the same zero divisor. At a matching ordinary Yoshida point, write $\alpha_g$ for the unit root of the smaller constituent. Then
\[
\beta^\circ=\alpha_g,\qquad
\gamma^\circ=\chi_g(p)/\alpha_g,\qquad
\mathfrak a=\chi_g(p)^{-1}\alpha_g,
\]
and both removed factors vanish. Their effect on the pole is described by Theorem~\ref{thm:criterion} and the identity below.

\begin{prop}\label{teo:interpolation}
On the respective three-variable domains, the functions of Definition~\ref{def:joint-functions} satisfy
\begin{equation}\label{eq:global-improvement}
\begin{split}
L_p^{(F)}=(1-\gamma^\circ\mathfrak a)L_p^{(F),\imp},\qquad
L_p^{(D)}=(1-\beta^\circ/\mathfrak a^\vee)L_p^{(D),\imp}.
\end{split}
\end{equation}
\end{prop}

\begin{proof}
At a good crystalline point in the regular interpolation range, the Fourier transport \eqref{eq:katz-transport} identifies the Eisenstein data. In region $(F)$ the classical comparison is \cite[\S\S7--8]{LPSZ1} and \cite[\S10]{LZ-BK}. In region $(D)$, Theorem~\ref{thm:analytic-algebraic}, \eqref{eq:modular-projection} and Proposition~\ref{prop:running-level} identify both scalar evaluations with the algebraic level-two pairing. Their automorphic realizations have the same coefficient quotient, periods and archimedean factors. The local Whittaker vectors are those fixed in Section~\ref{subsec:improved-construction}; in particular the elliptic ordinary vector is the one identified in Lemma~\ref{lem:elliptic-supported-control}.

Propositions~\ref{prop:local-F} and \ref{prop:local-D} compare the two local evaluations with factors $q_F$ and $q_D$, respectively. In region (D), the factor $p^{t_D}A_2$ in \eqref{eq:absolute-level-factor} occurs on both sides. Expressing the local parameters in motivic normalization gives $e_F$ and $e_D$ by \eqref{eq:local-to-global-F}--\eqref{eq:local-to-global-D}. The remaining local functionals are identical, so the equalities persist after taking their product. This proves \eqref{eq:global-improvement} at the classical points, including those where the global period vanishes.

To pass to families, multiply the difference between the two sides by $\ell(\mathbf t_\diamond)$. By Lemma~\ref{lem:regularized-eisenstein}, this gives an analytic function. In region (D), it vanishes at the classical points of Lemma~\ref{lem:three-variable-uniqueness}, and hence vanishes identically. In region (F), the same argument applies to the range $r_2\geq3$, $c\geq2$, $r_1-r_2+2-c\geq2$, with $r_1-r_2\geq3$. Indeed, for fixed $r_2,c$ there is an integral tail of possible $r_1$; applying one-variable uniqueness successively in $r_1,r_2,c$, and then the finite-flat argument of that lemma, proves the assertion.

Finally, Remark~\ref{prop:LR-comparison} identifies the function in region (D) with that of \cite{LR3}. Since the comparison factor is the same for the depleted and ordinary inputs, \eqref{eq:global-improvement} holds in that normalization as well.
\end{proof}

\begin{rk}
At crystalline classical interpolation points, these identities replace the degree-eight Euler polynomial $\mathcal B_\diamond$ in the ordinary formula by $\mathcal B_\diamond^{\imp}$, with the substitutions above. The periods, archimedean factors, tame factors and common local normalization factors are unchanged. As identities of meromorphic functions, they also determine the behaviour at the boundary, where the improved function may have a pole.
\end{rk}

The Euler-factor identity also computes the value of the original function at a matching Yoshida point. We retain the fixed Siegel point and the logarithmic elliptic-weight coordinate of Section~\ref{subsec:poles-comparison}.

\begin{cor}\label{cor:ordinary-value}
At a matching ordinary Yoshida point $(P_0,Q_0)$ satisfying the hypotheses of Theorem~\ref{thm:criterion}, one has
\[
L_p^{(\diamond)}(P_0,Q_0)
=-(1-p^{-1})\,\partial_z\log\mathfrak a(Q_0)\,B_\diamond(\mathbf1),
\qquad\diamond=F,D.
\]
Here $B_\diamond(\mathbf1)$ is defined in \eqref{eq:constant-pairing} and $\partial_z\log\mathfrak a=\mathfrak a^{-1}\partial_z\mathfrak a$. For the tame data of Theorem~\ref{thm:criterion}, this value is nonzero if and only if $\partial_z\log\mathfrak a(Q_0)\ne0$.
\end{cor}
\begin{proof}
Write $L_p^{(\diamond),\imp}(P_0,Q)=R_\diamond/z+O(1)$. Since $e_\diamond(P_0,Q_0)=0$, Proposition~\ref{teo:interpolation} gives $L_p^{(\diamond)}(P_0,Q_0)=\partial_z e_\diamond(P_0,Q_0)R_\diamond$. The character $\chi_\Sigma$ is fixed, and $\gamma^\circ\mathfrak a(Q_0)=\beta^\circ/\mathfrak a^\vee(Q_0)=1$, so
\[
\partial_z e_F(Q_0)=-\partial_z\log\mathfrak a(Q_0),\qquad
\partial_z e_D(Q_0)=\partial_z\log\mathfrak a(Q_0).
\]
Combining these derivatives with the two signs in Proposition~\ref{prop:residue} proves the formula. The last assertion follows from the nonvanishing of $B_\diamond(\mathbf1)$ for the data of Theorem~\ref{thm:criterion}.
\end{proof}

Under the hypotheses of Hida's formula recalled in Section~\ref{sec:GL2xGL2}, one has $\partial_z\log\mathfrak a(Q_0)=-\tfrac12\mathcal L(\operatorname{ad}^0 V_g)$, where $g$ is the smaller constituent. Passing to its contragredient preserves the adjoint representation and, since the nebentype is fixed, the logarithmic derivative.

To compare coordinates, write $t=1+T=\varepsilon(\gamma)\gamma^{k-1}$ in Hida's convention, where $k$ is the modular weight, $\gamma$ generates $1+p\mathbb Z_p$ and $\varepsilon$ is fixed. Then $\partial_k=\log_p(\gamma)t\partial_t$, so the formula $-2\log_p(\gamma)\mathfrak a^{-1}t\partial_t\mathfrak a$ of \cite[Theorem 2.4]{Hida:AdjointConstancy} gives $-2\partial_k\log\mathfrak a$. Since $c$ is the modular weight of our family and $z=c-c_0$, this is also the derivative with respect to $z$. Its identification with Greenberg's invariant at a classical point uses the hypotheses of Theorem 2.1 of that reference.

Corollary~\ref{cor:ordinary-value}, with this identification of the derivative, proves Corollary~\ref{intro:ordinary-value}.

\section{Factorization formulas}\label{sec:factorization}

The pole criterion identifies the smaller constituent of a Yoshida lift. A factorization of the $p$-adic $L$-function would say more: it would express the residue as a product of a Rankin value for the larger constituent and an improved Rankin residue for the smaller one. We describe the expected identities, starting from the complex factorization and its weight shifts.

\subsection{The complex factorization}

Let $\Pi$ be a Yoshida lift of Siegel weight $(k_1,k_2)$, with $k_2\geq3$, attached to normalized newforms $f,g$ of weights
\[
K_f=k_1+k_2-2,\qquad K_g=k_1-k_2+2.
\]
Write $h$ for the elliptic form of weight $c\geq2$ corresponding to $\Sigma$. In unitary automorphic normalization the transfer of $\Pi$ is $\Sigma_f\boxplus\Sigma_g$. Hence the primitive transfer $L$-function satisfies
\begin{equation}\label{eq:yoshida-transfer-factorization}
L_{\mathrm{tr}}(\Pi\times\Sigma_h,s)
=L(\Sigma_f\times\Sigma_h,s)L(\Sigma_g\times\Sigma_h,s).
\end{equation}
The subscript $\mathrm{tr}$ indicates the primitive finite function defined by the transfer. Under the local compatibility of Theorem~\ref{intro:complex}, it agrees with the finite part of the integral definition used in Section~3. We may also delete the factors at a common set of ramified primes throughout \eqref{eq:yoshida-transfer-factorization}.

The two Rankin functions have different motivic weights, since $K_f-K_g=2(k_2-2)$. Under the change of variable
\[
s'=s+\frac{k_1+k_2+c-4}{2}
\]
gives the two Rankin arguments
\[
j_f=s+\frac{K_f+c-2}{2}=s',\qquad
j_g=s+\frac{K_g+c-2}{2}=s'-(k_2-2).
\]

Thus the motivically normalized transfer function satisfies \begin{equation}\label{eq:yoshida-motivic-factorization} L_{\mathrm{tr}}^{\mathrm{mot}}(\Pi\times h,s') =L(f\otimes h,s')L(g\otimes h,s'-(k_2-2)). \end{equation}
This shift determines the Rankin specialization corresponding to each of our three-variable domains.

Write $g^*$ for the contragredient form, with the ordinary refinement fixed as in Section~\ref{sec:GL2xGL2}. At $h=g^*$, the completed unitary factorization gives
\[
\operatorname{Res}_{s=1}\Lambda(\Pi\times\Sigma_g^\vee,s)
=\Lambda(\Sigma_f\times\Sigma_g^\vee,1)\,
\operatorname{Res}_{s=1}\Lambda(\Sigma_g\times\Sigma_g^\vee,s).
\]
The first factor is nonzero by \cite{Sha81}; its finite motivic value is $L(f\otimes g^*,k_1)$, with $K_g\leq k_1\leq K_f-1$. The analogous $p$-adic formula should involve this larger Rankin factor and the residue of the improved function at the diagonal.

\subsection{Expected $p$-adic factorizations}

Let $\hf,\hg$ be ordinary families whose Yoshida lifts form a family $\underline\Pi$, and let $\hh$ be an ordinary family through $h$. We choose compatible central characters and refinements, with the weight relations above. Thus the factorization is considered along a Yoshida family, where both elliptic constituents vary with $\Pi$.

Each Rankin function is normalized using its dominant form. With that form written first, as in Section~\ref{sec:GL2xGL2}, the order of the factors is
\[
\begin{array}{c|c|c}
\text{region}&\text{weight inequalities}&\text{ordered Rankin pairs}\\ \hline
(A)&c>K_f>K_g&(h,f),\ (h,g)\\
(D)&K_f>c>K_g&(f,h),\ (h,g)\\
(F)&K_f>K_g>c&(f,h),\ (g,h).
\end{array}
\]
We match the Euler factors by deleting a common set of tame primes when necessary. The factors $C_\diamond$ account for the comparison of Siegel and Rankin periods and the remaining local normalizations; local factors depending on the cyclotomic variable are kept explicit. The conjecture includes the asserted variation of $C_\diamond$.

\begin{conj}\label{conj:fact}
On a compatible Yoshida family as above, with matching local normalizations, there are nonzero meromorphic functions $C_A,C_D,C_F$ of the three weight variables, independent of the cyclotomic variable, such that
\begin{align*}
L_p^{(A)}(k_1,k_2,c,s')
&=C_A\,L_p(\hh,\hf)(c,K_f,s')
             L_p(\hh,\hg)(c,K_g,s'-(k_2-2)),\\
L_p^{(D)}(k_1,k_2,c,s')
&=C_D\,L_p(\hf,\hh)(K_f,c,s')
             L_p(\hh,\hg)(c,K_g,s'-(k_2-2)),\\
L_p^{(F)}(k_1,k_2,c,s')
&=C_F\,L_p(\hf,\hh)(K_f,c,s')
             L_p(\hg,\hh)(K_g,c,s'-(k_2-2)).
\end{align*}
The identities are taken on the domains of the corresponding Siegel functions and restrict to the three-variable domains used here. Where a four-variable extension is needed, its existence is included in this conjectural formulation; Remark~\ref{prop:LR-comparison} concerns the original three-variable domain. As before, the weights denote points of the chosen eigenfamilies, and the shifts are interpreted using weight characters, including their finite-order parts.
\end{conj}

The ordering is consistent with the critical ranges. In motivic notation, the inequalities of Section~\ref{subsec:degree-eight} become
\[
\begin{array}{c|c|c}
&j_f&j_g\\ \hline
(A)&K_f\leq j_f\leq c-1&K_g\leq j_g\leq c-1\\
(D)&c\leq j_f\leq K_f-1&K_g\leq j_g\leq c-1\\
(F)&c\leq j_f\leq K_f-1&c\leq j_g\leq K_g-1.
\end{array}
\]
In each row, the degree-eight critical range is the intersection of the two shifted Rankin critical ranges. The complex factorization therefore compares critical values on both sides.

Liu proves a factorization of this form for a measure valued in automorphic forms on a Yoshida family: \cite[Proposition 4.2.3]{liu2} expresses it as a product of two Rankin measures and an explicit factor attached to the family. The construction in \cite[\S4.2]{liu2} uses a common finite discrete-series place and a compatible square root of the product of central characters. Comparing it with our conventions also entails the twist of the second constituent and the local vectors and periods of our coherent pairings.

For our coherent functions, the remaining comparison concerns the periods and local normalizations in families. Once these are identified, \eqref{eq:yoshida-motivic-factorization} and the interpolation formulas would give the equality at classical critical points, and uniqueness on the common domain would give the family identity. On the three-variable domains studied here, the cyclotomic argument is already determined by the weights.

\subsection{Relation with the improved functions}

At the boundary considered in Section~\ref{subsec:poles-comparison}, the matching form is $h=g^*$. The Rankin factor involving $g$ should therefore account for the pole. Its argument must first be compared with the endpoint used to define the improved Hida--Rankin function in Section~\ref{sec:GL2xGL2}.

The Rankin factor involving $g$ is evaluated at $j_g=s'-(k_2-2)$. On the two domains considered in Section~\ref{sec:improved}, this gives
\[
\begin{array}{c|c|c}
&s'&j_g\\ \hline
(F)&k_1-1&K_g-1\\
(D)&c+k_2-3&c-1.
\end{array}
\]
In (F), $g$ is dominant and the argument is its upper critical endpoint $K_g-1$; in (D), $h$ is dominant and the endpoint is $c-1$. The improved Rankin function in Section~\ref{sec:GL2xGL2} is defined at the lower endpoint. The functional equation relates these two choices.

For forms of weights $u>v$, the functional equation replaces the motivic argument $j$ by
\[
u+v-1-j,
\]
and the forms by their contragredients. It sends $j=u-1$ to $j=v$, exactly the required lower endpoint. Transporting the improved function through this equation requires the corresponding functional-equation factors, ordinary refinements and removed Euler factors. We retain these in the notation below.

Fix now $\Pi$, and hence $f$ and $g$, and let $h$ vary. If these comparisons hold, the resulting improved factorization should take the form
\[
L_p^{(\diamond),\imp}(\Pi\times h)
=
C_\diamond(h)\,
R_\diamond^{\imp}(g,h)\,
R_\diamond^{\mathrm{maj}}(f,h),
\qquad \diamond\in\{F,D\}.
\]
Here $R_\diamond^{\imp}(g,h)$ is the improved Rankin factor transported to the upper endpoint, with the ordering and contragredients prescribed by the functional equation. We write $R_\diamond^{\mathrm{maj}}(f,h)$ for the Rankin factor attached to the larger constituent. The factor $C_\diamond(h)$ includes the period and local comparisons, with the tame multiplier $\nu_K$ from Definition~\ref{def:joint-functions}.

Suppose this identity holds near $h_0=g^*$ and that $C_\diamond$ and $R_\diamond^{\mathrm{maj}}$ are regular there. Taking residues in the logarithmic elliptic-weight coordinate $z$ of \eqref{eq:coordinates} gives
\[
\operatorname{Res}_{z=0}
L_p^{(\diamond),\imp}(\Pi\times h)
=
C_\diamond(h_0)\,R_\diamond^{\mathrm{maj}}(f,h_0)\,
\operatorname{Res}_{z=0}R_\diamond^{\imp}(g,h).
\]
Thus the simple pole of the improved Rankin factor gives a simple pole of the improved Siegel function whenever
\[
C_\diamond(h_0)\,R_\diamond^{\mathrm{maj}}(f,h_0)\ne0.
\]
Section~\ref{subsec:poles-comparison} already proves the pole criterion through the distinguished period. The factorization would additionally determine the residue in terms of the two elliptic constituents, under the stated regularity and nonvanishing conditions.

\subsection{The Saito--Kurokawa analogy}

The same question also arises for Saito--Kurokawa lifts, whose spin parameter has an elliptic constituent and two one-dimensional constituents. In the classical scalar-weight convention, a lift $F$ of weight $(k,k)$ is attached to an elliptic form $f$ of weight $2k-2$. At full level and for even $k$,
\[
L_{\mathrm{spin}}(F,s') =L(f,s')\zeta(s'-k+1)\zeta(s'-k+2);
\]
see \cite[Introduction]{Schmidt:SK:Functoriality}. Tensoring the corresponding unramified parameters with those of $h$ gives, outside the ramified places,
\[
L_{\mathrm{spin}\otimes\mathrm{std}}(F\times h,s')
=L(f\otimes h,s')L(h,s'-k+1)L(h,s'-k+2).
\]
A $p$-adic analogue would thus involve one Hida--Rankin function and two elliptic $p$-adic $L$-functions. The Saito--Kurokawa representation is nongeneric, so this question requires a cohomological construction beyond the globally generic setting of the present paper.

\bibliographystyle{amsalpha}
\bibliography{refs}

@Unpublished{boxerpilloni20,
	author = {Boxer, George and Pilloni, Vincent},
	note   = {preprint},
	title  = {Higher {C}oleman theory},
	year   = {2020},
	url = {https://perso.ens-lyon.fr/vincent.pilloni/HigherColeman.pdf},
}

@article {CLJ,
    AUTHOR = {Cauchi, Antonio and Lemma, Francesco and Rodrigues Jacinto,
              Joaqu\'in},
     TITLE = {Algebraic cycles and functorial lifts from {$G_2$} to {$\mathrm
              PGSp_6$}},
   JOURNAL = {Algebra Number Theory},
  FJOURNAL = {Algebra \& Number Theory},
    VOLUME = {19},
      YEAR = {2025},
    NUMBER = {3},
     PAGES = {551--616},
      ISSN = {1937-0652,1944-7833},
   MRCLASS = {11F46 (11F67 14G10 14G35)},
  MRNUMBER = {4879361},
       DOI = {10.2140/ant.2025.19.551},
       URL = {https://doi-org.ezbusc.usc.gal/10.2140/ant.2025.19.551},
}

@misc{GrahamRockwood,
    AUTHOR = {Graham, Andrew and Rockwood, Rob},
    TITLE = {Nearly higher {C}oleman theory and p-adic L-functions for $\mathrm{GSp}(4) \times \mathrm{GL}(2)$ and $\mathrm{GSp}(4) \times \mathrm{GL}(2) \times \mathrm{GL}(2)$},
    YEAR = {2024},
    NOTE = {preprint, available at \url{https://arxiv.org/pdf/2411.04559.pdf}},
}

@article {hida2,
    AUTHOR = {Hida, Haruzo},
     TITLE = {A {$p$}-adic measure attached to the zeta functions associated
              with two elliptic modular forms. {II}},
   JOURNAL = {Ann. Inst. Fourier (Grenoble)},
  FJOURNAL = {Universit\'e{} de Grenoble. Annales de l'Institut Fourier},
    VOLUME = {38},
      YEAR = {1988},
    NUMBER = {3},
     PAGES = {1--83},
      ISSN = {0373-0956,1777-5310},
   MRCLASS = {11S40 (11F11 11F33 11F67 11F85)},
  MRNUMBER = {976685},
MRREVIEWER = {Neal\ Koblitz},
       DOI = {10.5802/aif.1141},
       URL = {https://doi.org/10.5802/aif.1141},
}

@article {KLZ,
	AUTHOR = {Kings, Guido and Loeffler, David and Zerbes, Sarah Livia},
	TITLE = {Rankin-{E}isenstein classes and explicit reciprocity laws},
	JOURNAL = {Camb. J. Math.},
	FJOURNAL = {Cambridge Journal of Mathematics},
	VOLUME = {5},
	YEAR = {2017},
	NUMBER = {1},
	PAGES = {1--122},
	ISSN = {2168-0930},
	MRCLASS = {11F85 (11F67 11G40 14G35 14H52)},
	DOI = {10.4310/CJM.2017.v5.n1.a1},
	URL = {http://dx.doi.org/10.4310/CJM.2017.v5.n1.a1},
}

@article {lemma20,
    AUTHOR = {Lemma, Francesco},
     TITLE = {Algebraic cycles and residues of degree 8 {$L$}-functions of
              {$\mathrm{GSp}(4)\times \mathrm{GL}(2)$}},
   JOURNAL = {Int. Math. Res. Not. IMRN},
  FJOURNAL = {International Mathematics Research Notices. IMRN},
      YEAR = {2020},
    NUMBER = {12},
     PAGES = {3614--3648},
      ISSN = {1073-7928,1687-0247},
   MRCLASS = {11F67 (14G35)},
  MRNUMBER = {4120306},
MRREVIEWER = {Yongxiong\ Li},
       DOI = {10.1093/imrn/rny114},
       URL = {https://doi.org/10.1093/imrn/rny114},
}

@article {LSZ22,
    AUTHOR = {Loeffler, David and Skinner, Christopher and Zerbes, Sarah
              Livia},
     TITLE = {Euler systems for {${\mathrm GSp}(4)$}},
   JOURNAL = {J. Eur. Math. Soc. (JEMS)},
  FJOURNAL = {Journal of the European Mathematical Society (JEMS)},
    VOLUME = {24},
      YEAR = {2022},
    NUMBER = {2},
     PAGES = {669--733},
      ISSN = {1435-9855,1435-9863},
   MRCLASS = {11F46 (11F67 11F80 11R23)},
  MRNUMBER = {4382481},
MRREVIEWER = {Meng\ Fai\ Lim},
       DOI = {10.4171/jems/1124},
       URL = {https://doi.org/10.4171/jems/1124},
}

@Unpublished{liu1,
  author = {Liu, Zheng},
  note   = {preprint},
  title  = {{$p$}-adic {$L$}-functions for {$\GSp(4)\times \GL(2)$}},
  year   = {2023},
}

@Unpublished{liu2,
  author = {Liu, Zheng},
  note   = {preprint},
  title  = {{$p$}-adic {$L$}-functions for {$\GSp(4)\times \GL(2)$} II},
  year   = {2024},
  eprint = {2406.05847},
  url = {https://arxiv.org/abs/2406.05847},
}

@misc{LZ-zeta,
    AUTHOR = {Loeffler, David},
    TITLE = {On some zeta-integrals for unramified representations of {$\operatorname{GSp}(4)$}},
    YEAR = {2021},
    NOTE = {preprint, available at \url{https://arxiv.org/pdf/2106.16076.pdf}},
}

@Article{LPSZ1,
  author  = {Loeffler, David and Pilloni, Vincent and Skinner, Christopher and Zerbes, Sarah Livia},
  journal = {Duke Math. J.},
  title   = {Higher {H}ida theory and {$p$}-adic {$L$}-functions for {$\operatorname{GSp}(4)$}},
  year    = {2021},
  number  = {18},
  pages   = {4033--4121},
  volume  = {170},
  doi     = {10.1215/00127094-2021-0049},
}

@article {LR2,
    AUTHOR = {Loeffler, David and Rivero, \'Oscar},
     TITLE = {Algebraicity of {$L$}-values for {$\mathrm GSp_4 \times GL_2$} and
              {$\mathrm GSp_4 \times GL_2 \times GL_2$}},
   JOURNAL = {Q. J. Math.},
  FJOURNAL = {The Quarterly Journal of Mathematics},
    VOLUME = {75},
      YEAR = {2024},
    NUMBER = {2},
     PAGES = {391--412},
      ISSN = {0033-5606,1464-3847},
   MRCLASS = {11F66},
  MRNUMBER = {4765775},
       DOI = {10.1093/qmath/haae016},
       URL = {https://doi.org/10.1093/qmath/haae016},
}

@article{LR3,
  author = {Loeffler, David and Rivero, \'Oscar},
  title = {On {$p$}-adic {$L$}-functions for {$\GSp(4)\times \GL(2)$}},
  journal = {Pacific J. Math.},
  volume = {335},
  number = {2},
  pages = {373--400},
  year = {2025},
  doi = {10.2140/pjm.2025.335.373},
}

@misc{LZ-BK,
    AUTHOR = {Loeffler, David and Zerbes, Sarah},
    TITLE = {On the Bloch--Kato conjecture for {$\operatorname{GSp}(4)\times\operatorname{GL}(2)$}},
    YEAR = {2021},
    NOTE = {preprint, available at \url{https://arxiv.org/pdf/2106.14511.pdf}},
}

@Unpublished{LZvista,
	author = {Loeffler, David and Zerbes, Sarah Livia},
	note   = {preprint},
	title  = {$p$-adic {L}-functions and diagonal cycles for $\operatorname{GSp}(4)\times\operatorname{GL}(2)\times\operatorname{GL}(2)$},
	year   = {2020},
}

@Unpublished{LZ25-poles,
	author = {Loeffler, David and Zerbes, Sarah Livia},
	note   = {preprint},
	title  = {Poles of {$p$}-adic Asai {$L$}-functions and distinguished representations},
	year   = {2025},
    eprint = {2501.11472}
}

@Article{Pil20,
  author  = {Pilloni, Vincent},
  journal = {Duke Math. J.},
  title   = {Higher coherent cohomology and {$p$}-adic modular forms of singular weights},
  year    = {2020},
  number  = {9},
  pages   = {1647--1807},
  volume  = {169},
  doi     = {10.1215/00127094-2019-0075},
}

@article {RR21,
    AUTHOR = {Rivero, \'Oscar and Rotger, Victor},
     TITLE = {Derived {B}eilinson-{F}lach elements and the arithmetic of the
              adjoint of a modular form},
   JOURNAL = {J. Eur. Math. Soc. (JEMS)},
  FJOURNAL = {Journal of the European Mathematical Society (JEMS)},
    VOLUME = {23},
      YEAR = {2021},
    NUMBER = {7},
     PAGES = {2299--2335},
      ISSN = {1435-9855,1435-9863},
   MRCLASS = {11F67 (11F85 11G18 19F27)},
  MRNUMBER = {4269414},
MRREVIEWER = {Matteo\ Longo},
       DOI = {10.4171/jems/1054},
       URL = {https://doi.org/10.4171/jems/1054},
}

@article {Su24,
    AUTHOR = {Su, Jun},
     TITLE = {Coherent cohomology of {S}himura varieties and automorphic
              forms},
   JOURNAL = {Ann. Sci. \'Ec. Norm. Sup\'er. (4)},
  FJOURNAL = {Annales Scientifiques de l'\'Ecole Normale Sup\'erieure.
              Quatri\`eme S\'erie},
    VOLUME = {57},
      YEAR = {2024},
    NUMBER = {4},
     PAGES = {1039--1099},
      ISSN = {0012-9593,1873-2151},
   MRCLASS = {14G35 (11F70 11G18 17B56)},
  MRNUMBER = {4773300},
MRREVIEWER = {Rolf\ Berndt},
       DOI = {10.24033/asens.2586},
       URL = {https://doi-org.ezbusc.usc.gal/10.24033/asens.2586},
}

@article {TU99,
    AUTHOR = {Tilouine, J. and Urban, E.},
     TITLE = {Several-variable {$p$}-adic families of {S}iegel-{H}ilbert
              cusp eigensystems and their {G}alois representations},
   JOURNAL = {Ann. Sci. \'Ecole Norm. Sup. (4)},
  FJOURNAL = {Annales Scientifiques de l'\'Ecole Normale Sup\'erieure.
              Quatri\`eme S\'erie},
    VOLUME = {32},
      YEAR = {1999},
    NUMBER = {4},
     PAGES = {499--574},
      ISSN = {0012-9593},
   MRCLASS = {11F41 (11F75 11F80 11F85)},
  MRNUMBER = {1693583},
MRREVIEWER = {Alexey\ A.\ Panchishkin},
       DOI = {10.1016/S0012-9593(99)80021-4},
       URL = {https://doi.org/10.1016/S0012-9593(99)80021-4},
}

@book{Weissauer:Endoscopy,
 author = {Weissauer, Rainer},
 title = {Endoscopy for {{\(\mathrm{GSp}(4)\)}} and the cohomology of {Siegel} modular threefolds},
 fseries = {Lecture Notes in Mathematics},
 series = {Lect. Notes Math.},
 issn = {0075-8434},
 volume = {1968},
 isbn = {978-3-540-89305-9; 978-3-540-89306-6},
 year = {2009},
 publisher = {Berlin: Springer},
 language = {English},
 doi = {10.1007/978-3-540-89306-6},
 zbMATH = {5490504},
 Zbl = {1273.11089}
}

@article{Moriyama:ZetaIntegrals,
 author = {Moriyama, Tomonori},
 title = {{{\(L\)}}-functions for {{\(\text{GSp}(2)\times \text{GL}(2)\)}}: archimedean theory and applications},
 fjournal = {Canadian Journal of Mathematics},
 journal = {Can. J. Math.},
 issn = {0008-414X},
 volume = {61},
 number = {2},
 pages = {395--426},
 year = {2009},
 language = {English},
 doi = {10.4153/CJM-2009-021-x},
 zbMATH = {5549889},
 Zbl = {1225.11071}
}

@article{Loeffler:Nov,
 author = {Loeffler, David},
 title = {On local zeta-integrals for {{\(\mathrm{GSp}(4)\)}} and {{\(\mathrm{GSp}(4) \times \mathrm{GL}(2)\)}}},
 fjournal = {The New York Journal of Mathematics},
 journal = {New York J. Math.},
 issn = {1076-9803},
 volume = {30},
 pages = {1--23},
 year = {2024},
 language = {English},
 url = {nyjm.albany.edu/j/2024/30-1.html},
 zbMATH = {7807584},
 Zbl = {1535.22056}
}

@article {Shapiro:Soudry,
    AUTHOR = {Piatetski-Shapiro, Ilya I. and Soudry, David},
     TITLE = {{$L$}\ and {$\varepsilon $}\ factors for {${ GSp}(4)$}},
   JOURNAL = {J. Fac. Sci. Univ. Tokyo Sect. IA Math.},
  FJOURNAL = {Journal of the Faculty of Science. University of Tokyo.
              Section IA. Mathematics},
    VOLUME = {28},
      YEAR = {1981},
    NUMBER = {3},
     PAGES = {505--530 (1982)},
      ISSN = {0040-8980},
   MRCLASS = {22E50 (10D30 10D40)},
  MRNUMBER = {656034},
MRREVIEWER = {Guy\ Henniart},
}

@book{Knapp:Representations,
 author = {Knapp, Anthony W.},
 title = {Representation theory of semisimple groups. {An} overview based on examples. {With} a new preface by the author.},
 edition = {Paperback ed.},
 isbn = {0-691-09089-0},
 year = {2001},
 publisher = {Princeton, NJ: Princeton University Press},
 language = {English},
 zbMATH = {1685565},
 Zbl = {0993.22001}
}

@article{Soudry,
 author = {Soudry, David},
 title = {The {L} and {{\(\gamma\)}} factors for generic representations of {GSp}(4,k){{\(\times GL(2,k)\)}} over a local nonarchimedean field k},
 fjournal = {Duke Mathematical Journal},
 journal = {Duke Math. J.},
 issn = {0012-7094},
 volume = {51},
 pages = {355--394},
 year = {1984},
 language = {English},
 doi = {10.1215/S0012-7094-84-05117-2},
 zbMATH = {3887813},
 Zbl = {0557.12012}
}

@article{Vogan:Zuckerman,
  author  = {Vogan, David A. and Zuckerman, Gregg J.},
  title   = {Unitary representations with non-zero cohomology},
  journal = {Compositio Mathematica},
  volume  = {53},
  number  = {1},
  year    = {1984},
  pages   = {51--90},
  url     = {https://www.numdam.org/item/CM_1984__53_1_51_0/}
}

@article{Boxer:Pilloni:Siegel,
 author = {Boxer, George and Pilloni, Vincent},
 title = {Higher {Hida} theory for {Siegel} modular forms},
 fjournal = {Inventiones Mathematicae},
 journal = {Invent. Math.},
 issn = {0020-9910},
 volume = {244},
 number = {1},
 pages = {45--141},
 year = {2026},
 language = {English},
 doi = {10.1007/s00222-025-01393-2},
 zbMATH = {8188948}
}

@article{Harris:Partial,
 author = {Harris, Michael},
 title = {Automorphic forms of {{\({\bar \partial}\)}}-cohomology type as coherent cohomology classes},
 fjournal = {Journal of Differential Geometry},
 journal = {J. Differ. Geom.},
 issn = {0022-040X},
 volume = {32},
 number = {1},
 pages = {1--63},
 year = {1990},
 language = {English},
 doi = {10.4310/jdg/1214445036},
 zbMATH = {4168820},
 Zbl = {0711.14012}
}

@article{AsgariShahidi,
 author = {Asgari, Mahdi and Shahidi, Freydoon},
 title = {Generic transfer from {GSp}(4) to {GL}(4)},
 journal = {Compositio Mathematica},
 volume = {142},
 number = {3},
 year = {2006},
 pages = {541--550},
 doi = {10.1112/S0010437X06001904},
 eprint = {math/0412505}
}

@article{Schmidt:SK:Functoriality,
 author = {Schmidt, Ralf},
 title = {The {Saito--Kurokawa} lifting and functoriality},
 journal = {American Journal of Mathematics},
 volume = {127},
 number = {1},
 year = {2005},
 pages = {209--240},
 url = {https://www.math.uni-sb.de/ag/schulze/rschmidt/papers/skfunc.pdf}
}

@article{GT11,
 author = {Gan, Wee Teck and Takeda, Shuichiro},
 title = {The local {Langlands} conjecture for {{\(\text{GSp}(4)\)}}},
 fjournal = {Annals of Mathematics. Second Series},
 journal = {Ann. Math. (2)},
 issn = {0003-486X},
 volume = {173},
 number = {3},
 pages = {1841--1882},
 year = {2011},
 language = {English},
 doi = {10.4007/annals.2011.173.3.12},
 zbMATH = {5960694},
 Zbl = {1230.11063}
}

@article{JPSS83,
 author = {Jacquet, Herv{\'e} and Piatetski-Shapiro, Ilya I. and Shalika, Joseph},
 title = {Rankin-{Selberg} convolutions},
 fjournal = {American Journal of Mathematics},
 journal = {Amer. J. Math.},
 issn = {0002-9327},
 volume = {105},
 pages = {367--464},
 year = {1983},
 language = {English},
 doi = {10.2307/2374264},
 url = {semanticscholar.org/paper/8c35a3abc5d201e32e6720703f58964e1d572d0b},
 zbMATH = {3832249},
 Zbl = {0525.22018}
}

@article{Del74,
 author = {Deligne, Pierre},
 title = {La conjecture de {Weil}: {I}},
 journal = {Publications Math{\'e}matiques de l'IH{\'E}S},
 volume = {43}, year = {1974}, pages = {273--307},
 doi = {10.1007/BF02684373}
}

@article{Fox03,
 author = {Fox, Glenn J.},
 title = {A method of {Washington} applied to the derivation of a two-variable {$p$}-adic {$L$}-function},
 journal = {Pacific Journal of Mathematics},
 volume = {209}, number = {1}, year = {2003}, pages = {31--40}
}

@article{Horawa:Prasanna,
 author = {Horawa, Aleksander and Prasanna, Kartik},
 title = {Motivic action for {Siegel} modular forms},
 journal = {Journal of the European Mathematical Society},
 year = {2025}, note = {Published online},
 doi = {10.4171/JEMS/1714}, url = {https://arxiv.org/abs/2307.04115}
}

@article{Pil13,
 author = {Pilloni, Vincent},
 title = {Overconvergent modular forms},
 journal = {Annales de l'Institut Fourier}, volume = {63}, number = {1},
 year = {2013}, pages = {219--239}, doi = {10.5802/aif.2759}
}

@article{Sha81,
 author = {Shahidi, Freydoon}, title = {On certain {$L$}-functions},
 journal = {American Journal of Mathematics}, volume = {103}, number = {2},
 year = {1981}, pages = {297--355}, doi = {10.2307/2374219}
}

@book{Wal92,
 author = {Wallach, Nolan R.}, title = {Real Reductive Groups {II}},
 series = {Pure and Applied Mathematics}, volume = {132-II},
 publisher = {Academic Press}, year = {1992}
}

@article{BP:Modular,
 author = {Boxer, George and Pilloni, Vincent},
 title = {Higher {Hida} and {Coleman} theories on the modular curve},
 journal = {{\'E}pijournal de G{\'e}om{\'e}trie Alg{\'e}brique},
 volume = {6}, year = {2022}, doi = {10.46298/epiga.2022.6112},
 url = {https://arxiv.org/abs/2002.06845}
}

@misc{RC:EichlerShimura,
 author = {Rodr{\'i}guez Camargo, Juan Esteban},
 title = {{$p$}-adic {Eichler--Shimura} maps for the modular curve},
 year = {2023}, note = {arXiv:2102.13099v3},
 url = {https://arxiv.org/abs/2102.13099}
}

@article{Hida:Adjoint,
 author={Hida, Haruzo},title={{Greenberg}'s {$\mathcal L$}-invariants of adjoint square {Galois} representations},
 journal={Int. Math. Res. Not.},year={2004},number={59},pages={3177--3189},
 doi={10.1155/S1073792804133126}
}

@article{Hida:AdjointConstancy,
 author={Hida, Haruzo},title={Constancy of adjoint {$\mathcal L$}-invariant},
 journal={J. Number Theory},volume={131},year={2011},pages={1331--1346},
 doi={10.1016/j.jnt.2011.02.001}
}
\enlargethispage{4\baselineskip} % Keep the author addresses on the final page.

\end{document}